\documentclass[12pt,a4paper,final]{amsart}
\usepackage{version}
\usepackage[notcite,notref]{showkeys}
\usepackage{amsmath,amsthm,amssymb,latexsym,epsfig,graphicx,subfigure}
\numberwithin{equation}{section}

\newtheorem{thm}{Theorem}[section]
\newtheorem{lm}[thm]{Lemma}

\newtheorem{cor}[thm]{Corollary}
 
\newtheorem{pr}[thm]{Proposition}
\theoremstyle{definition}
\newtheorem{df}[thm]{Definition}
\theoremstyle{definition}

\newtheorem{ex}[thm]{Example}
\newtheorem{conj}[thm]{Conjecture}

\newcommand{\Rn}{\mathbb{R}^{n}}
\newcommand{\hn}{\mathbb{R}^{n-1}}
\newcommand{\R}{\mathbb{R}}

\newcommand{\Z}{\mathbb{Z}}

\renewcommand{\hn}{\mathbb{H}^n}
\newcommand{\C}{\mathbb{C}}

\renewcommand{\H}{\mathcal{H}}

\newcommand{\h}{\mathbb{H}}

\newcommand{\Lip}{\text{Lip}}

\newcommand {\grtrsim} {\ {\raise-.5ex\hbox{$\buildrel>\over\sim$}}\ }
\newcommand{\e}{\varepsilon}

\newcommand{\khii}{\text{\lower -.4ex\hbox{$\chi$}}}
\DeclareMathOperator{\spt}{spt}

\DeclareMathOperator{\tanm}{Tan} \DeclareMathOperator{\di}{div}
\DeclareMathOperator{\card}{card} \DeclareMathOperator{\sing}{Sing}
 \DeclareMathOperator{\reg}{Reg}
 \DeclareMathOperator{\apTan}{apTan} 
\DeclareMathOperator{\curl}{curl}\DeclareMathOperator{\biLip}{biLip}
\DeclareMathOperator{\Ric}{Ric}

\renewcommand{\a}{\alpha}

\newcommand{\restrict}{\begin{picture}(12,12)
                       \put(2,0){\line(1,0){8}}
                       \put(2,0){\line(0,1){8}}
                      \end{picture}}

\begin{document}
\title {Rectifiability; a survey}
\author{Pertti Mattila}


\begin{abstract} 
This is a survey on rectifiability. I discuss basic properties of rectifiable sets, measures, currents and varifolds and their role in complex and harmonic analysis, potential theory, calculus of variations, PDEs and some other topics.
\end{abstract}

\maketitle
{\small\tableofcontents}
\section{Introduction}
Rectifiable sets, measures, currents and varifolds are basic objects of geometric measure theory. In particular during the last four decades they have spread out to many areas of analysis and geometry. One of the goals of this survey is to show how rectifiability unifies surprisingly many different topics.  Starting from the beginning and basic theories, I shall briefly describe many of these appearances. The table of contents should give a pretty good idea of what will follow. Here I just mention some of the milestones.

The pre-beginning is the right generalization of length, area, etc. This was provided by Carath\'eodory in 1914. Lebesgue with his measure had given a generalization of volume in the Euclidean $n$-space and Carath\'eodory continued from this to define for integers $0<m<n$ an outer measure generalizing the $m$-dimensional area. This measure $\H^m$ is now called $m$-dimensional Hausdorff measure, since Hausdorff generalized it further in 1919 to non-integral values of $m$. Once equipped with this tool Besicovitch began in the 1920s to study properties of $\H^1$ measurable planar sets $E$ with $\H^1(E)<\infty$. In three papers
 he was able to reveal an amazing amount of structure. Such a set splits into a rectifiable and purely unrectifiable part, the former has properties similar to those of smooth surfaces and the latter completely opposite properties. Federer generalized most of Besicovitch's theory in 1947 to $m$-dimensional sets in $\Rn$. After the founding work of Besicovitch and Federer the ingenious ideas of Marstrand and Preiss have been the most influential for the basic theory of rectifiability.

In the 1950s De Giorgi described the structure of sets of finite perimeter in terms of rectifiable sets.
This gave rectifiability a permanent place in the calculus of variations.
First for sets of codimension one and then via Federer and Fleming's theory of normal and integral currents in 1960 for all dimensions. These are generalized surfaces and another class of them of fundamental importance over the years, rectifiable varifolds, was introduced by Almgren in the 1960s and developed by Allard in the 1970s. In the 1990s Simon proved some fundamental results on the rectifiability of the singularities of minimal currents and harmonic maps.

Rectifiability has played a big role in complex and harmonic analysis. In the 1950s Vitushkin anticipated this for the geometric description of removable sets of bounded complex analytic functions, fully confirmed much later in 1998 by David. In the 1990s, continuing from 
Jones's analyst's traveling salesman theorem, David and Semmes established the theory of uniform rectifiability and its connections to singular integrals and other topics of harmonic analysis.

Rectifiability has found a prominent place also outside Euclidean spaces. After some important work in the 1990s by Ambrosio, by Preiss and Tiser and by  Kirchheim, Ambrosio and Kirchheim developed the theory of rectifiable sets and currents in metric spaces in 2000. In 2001 Franchi, Serapioni and Serra Cassano introduced the right notions of rectifiability in Heisenberg groups, which has lead to an extensive theory in general Carnot groups.

This survey covers many topics, but all of them briefly and only scratching the surface. I am trying to give the  reader a flavour of each of them without detailed proofs, often with ideas of the proofs, and often just presenting the results. There probably are many topics and articles which I have not, but should have, mentioned. And certainly this text is biased, I have concentrated more on topics I know best, including my own work.

{\bf Acknowledgements.} I would like to thank David Bate, Vasilis Chousionis, Damian Dabrowski, Antti K\"aenm\"aki, Ulrich Menne, David Preiss and Xavier Tolsa for useful comments.

\section{Preliminaries}\label{Preli}
\subsection{Notation}\label{Notation}
We denote by $\mathcal L^n$ the Lebesgue measure in the Euclidean $n$-space $\Rn$. In a metric space $X$,\ $d(A)$ stands for the diameter of $A$,\ $d(A,B)$ the minimal distance between the sets $A$ and $B$, $d(x,A)$ the distance from a point $x$ to a set $A$. The closed ball with centre $x\in X$ and radius $r>0$ is denoted by $B(x,r)$ and the open ball by $U(x,r)$. In $\Rn$ we sometimes denote $B^n(x,r)$. The unit sphere in $\Rn$ is $S^{n-1}$.  
The Grassmannian manifold of linear $m$-dimensional subspaces of $\Rn$ is $G(n,m)$. It is equipped with an orthogonally invariant Borel probability measure $\gamma_{n,m}$. For $V\in G(n,m)$ we denote by $P_V$ the orthogonal projection onto $V$.

For $A\subset X$  we denote by $\mathcal M(A)$ the set of non-zero finite Borel measures $\mu$ on $X$ with support $\spt\mu\subset A$. We shall denote by $f_{\#}\mu$ the push-forward of a measure $\mu$  under a map  $f: f_{\#}\mu(A)= \mu(f^{-1}(A))$. The restriction of $\mu$ to a set $A$ is defined by $\mu\restrict A(B)=\mu(A\cap B)$. The notation $\ll$ stands for absolute continuity.

The characteristic function of a set $A$ is $\chi_A$. By the notation $M\lesssim N$ we mean that $M\leq CN$ for some constant $C$. The dependence of $C$ should be clear from the context. The notation $M \sim N$ means that $M\lesssim N$ and $N\lesssim M$. By  $c$ and $C$ we mean positive constants with obvious dependence on the related parameters. 

\subsection{Hausdorff measures}
For $m\geq 0$ the $m$-dimensional Hausdorff measure $\H^m=\H^m_d$ in a metric space $(X,d)$ is defined by

$$\H^m(A)=\lim_{\delta\to 0}\inf\{\sum_{i=1}^{\infty}\a(m)2^{-m}d(E_i)^m: A\subset \bigcup_{i=1}^{\infty}E_i, d(E_i)<\delta\}.$$
Then $\H^0$ is the counting measure. Usually $m$ will be a positive integer and then $\a(m)=\mathcal L^m(B^m(0,1))$, from which it follows that $\H^m=\mathcal L^m$ in $\R^m$. For non-integral values of $m$ the choice of $\a(m)$ does not really matter. We denote by $\dim$ the Hausdorff dimension. The \emph{spherical Hausdorff measure} $\mathcal S^m$ is defined in the same way but using only balls as covering sets.

The lower and upper $m$-densities of $A\subset X$ are defined by

$$\Theta^m_{\ast}(A,x)=\liminf_{r\to 0}\a(m)^{-1}r^{-m}\mathcal H^m(A\cap B(x,r)),$$ 
$$\Theta^{\ast m}(A,x)=\limsup_{r\to 0}\a(m)^{-1}r^{-m}\mathcal H^m(A\cap B(x,r)).$$ 
The density $\Theta^m(A,x)$ is defined as their common value if they are equal.

We have 
\begin{thm}\label{dens}If $A$ is $\H^m$ measurable and $\mathcal H^m(A)<\infty$, then
$$2^{-m}\leq\Theta^{\ast m}(A,x)\leq 1\ \text{for}\ \mathcal H^m\ \text{almost all}\ x\in A,$$
$$\Theta^{\ast m}(A,x)=0\ \text{for}\ \mathcal H^m\ \text{almost all}\ x\in \Rn\setminus A.$$
\end{thm} 

When $m\leq 1$ the constant $2^{-m}$ is sharp, for $m>1$ the best constants are not known.

We also have

\begin{thm}\label{densdiam}If $A\subset X$ is $\H^m$ measurable and $\mathcal H^m(A)<\infty$, then
$$\lim_{\delta\to 0}\sup\{d(B)^{-m}\H^m(A\cap B): x\in B, d(B)<\delta\}= 1\ \text{for}\ \mathcal H^m\ \text{almost all}\ x\in A.$$
\end{thm}

For general measures we have

\begin{thm}\label{densmeas}Let $\mu\in\mathcal M(X), A\subset X,$  and $0<\lambda<\infty$.
\begin{itemize}
\item[(1)] If $\Theta^{\ast m}(A,x)\leq\lambda$ for $x\in A$, then $\mu(A)\leq 2^m\lambda \H^m(A).$
\item[(2)] If $\Theta^{\ast m}(A,x)\geq\lambda$ for $x\in A$, then $\mu(A)\geq \lambda \H^m(A).$
\end{itemize}
\end{thm}

For the above results see \cite[2.10.17-19]{Fed69}, \cite[Section 2.2]{Fal85}, or \cite[Chapter 6]{Mat95}.

We say that a closed set $E$ is AD-$m$-regular (AD for Ahlfors and David) if there is a positive number $C$ such that
$$r^m/C \leq \H^m(E\cap B(x,r)) \leq Cr^m\ \text{for}\ x\in E, 0<r<d(E).$$
A measure $\mu$ is said to be AD-$m$-regular if
$$r^m/C \leq \mu(B(x,r)) \leq Cr^m\ \text{for}\ x\in \spt\mu, 0<r<d(\spt\mu),$$
which means that $\spt\mu$ is an AD-$m$-regular set.

\subsection{Lipschitz maps}\label{Lipmaps} 
Since Lipschitz maps are at the hearth of rectifiability, I state some basic well-known facts about them. We say that a map $f:X\to Y$ between metric spaces $X$ and $Y$ is \emph{Lipschitz} if there is a positive number $L$ such that 
$$d(f(x),f(y))\leq Ld(x,y)\ \text{for}\ x,y\in X.$$
The smallest such $L$ is the Lipschitz constant of $f$, it is denoted by $\Lip(f)$.

Euclidean valued Lipschitz maps $f:A\to\R^k, A\subset X,$ can be extended: there is a Lipschitz map $g:X\to\R^k$ such that $g|A=f$, see \cite[2.10.43-44]{Fed69} of \cite[Chapter 7]{Mat95}.

Any Lipschitz map $g:\R^m\to\R^k$ is almost everywhere differentiable by Rademacher's theorem, see \cite[3.1.6]{Fed69} or \cite[7.3]{Mat95}.

There is the Lusin type property: if $f:A\to\R^k, A\subset\R^m$, is Lipschitz, then for every $\e>0$ there is a $C^1$ map $g:\R^m\to\R^k$ such that
\begin{equation}\label{Lusin}
\mathcal L^m(\{x\in A: g(x)\neq f(x)\})<\e,\end{equation}
see \cite[3.1.16]{Fed69}. 

\section{Rectifiable curves}\label{curves}

Let us first have a quick look at  rectifiable curves concentrating on some facts which are relevant for the more general rectifiable sets.

By a curve in $\Rn$ we mean a continuous image of a line segment, $C=f([a,b])$. The curve is rectifiable if you can rectify it, that is, take hold of the end-points and pull it straight. This is the same as to say that the curve has finite length. A standard definition of the length is that it is the total variation of $f$:

$$V_f(a,b)= \sup\{\sum_{j=1}^k|f(x_j)-f(x_{j-1})|: a=x_0<x_1<\dots <x_k=b\}.$$

However, this depends on $f$ since $f$ could travel through some parts of $C$ several times. For us the length of $C$ will be the 1-dimensional Hausdorff measure $\mathcal H^1(C)$ of $C$. It agrees with $V_f(a,b)$ if $f$ is injective, or more generally if the set points of $C$ which are covered more than once has zero $\mathcal H^1$ measure. 

If $V_f(a,b)<\infty$, then $f$ is a function of bounded variation. Such functions have many well-known nice properties, but for us it is important to know that we can do better: if $f$ is the arc-length parametrization of $C$, then $f$ is Lipschitz. This is essential in particular in the case of higher dimensional rectifiable sets.

Let now $C=f([a,b])$ be a rectifiable curve in $\Rn$ with a Lipschitz parametrization $f$. Here are some of its basic easily verifiable properties:

\begin{equation}\label{length}
Area\ formula: \int_C N(f,y)\,d\mathcal H^1 y = \int_a^b \sqrt{f'_1(x)^2+\dots + f'_n(x)^2}\,dx,
\end{equation}
where $N(f,y)$ is the number of points $x\in [a,b]$ with $f(x)=y$. So in particular
$$\mathcal H^1(C)=\int_a^b \sqrt{f'_1(x)^2+\dots +f'_n(x)^2}\,dx,$$
if $f$ is injective. The key for the proof is Rademacher's theorem (or Lebesgue's in the one-dimensional case); Lipschitz mappings are almost everywhere differentiable. Using also that $\sqrt{f'_1(x)^2+\dots+f'_n(x)^2}$ tells us how the derivative of $f$ at $x$ changes length a rather elementary proof can be given. As expected, a higher dimensional version also is valid, and will be presented later. Hence the name area formula.

With the help of the area formula and again Rademacher's theorem the following two properties are not too hard to verify:

\begin{equation}\label{tangents}
Tangents: C\ \text{has a tangent at}\ \mathcal H^1\ \text{almost every point}\ x\in C.
\end{equation}

\begin{equation}\label{density}
Density: \lim_{r\to 0}\frac{\mathcal H^1(C\cap B(x,r))}{2r} = 1\ \text{for}\ \mathcal H^1\ \text{almost all}\ x\in C.
\end{equation}

In the plane the length can be computed by counting the intersection points with lines:

\begin{equation}\label{Croft}
Crofton\ formula: 2\mathcal H^1(C)= \int\text{card} (C\cap L)\,dL.
\end{equation} 

Here the measure $dL$ on lines can be obtained by parametrizing the lines as $\{te+a:t\in\R\}, e\in S^{1}, a\in e^{\perp}$, and integrating over $e$ and $a$. 

This formula is trivially checked when $C$ is a line segment. The general case can be done using Rademacher's theorem and approximation by polygonal curves. Crofton proved this in 1868, which marked the beginning of integral geometry. Unless you want to start at 1777 with Count Buffon and his needle.

In the beginning I said that a curve is rectifiable if it has finite length. But if we take Hausdorff measure as length can we get from its finiteness the Lipschitz parametrization which was used above? Yes, we can, even in general metric spaces:

\begin{thm}\label{cont}
If $X$ is a metric space and $C\subset X$ is a compact connected set with $\mathcal H^1(C)<\infty$, then there is a Lipschitz mapping $f:[0,1]\to X$ with $f([0,1])=C$.
\end{thm} 

For a rather easy proof in $\Rn$, see \cite[Theorem I.1.8]{DS93}, and in the Hilbert space, \cite[Lemma 3.7]{Sch07}. Here are some ideas. For each $\delta>0$ choose a maximal $\delta$ separated subset $A_{\delta}$ of $C$. Connect with line segments all those pairs of points of $A_{\delta}$ which have distance at most $2\delta$ and let $C_{\delta}$ be the union of these segments. Playing with some graphs show that $C_{\delta}$ is a continuum which can be parametrized by a Lipschitz map $f_{\delta}:[0,1]\to \Rn$ with 
$\Lip(f_{\delta})\lesssim \H^1(C)$. Finally use the Arzela-Ascoli theorem to get $f$ as the limit of some sequence $(f_{\delta_j})$.

Theorem \ref{cont} was proved by Eilenberg and Harrold in \cite{EH43}. It is one of the reasons why the rectifiability theory often is much easier for one-dimensional sets. Another reason is compactness and lower semicontinuity:

\begin{thm}\label{cont1}
If $C_k\subset B^n(0,1), k=1,2,\dots,$ are continua, then there is a subsequence $C_{k_j}$ converging in the Hausdorff distance to a continuum $C$ with $\mathcal H^1(C) \leq \liminf_{j\to\infty}\mathcal H^1(C_{k_j})$.
\end{thm} 

Again the proof is rather easy, see \cite[Theorems 3.16 and 3.18]{Fal85}.

\section{One-dimensional rectifiable sets}\label{Onedim}
Almost all of the theory of one-dimensional rectifiable sets in the plane was developed by Besicovitch in the three papers \cite{Bes28}, \cite{Bes38} and \cite{Bes39}. Some generalizations are due to Morse and Randolph \cite{MR44} and Moore \cite{Moo50}. For most of the proofs of the results in this chapter, see \cite{Fal85} and \cite{Deg81}. Several ideas extend to higher dimensions and in the next chapter we shall mainly discuss the new ideas needed. We consider here only subsets of the plane although essentially everything holds for one-dimensional rectifiable subsets of $\Rn$. 

\subsection{Definitions and tangents}\label{rect1def}

One-dimensional rectifiable sets generalize rectifiable curves and they have many properties which make this class flexible in measure theoretic sense and useful in many applications. In particular,\\
- every rectifiable curve is a rectifiable set,\\
- countable unions of rectifiable sets are rectifiable,\\
- subsets of  rectifiable sets are rectifiable,\\
- sets of zero $\H^1$ measure are rectifiable,\\
- if for every $\e>0$  there is a rectifiable set $F\subset E$ with $\H^1(E\setminus F)<\e$, then $E$ is rectifiable.

These properties offer an immediate definition: $E$ is 1-rectifiable if $\H^1$ almost of it can be covered with rectifiable curves. I state it a bit differently:

\begin{df}\label{1-rect}
A set $E\subset\R^2$ is \emph{$1$-rectifiable} if there are Lipschitz maps $f_i:\R\to \R^2, i=1,2,\dots,$ such that
$$\H^1(E\setminus\bigcup_{i=1}^{\infty}f_i(\R))=0.$$
\end{df}

Next we shall see that rectifiable sets have, in an appropriate sense, the properties of rectifiable curves which we presented in the previous chapter. This is actually rather easy, but more essentially only rectifiable sets have these properties. That is, each of the properties allows a converse statement.

We shall begin with tangents. Consider the following example. Let $q_i, 1=1,2,\dots$, be all the points of the plane with rational coordinates and set $E=\bigcup_{i=1}^{\infty}S_i$ with $S_i=\{x\in\R^2:|x-q_i|=2^{-i}\}$. Then $E$ is 1-rectifiable and $\H^1(E)<\infty$. But it is dense in $\R^2$. So how can it have any tangents? Obviously it cannot in the ordinary sense and we have to give a measure theoretic definition of tangent. For $a\in\R^2, s>0, L\in G(2,1)$, define the cone (angular sector)

$$X(a,L,s)=\{x\in\R^2: d(x,L+a)<s|x-a|\}.$$

\begin{df}\label{apprtan}
A line $L\in G(2,1)$ is an \emph{approximate tangent line} of a set $E\subset\R^2$ at a point $a\in\R^2$ if $\Theta^{\ast 1}(E,a)>0$ and for every $s>0$,
$$\lim_{r\to 0}r^{-1}\H^1(\{x\in E\cap B(a,r)\setminus X(a,L,s)\})=0.$$
\end{df}

It is convenient to define the approximate tangents as lines through the origin. Then the geometric  approximate tangent at $a$ is the translate by $a$.

Now we immediately have that the above union $E$ of the circles $S_i$ has an approximate tangent at almost all of its points: each $S_i$ has an ordinary tangent at all of its points and at almost all points it is an approximate tangent of $E$ by the density theorem \ref{dens}. Since rectifiable curves have tangents almost everywhere, essentially the same argument shows that every 1-rectifiable set with finite $\H^1$ measure has an approximate tangent at almost all of its points.

The converse is not very difficult either. I give the idea assuming that $E$ has ordinary tangents almost everywhere. Then we can write $E$ as a countable union of a set of measure zero and sets $F$ for which there exist $s>0$ and $L\in G(2,1)$ such that $F\setminus X(a,L,s)=\emptyset$ for $a\in F$. The set $F$ is such that the tangents in its points are close to a fixed line $L$. This implies that the restriction of the projection $P_L$ to $F$ is one-to-one with a Lipschitz inverse. Hence $F=(P_L|F)^{-1}(P_L(F))$ is rectifiable.

The case of approximate tangents causes technical difficulties, but the main idea is the same. So we have

\begin{thm}\label{tanthm1}
If $E$ is $\H^1$ measurable and $\mathcal H^1(E)<\infty$, then $E$ is 1-rectifiable if and only it has an approximate tangent line at almost all of its points.
\end{thm}

\begin{df}\label{1-unrect}
A set $E\subset\R^2$ is \emph{purely $1$-unrectifiable} if $\H^1(E\cap F)=0$ for every 1-rectifiable set $F\subset \R^2$ (or, equivalently, $\H^1(E\cap C)=0$ for every rectifiable curve $C$). 
\end{df}

The proof of the following proposition is an easy exercise:

\begin{pr}
If $E$ is $\H^1$ measurable and $\mathcal H^1(E)<\infty$, then $E=R\cup P$ where $R$ is 1-rectifiable and $P$ is purely 1-unrectifiable.
\end{pr}

\begin{ex}\label{Cantor}
A standard example of a  purely 1-unrectifiable compact subset of the plane is the self-similar Cantor set $C$ obtained by choosing the squares of side-length 1/4 in the corners of the unit square and continuing indefinitely. That is, $C=C_1\times C_1$ where $C_1$ is a symmetric Cantor set on the line of Hausdorff dimension $1/2$. Then $0<\H^1(C)<\infty$. It is easy to check that $C$ has no approximate tangents.
\end{ex}

Obviously all properties of rectifiable sets correspond to properties of purely unrectifiable sets, and vice versa. In particular, $E$ is purely unrectifiable if and only it has no approximate tangent at almost all of its points. But looking at the above argument we can say more. We did not need large sectors $B(a,r)\setminus X(a,L,s)$ in the complement of $F$, arbitrarily narrow ones were enough. This leads to, although not quite trivially, see \cite{Fal85}, Theorem 3.29, or \cite{Mat95}, Corollary 15.16, 

\begin{thm}\label{unrectsect}
If $E$ is $\H^1$ measurable and purely 1-unrectifiable  and $\mathcal H^1(E)<\infty$, then  for every $L\in G(2,1)$ and every $s>0$,
$$\limsup_{r\to 0}r^{-1}\H^1(E\cap B(a,r)\cap X(a,L,s))\gtrsim s$$
for $\H^1$ almost all $a\in E$.
\end{thm}

\subsection{Densities}\label{dens1}

Next we look at density properties of rectifiable sets. We again get easily from the area formula for Lipschitz maps that rectifiable sets have density 1 almost everywhere. Also the converse holds:

\begin{thm}\label{densthm1}
If $E$ is $\H^1$ measurable and $\mathcal H^1(E)<\infty$, then $E$ is 1-rectifiable if and only $\Theta^1(E,x)=1$ for $\H^1$ almost all $x\in E$.
\end{thm}

This is where Besicovitch used connectivity essentially. He did it in terms of circle pair regions 
$$U(x,y)=U(x,|x-y|)\cap U(y,|x-y|),\ x,y\in \R^2.$$ To see how they appear notice that if $C\subset\R^2$ is compact and $U(x,y)\cap C\not=\emptyset$ for $x,y\in C$, then $C$ is connected. If not, $C$ is the union of two disjoint compact sets $C_1$ and $C_2$. Let $x_i\in C_i$ with $|x_1-x_2|=d(C_1,C_2)$. Then $C\cap U(x_1,x_2)=\emptyset$, which contradicts our assumption.

Besicovitch took the density 1 property as the starting point and defined a set to be regular if the conclusion of Theorem \ref{densthm1} holds. The term rectifiable was introduced by Federer in \cite{Fed47a}, although his terminology too is somewhat different from mine, see \cite[3.2.14]{Fed69}.

With some more work one can then prove, see \cite[ Lemma 3.22]{Fal85},

\begin{lm}\label{circpair}
Let $E$ be $\H^1$ measurable and $\mathcal H^1(E)<\infty$. Suppose that $\alpha>0$ and $F\subset E$ is compact for which $\H^1(F)>0$ and $\H^1(E\cap U(x,y))>\alpha |x-y|$ for $x,y\in F$. Then there is a continuum $C$ such that $0<\H^1(C\cap E)\leq \H^1(C)<\infty$. In particular, $E$ is not purely 1-unrectifiable.
\end{lm}


\begin{thm}\label{3/4-dens}
If $E$ is $\H^1$ measurable, $\mathcal H^1(E)<\infty$ and $\Theta^1_{\ast}(E,x) > 3/4$ for $\H^1$ almost all $x\in E$, then $E$ is 1-rectifiable.
\end{thm}

This immediately implies Theorem \ref{densthm1}.

Once we have Lemma \ref{circpair} available the proof of Theorem \ref{3/4-dens} is simple: If it is false, we can find, using also Theorem \ref{densdiam}, a compact subset $F$ of $E$ with $\H^1(F)>0$ and positive numbers $\alpha$ and $r_0$ such that 
$$\H^1(E\cap U(x,r)) > (\tfrac{3}{4}+\a)2r\ \text{for}\ x\in F, 0<r<r_0,$$
and  
$$\H^1(E\cap B) \leq (1+\a)d(B)\ \text{whenever}\ F\cap B\not=\emptyset\ \text{and}\ d(B)<r_0.$$
If $d(F)<r_0/3$, then for $x,y\in F$ with $r=|x-y|$,
\begin{align*}
&\H^1(E\cap U(x,y))\\
&=\H^1(E\cap U(x,r))+\H^1(E\cap U(y,r))-\H^1(E\cap (U(x,r)\cup U(y,r))\\
&>2\cdot (\tfrac{3}{4}+\a)2r - (1+\a)3r = \a r = \a|x-y|.
\end{align*}
This shows by Lemma \ref{circpair} that $E$ is not purely unrectifiable, which is enough.

The bound $3/4$ is not sharp; Preiss and Tiser \cite{PT92} improved it by small $\e>0$. They did it in general metric spaces, we shall return to this in Chapter \ref{metricspaces}. The best bound is not known but it is expected to be $1/2$. This is the famous Besicovitch 1/2-problem:

\begin{conj}\label{1/2}
If $E$ is $\H^1$ measurable, $\mathcal H^1(E)<\infty$ and $E$ is purely 1-unrectifiable, then $\Theta^1_{\ast}(E,x)\leq 1/2$ for $\H^1$ almost all $x\in E$.
\end{conj}

Partial results, with some extra conditions on $E$, were obtained by Farag in \cite{Far00a} and \cite{Far02}.

Besicovitch used circle pairs to prove a stronger form of Theorem \ref{densthm1}: rectifiability follows from the existence of density. So we have

\begin{thm}\label{densthm1b}
Let $E$ be $\H^1$ measurable with $\mathcal H^1(E)<\infty$. Then the following are equivalent 
\begin{itemize}
\item[(1)] $E$ is 1-rectifiable.
\item[(2)] $\Theta^1(E,x)=1$ for $\H^1$ almost all $x\in E$.
\item[(3)] $\Theta^1(E,x)$ exists for $\H^1$ almost all $x\in E$.
\end{itemize}
\end{thm}

\subsection{Projections}

When discussing tangents we saw that rectifiable sets have projections of positive length on tangent lines. A closer look at that argument shows that if $E$ is 1-rectifiable with $\H^1(E)>0$, then there is at most one line $L\in G(2,1)$ such that $\H^1(P_L(E))=0$. We can characterize rectifiable sets via projections, but it is more natural to state this for purely unrectifiable sets:

\begin{thm}\label{proj1}
If $E$ is $\H^1$ measurable and $\mathcal H^1(E)<\infty$, then $E$ is purely 1-unrectifiable if and only $\H^1(P_L(E))=0$ for almost all $L\in G(2,1)$.
\end{thm}

What is left to prove is that any purely 1-unrectifiable $\H^1$ measurable set $E$ with $\mathcal H^1(E)<\infty$ projects to measure zero in almost all directions. This is the famous Besicovitch projection theorem. Here are some ideas for the proof.

The essence of the proof is to study how $E$ behaves close to lines. Fix first $L\in G(2,1)$. Let $\delta>0$ and split $E$ into three parts. The first $E_1(\delta,L)$ consists of points $a\in E$ such that only a little  measure is approaching $a$ in the direction $L$: 
$$\limsup_{s\to 0}\sup_{0<r<\delta}(sr)^{-1}\H^1(E\cap B(a,r) \cap X(a,L,s)) = 0.$$
The second set $E_2(\delta,L)$ consists of points $a\in E$ such that a lot of measure is approaching $a$ in the direction $L$: 
$$\limsup_{s\to 0}\sup_{0<r<\delta}(sr)^{-1}\H^1(E\cap B(a,r) \cap X(a,L,s)) = \infty.$$
It might be that nothing else is needed. Maybe almost all pairs $(a,L)$ must satisfy one of these properties. But that is not known. What saves us is that it suffices to add another property. The third set $E_3(L)$ consists of points $a\in E$ such that $E\cap (L+a)\cap B(a,\delta) \not=\emptyset$ for all $\delta>0$. 

The last property is easy to deal with because it is rather simple to show that for any set $F\subset\R^2$ with $\H^1(F)<\infty$,\ $\card(F\cap (L+x)) < \infty$ for almost all $x\in L^{\perp}$. Hence $\H^1(P_{L^{\perp}}(E_3(L))=0$.

The second case also is fairly easy. It is only slightly more complicated than the case where the cones $X(a,L,s)$ are replaced by cylinders $C(a,L,r,s)=\{x\in B(a,r): d(x,L+a)<sr\}$. If many of these have big measure as compared to $sr$, then applying Vitali's covering theorem on $L^{\perp}$, we find that 
$\H^1(P_{L^{\perp}}(E_3(L)))=0$.

We have already dealt with the first set: by Theorem \ref{unrectsect} $\H^1(E_1(\delta,L))=0$. Then after some technicalities we are left to show that for almost all $(a,L)$ if $a$ is the only point of $E$ on $E\cap B(a,\delta)\cap (L+a)$, then 
$$\limsup_{s\to 0}\sup_{0<r<\delta}(sr)^{-1}\H^1(E\cap B(a,r) \cap X(a,L,s))$$
is either 0 or $\infty$. 

For this we can use a general lemma due to Mickle and Rad\'o \cite{MR58}:

\begin{lm}\label{mickle}

If $\psi$ is an outer measure on $\Rn$ and $A$ is Lebesgue measurable with $\psi(A)=0$, then for $\mathcal L^n$ almost all $x\in A$,\ $\limsup_{r\to 0} r^{-n}\psi(B(x,r))$ is either 0 or $\infty$.
\end{lm}

The proof is a simple covering argument combined with Lebesgue density theorem. This lemma is applied on $G(2,1)$, identified with $S^1$, and I leave it to the reader to figure out how $\psi$ is chosen, or see \cite[Chapter 18]{Mat95}. Notice that we do not, and we cannot, assume that Borel sets are $\psi$ measurable.

The proof sketched above is the original Besicovitch's proof from \cite{Bes39}, but modified and polished by several people. No essentially different proof is known. Tao \cite{Tao09} proved a quantitative multiscale version of Besicovitch's projection theorem, but he did this quantizing Besicovitch's proof. Davey and Taylor \cite{DT21b} proved a non-linear version of Tao's result.

As mentioned above the third alternative $E_3(L)$ may not be needed at all, the following might be true:

\begin{conj}\label{unrect-rad}
If $E$ is $\H^1$ measurable,  purely 1-unrectifiable and $\mathcal H^1(E)<\infty$, then for $\H^1$ almost all $a\in E$ almost all lines through $a$ meet $E$ only at $a$.
\end{conj}

Although this is not known Cs\"ornyei and Preiss \cite{CP07} disproved a related conjecture which asked whether it is true that for every $\H^1$ measurable $E$ with  $\mathcal H^1(E)<\infty$ almost all lines through almost all points of $E$ meet $E$ in a finite set. Their set contains non-trivial rectifiable and purely unrectifiable parts.

The analogue of Crofton's formula \eqref{Croft} extends from rectifiable curves to 1-rectifiable sets with routine arguments.

\subsection{Analyst's traveling salesman problem}\label{tsp}

This topic is discussed in the books \cite{BP17}, \cite{Paj02} and \cite{Tol14}.

An analogue of the classical traveling salesman problem for general compact planar sets is: when can all points of a compact set $F\subset\R^2$ be traversed via a rectifiable curve and how long must such a curve be? Jones proved in \cite{Jon90} the following theorem. Let $F\subset\R^2$ be compact. For any square $Q$ set
$$\beta_F(Q)=\inf_{L\ \text{a line}}\sup\{d(x,L)/d(Q): x\in F\cap Q\}.$$
Let $\mathcal D$ be the family of all dyadic squares in the plane. 

\begin{thm}\label{jones}
A compact set $F\subset\R^2$ is contained in a rectifiable curve if and only if the Jones square function
$$\beta(F):=\sum_{Q\in \mathcal D}\beta_F(3Q)^2d(Q)< \infty.$$
The length of the shortest such curve is comparable to $\beta(F) + d(F)$.
\end{thm}

If $\beta(F)<\infty$ the key to the construction of the curve is Pythagoras theorem. Suppose that a line segment $I$ is a good approximation of $F$ inside a rectangle of side lengths $d$ and $\beta d$. To get good approximations in the next smaller scale, we might need to replace I by two line segment $I_1$ and $I_2$ with one common endpoint and the other endpoints common with those of $I$. Then the increase of length is roughly $d(\sqrt{1+\beta^2}-1)\sim \beta^2d$. The complete construction is not easy, but it is easier than the proof of the converse statement.

Theorem \ref{jones} was extended to $\Rn$ by Okikiolu \cite{Oki92}. Schul generalized it in \cite{Sch07} to infinite-dimensional Hilbert spaces. This means that he showed that the constants are independent of $n$. A different proof of Theorem \ref{jones} and a nice treatment of this topic can be found in the book of Bishop and Peres \cite{BP17}. They use arguments involving Crofton's formula.

Jones's theorem can easily be stated in integral form which leads to the following characterization of 1-rectifiable sets. Since sets of measure zero can affect greatly the $\beta$ numbers, it is better to state it for purely unrectifiable sets. Let 

\begin{equation}\label{beta}
\beta_E(x,r) = \inf_L\sup\{d(x,L)/r: x\in E\cap B(x,r)\}.
\end{equation}

\begin{thm}\label{betathm}
If $E\subset\R^2$ is $\H^1$ measurable and $\mathcal H^1(E)<\infty$, then $E$ is purely 1-unrectifiable if and only
$$\int_{\R^2}\int_0^{\infty}\beta_F(x,r)^2r^{-2}\,dr\,dx = \infty$$
for every Borel set $F\subset E$ with $\mathcal H^1(F)>0$.
\end{thm}

For a pointwise characterization in higher dimensions, see Theorem \ref{Azzam}.


I shall briefly mention another related way to characterize rectifiability in terms of \emph{Menger curvature} $c(x,y,z)$, which also measures approximation by lines. For three points $x,y,z$ in the plane it is defined as the reciprocal of the radius of the circle passing through these points. Thus $c(x,y,z)=0$ if and only if the points are collinear.  More quantitatively by elementary geometry,
\begin{equation}\label{menger}
c(x,y,z)=\frac{2d(x,L_{y,z})}{|x-y||x-z|}=\frac{4A(x,y,z)}{|x-y||x-z||y-z|},
\end{equation}
where $L_{y,z}$ is the line through $y$ and $z$ and $A(x,y,z)$ is the area of the triangle with vertices $x,y$ and $z$.

\begin{thm}\label{davleg}
If $E\subset\R^2$ is $\H^1$ measurable, $\mathcal H^1(E)<\infty$ and 
\begin{equation}\label{davlegeq}
\int_E\int_E \int_E c(x,y,z)^2\,d\H^1x\,d\H^1y\,d\H^1z < \infty,\end{equation}
then $E$ is 1-rectifiable
\end{thm}

This was first proved by David (unpublished) and then by L\'eger in \cite{Leg99}. 
In \cite{Tol05} Tolsa gave a different proof based on the analytic capacity.

The main problem for proving Theorem \ref{davleg} is that we do not have AD-regularity, or even positive lower density. With such assumptions the proof would be much easier. To get a vague idea, suppose that $E$ is AD-1-regular $x_i\in E, i=1,2,3,  |x_i-x_j|=2r$, when $i\neq j$, and $\beta_E(x_i,r)\sim 1$. Then with some $c>0$ and\\ $A_i=E\cap B(x_i,r)\setminus B(x_i,cr)$,
$$\int_{A_1}\int_{A_2} \int_{A_3} c(x,y,z)^2\,d\H^1x\,d\H^1y\,d\H^1z \sim r.$$
Clearly this cannot happen too often if \eqref{davlegeq} holds.

Lerman and Whitehouse \cite{LW11} and  Meurer \cite{Meu18} proved higher dimensional versions.

We shall return to this concept in several later chapters.

Here is another related square function: Let $\Gamma$ be a Jordan curve in the plane and let $\Omega^+$ and $\Omega^-$ be its interior and exterior domains. For $x\in\Gamma, r>0,$ let $I^+(x,r)$ and $I^-(x,r)$ be the longest arcs of $\partial B(x,r)$ contained in $\Omega^+$ and $\Omega^-$, respectively. Set
$$\e(x,r)= \max\{|\pi r-\H^1(I^+(x,r))|,|\pi r-\H^1(I^-(x,r))|\}/r$$
and define 
$$\mathcal E(x)^2 = \int_0^1\e(x,r)^2/r\,dr.$$
Jaye, Tolsa and Villa \cite{JTV21} proved the following theorem (and more) solving an old $\e^2$ conjecture of Carleson:

\begin{thm}
Let $\Gamma$ be a Jordan curve. Then for $\H^1$ almost all $x\in\Gamma$,\  $\Gamma$ has a tangent at $x$ if and only if $\mathcal E(x)<\infty.$
\end{thm}

That $\mathcal E(x)<\infty$ at almost every tangential point $x$ of $\Gamma$ is classical, so these authors proved the converse. The proof uses methods developed in connection of Theorems \ref{jones} and \ref{davleg}, but a lot of new ideas and techniques are also required.

\section{Higher dimensional rectifiable sets}\label{Higherdim}

Federer generalized most of Besicovitch's theory to higher dimensions in \cite{Fed47a}. Most of the proofs, or sketches with further references, for the results of this section can be found in \cite{Fed69}, \cite{LY02}  and \cite{Mat95}.

\subsection{Definitions and area and coarea formulas}

We now define

\begin{df}\label{m-rect}
A set $E\subset\R^n$ is \emph{$m$-rectifiable} if there are Lipschitz maps $f_i:\R^m\to \R^n, i=1,2,\dots,$ such that
$$\H^m(E\setminus\bigcup_{i=1}^{\infty}f_i(\R^m))=0.$$
A set $E\subset\R^n$ is \emph{purely $m$-unrectifiable} if $\H^m(E\cap F)=0$ for every $m$-rectifiable set $F\subset \R^n$.
\end{df}

Usually $m$ and $n$ will be integers with $0<m<n$, but sometimes $m$ can be 0, then 0-rectifiable means countable. We shall often also consider rectifiability of measures:

\begin{df}\label{m-rectmeas}
A measure $\mu$ on $\R^n$ is \emph{$m$-rectifiable} if there are Lipschitz maps $f_i:\R^m\to\Rn$ such that 
$$\mu(\Rn\setminus\bigcup_{i=1}^{\infty}f_i(\R^m))=0.$$

$\mu$ is \emph{purely $m$-unrectifiable} if $\mu(f(\R^m))=0$ for every Lipschitz map $f:\R^m\to\Rn$.
\end{df}
Often the condition $\mu\ll\H^m$ is added, but in may places later it is better to have the definition without it. In some texts it is only  required that $\mu(\Rn\setminus E)=0$ for some $m$-rectifiable set $E$. 

So the $m$-rectifiability of $E$  means that $\H^m\restrict E$ is $m$-rectifiable. We can uniquely decompose any $\mu\in\mathcal M(\R^n)$ as $\mu=\mu_r+\mu_u$, where $\mu_r$ is $m$-rectifiable and $\mu_u$ is purely $m$-unrectifiable. 

It often is useful to use other sets in place of Lipschitz images. The following alternatives give equivalent definitions:

\begin{itemize}
\item[(1)] Lipschitz images of arbitrary or compact subsets of $\R^m$.
\item[(2)] $C^1$ images of $\R^m$, or of arbitrary or compact subsets of $\R^m$.
\item[(3)] Lipschitz (or $C^1$) graphs over subsets of $m$-planes.
\item[(4)] $m$-dimensional $C^1$ submanifolds of $\Rn$.
\item[(5)] Level sets of regular $C^1$ mappings $f:\Rn\to\R^{n-m}$, that is, sets $\{x\in A: f(x)=y\}$ where $f:\Rn\to\R^{n-m}$ is $C^1$ with $\nabla f(x)\not=0$ for $x\in A$.
\end{itemize} 

The proofs are routine verifications applying classical theorems of analysis;  Whitney's extension theorem and implicit function theorem. And in particular Rademacher's theorem according to which a Lipschitz map $f$ is differentiable at almost every point $x$. That is, there is a linear map $Df(x)$ such that
\begin{equation}\label{diff}
f(x+y) - f(x) = Df(x)(y-x) + |x-y|\e(|x-y|),\ \lim_{h\to 0}\e(h)=0.
\end{equation}
From this we can often go from properties of linear maps to Lipschitz maps $f:\R^m\to\Rn, m\leq n$. For example, there is a Jacobian $Jf(x)$ defined in terms of the partial derivatives of $f$ such that $\H^m(Df(x)(A))=Jf(x)\mathcal L^m(A)$ for Lebesgue measurable sets $A\subset\R^m$.  This leads (but not trivially, see \cite[Section 3.2]{Fed69}  or \cite[Section 3.3]{EG92}) to
\begin{equation}\label{area}
Area\ formula:\ \int \card A\cap f^{-1}\{y\}\,d\mathcal H^m y = \int_AJf(x)\,d\mathcal L^mx.
\end{equation}
The proof consists of splitting $A=\cup_{i=1}^{\infty}A_i\cup B$. In each $A_i$ the $Df(x)$ is injective and close to a linear map $L_i$. In $B$,\ $Jf(x)=0$ and both sides of \eqref{area} are $0$ for $A=B$.

There is also the Fubini-type coarea formula for Lipschitz maps $f:\R^n\to\R^m, m\leq n$:
\begin{equation}\label{coarea}
Corea\ formula:\ \int \H^{n-m}(A\cap f^{-1}\{y\})\,d\mathcal L^m y = \int_AJf(x)\,d\mathcal L^nx.
\end{equation}
Here the Jacobian $J(x)$ of $f$ at $x$ again is a kind of determinant of $Df(x)$ determined by the property that the coarea formula is valid for linear maps. See \cite[Section 3.2]{Fed69} or \cite[Section 3.4]{EG92} for the linear algebraic definition and the proofs.

With some analysis this gives the rectifiability of level sets, see \cite[Theorem 3.2.15]{Fed69} or \cite[Remark 12.8]{Sim83}:

\begin{thm}\label{rectslice}
If $f:\R^n\to\R^m, m\leq n$, is Lipschitz, then  $f^{-1}\{y\}$ is $(n-m)$-rectifiable for $\mathcal L^m$ almost all $y\in\R^m$.
\end{thm}

Here 0-rectifiable means finite or countable.


\subsection{Tangent planes}\label{subtanplanes}

For the tangential properties we again define the cones

$$X(a,V,s)=\{x\in\R^n: d(x,V+a)<s|x-a|\},\ a\in\Rn, V\in G(n,m), s>0,$$
and the approximate tangent planes

\begin{df}\label{m-apprtan}
A plane $V\in G(n,m)$ is an \emph{approximate tangent plane} of a set $E\subset\R^n$ at a point $a\in\R^n$ if $\Theta^{\ast m}(E,a)>0$ and for every $s>0$,
$$\lim_{r\to 0}r^{-m}\H^m(E\cap B(a,r)\setminus X(a,V,s))=0.$$
We then denote $V=\apTan(E,a)$. 
\end{df} 

Observe that as Hausdorff measure also approximate tangent plane is a metric concept. We shall use exactly the same definition for other metrics.

The characterization of rectifiability follows by similar arguments as in the one-dimensional case:

\begin{thm}\label{tanthm}
If $E$ is $\H^m$ measurable and $\mathcal H^m(E)<\infty$, then $E$ is $m$-rectifiable if and only it has an approximate tangent plane at $\H^m$ almost all of its points.
\end{thm}

So $E$ is  $m$-rectifiable if for almost all $a\in E$ there is an $m$-plane which approximates $E$ well in all small balls $B(a,r)$. But what if we only have such an approximation in the weaker sense that the approximating plane is allowed to depend on the scale $r$? This is not enough as easy examples show, not even if $E$ would have positive lower density. One such example can be constructed as a subset of the modified von Koch snowflake curve where the angles go to zero but not too fast, see \cite[Section 20]{DS91}. But if we add the assumption  that $E$ has positive lower density at almost all of its points and the approximation is bilateral; not only the points of $E$ in $B(a,r)$ are close to a plane $W$ but also the points of $W\cap B(a,r)$ are close to $E$, then this weaker approximation implies rectifiability. This was proved by Marstrand in \cite{Mar61} for two-dimensional sets in $\R^3$ and generalized in \cite{Mat75} using Marstrand's fundamental ideas. I state this result below in terms of tangent measures.

\subsection{Tangent measures}\label{tanmeas}

Tangent measures were introduced by Preiss in \cite{Pre87} to solve the density characterization of rectifiability, which we shall discuss below. They have turned out to be useful in many other occasions, too.

Define

$$T_{a,r}(x)=(x-a)/r,\ x,a\in\Rn, r>0.$$

So $T_{a,r}$ blows up the ball $B(a,r)$ to the unit ball. Now we also blow up measures.

\begin{df}\label{tanmeasdf}
Let $\mu$ be a Radon measure on $\Rn$. A non-zero Radon measure $\nu$ is called a \emph{tangent measure} of $\mu$ at $a\in\Rn$ if there are sequences $(c_i)$ and $(r_i)$ of positive numbers such that $r_i\to 0$ and $c_iT_{a,r_i\#}\mu\to\nu$ weakly. We denote the set of tangent measures of $\mu$ at $a$ by $\tanm(\mu,a)$.
\end{df}

Tangent measures tell us how the measure looks locally.

Notice that this definition requires rather little structure. It is enough to have a locally compact metric group $(G,d)$ in place of $\Rn$ with dilations $\delta_r, r>0$, which are group homomorphisms such that $\delta_1$ is identity, $\delta_{rs}=\delta_r\circ\delta_s$ and $d(\delta_r(x),\delta_r(y))=rd(x,y)$.  Then we can use exactly the same definition for tangent measures. The following result was proved in \cite{Mat05} in this setting:

\begin{thm}\label{M2}
Let $\mu$ be a Radon measure on $G$. Then the following are equivalent:
\begin{itemize}
\item[(1)] For $\mu$ almost all $a\in G$\ $\mu$ has a unique (up to multiplication by a constant) tangent measure at $a$.
\item[(2)] For $\mu$ almost all $a\in G$\ there is a closed subgroup $H_a$ of $G$ which is invariant under the $\delta_r$ such that $\tanm(\mu,a)=\{c\lambda_a: 0<c<\infty\}$ where $\lambda_a$ is a left Haar measure of $H_a$.
\end{itemize}
\end{thm}

The point here with respect to rectifiability is that if something is rectifiable, it should look the same at all small scales around typical points. This theorem then tells us how it should look. We shall come back to this in Section \ref{Parabolic} and Chapter \ref{heisenberg}.

The following is essentially a restatement of Theorem \ref{tanthm}:

\begin{thm}\label{tanmthm}
If $E$ is $\H^m$ measurable and $\mathcal \H^m(E)<\infty$, then $E$ is $m$-rectifiable if and only for $\H^m$ almost all $a\in E$ there is $V_a\in G(n,m)$ such that every measure in $\tanm(\H^m\restrict E,a)$ is of the form $c\H^m\restrict V_a$ for some $0<c<\infty$.
\end{thm}

Let us call measures of the form $c\H^m\restrict V$ for some $V\in G(n,m), 0<c<\infty$, \emph{$m$-flat}. The following is the bilateral approximation criterion mentioned above, I state it for general measures:

\begin{thm}\label{tanmthm1}
Let $\mu\in\mathcal M(\Rn)$ be such that $0<\Theta_{\ast}^m(\mu,x)\leq\Theta^{\ast m}(\mu,x)<\infty$ for $\mu$ almost all $x\in \Rn$.  If for $\mu$ almost all $a\in \Rn$  every measure in $\tanm(\mu,a)$ is \emph{$m$-flat}, then $\mu$ is $m$-rectifiable.
\end{thm}

The proof is a bit tricky, but here is an idea in the plane for $\mu=\H^1\restrict E$. If $E$ were purely unrectifiable, it would project into a set of measure zero in almost all directions. In fact, we would not have to use the projection theorem since our assumptions combined with pure unrectifiability imply rather easily that all projections have measure zero, see \cite[Lemma 16.1]{Mat95}. So given $L\in G(2,1)$ with good approximation at some scale, it suffices to show that $\H^1(P_L(E))>0$. Suppose not and suppose $F\subset E$ is compact. Then we can find $a\in F$ such that all of $F$ lies in a half-plane with $a$ on its boundary, which is  orthogonal to $L$. If we had a good bilateral approximation at $a$ with some line, this line ought to be almost orthogonal to $L$. Of course, $a$ now is very special, but choosing $F$ suitably and using $\H^1(P_L(F))=0$, it is possible to find many such points. This idea combined with several others can be used to derive a contradiction with the fact that the upper density is bounded. Preiss gave a different proof in \cite{Pre87}.


In Theorem \ref{tanmthm1} the assumption on positive lower density cannot be dropped; Preiss \cite[5.9(2)]{Pre87} constructed an example of a purely 1-unrectifiable Borel set $A\subset\R^2$ with finite $\H^1$ measure for which all tangent measures are 1-flat at $\H^1$ almost all points.

Some more recent interesting results on tangent measures were proven by Kenig, Preiss and Toro in \cite{KPT09}. 

\subsection{Densities}\label{Densities}

Now we give density criteria for rectifiability. The proofs in the higher dimensional case are quite different from Besicovitch's arguments. We cannot use connectivity when $m>1$.  For instance  an $m$-dimensional analogue $F_m$ of Example \ref{Cantor} is contained in a continuum $C$ with $\H^m(C)<\infty$. This is easily seen by approximating $F_m$ with unions of small cubes and connecting them with line segments. A more interesting case is explained in \cite[4.2.25]{Fed69}.

\begin{thm}\label{densthm1a}
Let $E\subset\Rn$ be $\H^m$ measurable with $\mathcal H^m(E)<\infty$. Then the following are equivalent 
\begin{itemize}
\item[(1)] $E$ is $m$-rectifiable.
\item[(2)] $\Theta^m(E,x)=1$ for $\H^m$ almost all $x\in E$.
\item[(3)] $\Theta^m(E,x)$ exists for $\H^m$ almost all $x\in E$.
\end{itemize}
\end{thm}

That (1) implies (2), and hence (3), can be proven with the help of Rademacher's theorem, area formula and Lebesgue density theorem in $\R^m$. Of course, (2) implies (1) is a special case of (3) implies (1), but since the proof of the first is easier, although not easy, and it is based on different ideas, I say something about it, too. That (2) implies (1) was proved by Marstrand in \cite{Mar61} for $m=2, n=3$, and generalized in \cite{Mat75} relying heavily on Marstrand's ideas. The proof of (3) implies (1) is due to Preiss \cite{Pre87}. For this purpose he introduced the tangent measures. In \cite{Del08} De Lellis gives a very nice exposition of this proof. 

Chlebik \cite{Chl} generalized the part (2) implies (1) showing that there is $c(m)<1$ such that rectifiability already follows from $\Theta^m_{\ast}(E,x)>c(m)$ for $\H^m$ almost $x\in E$. Notice that $c(m)$ is independent of $n$. In fact, his proof also works in infinite-dimensional Hilbert spaces, but the implication $(3)\Rightarrow (1)$ in Theorem \ref{densthm1a} is false in infinite-dimensional Hilbert spaces, see Section \ref{dens2}. Preiss proved that the existence of density in (3) can be relaxed to $\Theta^m_{\ast}(E,x)>c(n,m)\Theta^{\ast m}(E,x)$ for $\H^m$ almost $x\in E$ where $c(n,m)<1$. Very little is known of these constants except when $m=1$, recall Besicovitch's $1/2$-problem from the previous chapter.

Preiss's result $(3)\Rightarrow (1)$ in Theorem \ref{densthm1a} can be stated more generally, but it actually is easily seen to be equivalent: 

\begin{thm}\label{preiss}
If $\mu\in\mathcal M(\Rn)$ and the positive and finite limit $\lim_{r\to 0}r^{-m}\mu(B(x,r))$ exists for $\mu$ almost all $x\in \Rn$, then $\mu$ is $m$-rectifiable.
\end{thm}
The equivalence to Theorem \ref{densthm1a} follows from the fact that $\mu$ and the restriction of $\H^m$ to $\{x:0<\lim_{r\to 0}r^{-m}\mu(B(x,r))<\infty\}$ are mutually absolutely continuous.

How do we get to rectifiability from density 1? The first step is Marstrand's reflection lemma.

\begin{lm}\label{refl}
Let $E\subset\Rn$ be $\H^m$ measurable with $\mathcal H^m(E)<\infty$ and $\Theta^m(E,x)=1$ for $\H^m$ almost all $x\in E$. If $\delta>0$, then there are $r_0>0$ and a compact subset $F$ of $E$ such that $\H^m(E\setminus F)<\delta$ and $d(2a-b,E)<\delta|a-b|$ whenever $a,b\in F$ and $|a-b|<r_0$.
\end{lm}

So for any such pair $a,b$ the symmetric point of $b$ with respect to $a$ is close to $E$. For one-dimensional sets this begins to look like rectifiability, but one can proceed from it also when $m>1$, although with many complications.

I explain the idea behind Lemma \ref{refl} when $m=1$. For the complete proof, see \cite{Mat75} or \cite{LY02}. Suppose for simplicity that we have found $F$ such that $\H^1(E\cap B(x,r)) = 2r$ when $x\in F$ and $r<r_0$. Let $a,b\in F$ with $r:=|a-b|<r_0$ and let $\e>0$ be much smaller than $\delta$. Let $B_1=B(a,(1-\e)r), B_2=B(b,\e r)$ and $B_3=B(2a-b,\delta r)$. Then $d(B_1\cup B_2 \setminus B_3)\leq 2(1-\e)r$, so we may (almost) assume by Theorem \ref{densdiam} that $\H^1(E\cap(B_1\cup B_2 \setminus B_3))\leq 2(1-\e)r$. If $E\cap B_3$ were empty, we would have $E\cap (B_1\cup B_2)\subset E\cap (B_1\cup B_2 \setminus B_3)$, whence
$$2r = 2(1-\e)r+2\e r = \H^1(E\cap B_1) + \H^1(E\cap B_2) \leq \H^1(E\cap(B_1\cup B_2 \setminus B_3))\leq 2(1-\e)r.$$
Lemma would follow from this contradiction.

Here are some basics behind Preiss's theorem (3) implies (1), that is, Theorem \ref{preiss}. 
The power of tangent measures is that they turn limiting conditions to uniform equations or inequalities. In this case, for $\mu$ almost all $a\in\Rn$ every $\nu\in\tanm(\mu,a)$ is \emph{$m$-uniform}, that is,
$$\nu(B(x,r))=cr^m\ \text{for all}\ x\in\spt\nu, r>0.$$
If we could show that all $m$-uniform measures are $m$-flat, we would be done by Theorem \ref{tanmthm1}. This is true for $m=1,2$, but not easy to show, in particular for $m=2$. However, it is false for $m>2$. Preiss observed that $\H^3$ restricted to the cone $\{x\in\R^4: x_4^2=x_1^2+x_2^2+x_3^2\}$ is 3-uniform. To overcome this problem when $m>2$ Preiss showed that if an $m$-uniform measure is not flat, then in a certain precise sense it is far from flat. Moreover, $\tanm(\mu,a)$ is connected. Hence it is enough to show that at almost all points $\mu$ has some flat tangent measures. This is an easier part of the argument and was already essentially proved by Marstrand in \cite{Mar64}.  

Suppose then that $\nu$ is an $m$-uniform measure such that $0\in\spt\nu$ and 
\begin{equation}\label{uni1}
\nu(B(x,r))=\a(m)r^m\ \text{for}\ x\in\spt\nu, r>0.
\end{equation}
Then we have the identities
\begin{equation}\label{uni12}
\int_{B(x,r)}(r^2-|x-y|^2)^2\,d\nu y = \int_{B(0,r)}(r^2-|y|^2)^2\,d\nu y,\ x\in\spt\nu, r>0,
\end{equation}
which are used to study $b_r$ and $Q_r$ defined by
\begin{equation}\label{uni2}
b_r\cdot v = \int_{B(0,r)}(r^2-|y|^2)(v\cdot y)\,d\nu y\big/\int_{B(0,r)}(r^2-|y|^2)\,d\nu y,\ v\in\Rn,
\end{equation}
\begin{equation}\label{uni3}
Q_r(v) = \int_{B(0,r)}(v\cdot y)^2\,d\nu y\big/\int_{B(0,r)}(r^2-|y|^2)\,d\nu y,\ v\in\Rn.
\end{equation}
It is shown that they have convergent subsequences $b_{r_i}\to b$ and $Q_{r_i}\to Q$ for which
$$\spt\nu\subset K:=\{x\in \Rn: Q(x) - |x|^2 +2b\cdot x = 0\}.$$
The proof for $m=1,2$ can be completed from this, but for $m>2$ much more is needed. In particular, Preiss performed a detailed study of higher order moments $\int (v\cdot y)^ke^{-s|y|^2}\, d\nu y, s>0,$ and their Taylor expansions.

Due to Theorem \ref{densdiam} in the language of tangent measures the assumption of density 1  corresponds to assuming in addition to \eqref{uni1} that $\nu(B(x,r))\leq \a(m)r^m$ for $x\in\Rn$ and $r>0$. Then one can show that $b=0$ and $\spt\nu =K$, which is an $m$-plane. This gives a different proof for (2) implies (1) in Theorem \ref{densthm1a}.

The structure of $m$-uniform measures is a very interesting and for a large part open problem in itself. Christensen \cite{Chr70} introduced the more general class of \emph{uniformly distributed} measures; $\nu(B(x,r))=\nu(B(y,r))$ for $x,y\in\spt\nu, r>0$. He showed that if such a measure has compact support, then it is contained in a sphere. Kirchheim and Preiss \cite{KP02} proved that the support of any uniformly distributed measure in $\Rn$ is an analytic variety. We already noted that $\nu=\H^3\restrict\{x\in\R^4: x_4^2=x_1^2+x_2^2+x_3^2\}$ is $3$-uniform. Kowalski and Preiss proved in \cite{KP87} that in addition to flat measures this is the only (up to translations and rotations) example in $\R^4$, and in $\Rn, n>4,$ all $(n-1)$-uniform measures are either flat or of the form  $\nu=\H^{n-1}\restrict\{x\in\R^n: x_4^2=x_1^2+x_2^2+x_3^2\}$. For a long time no non-flat $m$-uniform measures were known when $2<m<n-1$. Recently Nimer \cite{Nim18} produced many interesting examples. See also \cite{Tol15a}, \cite{Nim17} and \cite{Nim19} for other results.

I shall present a couple of results related to Theorem \ref{preiss}. In addition to that result Preiss's paper contains a lot of information about  general measures, including other characterizations of rectifiability. For example, see \cite[Theorem 4.11]{Pre87}:

\begin{thm}\label{preiss1}
Let $\mu\in\mathcal M(\Rn)$. Then 
\begin{equation}\label{preisseq}
\lim_{r\to 0}\mu(B(x,2r))/\mu(B(x,r))\ \text{exists for}\ \mu\ \text{almost all}\ x\in\Rn \end{equation}
if and only if all tangent measures of $\mu$ at $x$ are flat for $\mu$ almost all $x\in\Rn$.
\end{thm}

Here the flat measures can be of any dimension, but under the density assumptions of the next corollary, they are $m$-flat and we can use Theorem \ref{tanmthm1} to get an extension of Theorem \ref{preiss}:

\begin{cor}\label{preiss2}
Let $\mu\in\mathcal M(\Rn)$. If $0<\Theta_{\ast}^m(\mu,x)\leq \Theta^{\ast m}(\mu,x)<\infty$ for $\mu$ almost all $x\in\Rn$, then $\mu$ is $m$-rectifiable if and only if \eqref{preisseq} holds.
\end{cor}

Tolsa and Toro proved in \cite{TT15} a related result:

\begin{thm}\label{TT}
Let $\mu\in\mathcal M(\Rn)$ be such that $0<\Theta_{\ast}^m(\mu,x)\leq\Theta^{\ast m}(\mu,x)<\infty$ for $\mu$ almost all $x\in \Rn$. Then the following are equivalent:
\begin{itemize}
\item[(1)]$\mu$ is $m$-rectifiable.
\item[(2)]$\int_0^1\left|\frac{\mu(B(x,r))}{r^m}-\frac{\mu(B(x,2r))}{(2r)^m}\right|^2r^{-1}dr < \infty$ for $\mu$ almost all $x\in\Rn.$
\item[(3)]$\lim_{r\to 0}\left(\frac{\mu(B(x,r))}{r^m}-\frac{\mu(B(x,2r))}{(2r)^m}\right)=0$ for $\mu$ almost all $x\in\Rn.$
\end{itemize}
\end{thm}

The main contribution here is the implication $(1)\Rightarrow (2)$. Its proof uses Calder\'on-Zygmund techniques. That (3) implies (1) follows from Corollary \ref{preiss2}, because (3) clearly implies \eqref{preisseq} under the density assumptions.

The uniform rectifiability version of the equivalence of (1) and (2) was proved earlier by Chousionis, Garnett, Le and Tolsa in \cite{CGLT16}, see Theorem  \ref{CGLT}. In \cite{Tol17} Tolsa showed that when $m=1$ the equivalence of (1) and (2) holds assuming only $\Theta^{\ast 1}(\mu,x)>0$ for $\mu$ almost all $x\in \Rn$. Hence it characterizes 1-rectifiability of general $\H^1$ measurable sets with $\H^1(E)<\infty$.

We shall discuss the rectifiablity - densities question in metric spaces in Chapter \ref{metricspaces} but let us briefly consider the case where $\Rn$ is equipped with the $l^{\infty}$ norm. So the balls are cubes $Q(x,r)$ with sides parallel to the coordinate axis. This is widely open. Lorent \cite{Lor03} proved a partial result with a very complicated argument: if for a Radon measure $\mu$ in $\R^3$ we have  $\mu(Q(x,r))=r^2$ for all $x\in\spt\mu, r>0$, then $\mu$ is 2-rectifiable.

I still say a few words about generalized densities. First, Marstrand proved in \cite{Mar64} that if $s>0$ and for some $n$ there is $\mu\in\mathcal M(\Rn)$ such that the positive and finite limit $\lim_{r\to 0}r^{-s}\mu(B(x,r))$ exists for $\mu$ almost all $x\in\Rn$, then $s$ is an integer. Preiss made a thorough and deep investigation of the corresponding and related questions for general density functions. The following is a special case of \cite[Theorem 6.5]{Pre87}: 

\begin{thm}\label{preiss3}
Let $h:(0,\infty)\to (0,\infty)$ be such that the limit $\lim_{r\to 0}h(tr)/h(r)$ exists for all $t>0$. Then for some $n$ there is $\mu\in\mathcal M(\Rn)$ such that the positive and finite limit $\lim_{r\to 0}\mu(B(x,r))/h(r)$ exists for $\mu$ almost all $x\in\Rn$ if and only if
\begin{itemize}
\item[(1)] there is an integer $m, 0\leq m\leq n,$ such that $0<\lim_{r\to 0}r^{-m}h(r)<\infty$, or
\item[(2)] there is an integer $m, 1\leq m\leq n-1,$ such that $\lim_{r\to 0}r^{-m}h(r)=0$,\ $\lim_{r\to 0}h(tr)/h(r)= t^m$ for all $t>0$, and 
$\lim_{r\to 0}\sup_{t\in (0,1]}h(tr)/h(r)= 1$.
\end{itemize}
\end{thm}

By an example in \cite[Proposition 6.9]{Pre87} the first condition on $h$ is needed. Preiss called functions $h$ as above exact density functions. A surprising consequence of Theorem \ref{preiss3} is that $r/|\log r|$ is an exact density function but $r|\log r|$ is not. For rectifiability criteria with general density functions, see \cite[Corollary 5.4]{Pre87}

\subsection{Projections}\label{projections}

Federer extended Besicovitch's projection theorem to general dimensions in \cite{Fed47a}:

\begin{thm}\label{profed}
Let $E\subset\Rn$ be $\H^m$ measurable with $\mathcal H^m(E)<\infty$. Then $E$ is purely $m$-unrectifiable if and only if $\H^m(P_V(E))=0$ for almost all $V\in G(n,m)$.
\end{thm}

Federer proved this first for $m=n-1$ using Besicovitch's three alternatives method, which we explained in the previous chapter. Then he used downward induction on $m$ and some integral geometry to get the general case. 


White \cite{Whi98} gave a different proof. He took Besicovitch's result in the plane for granted and showed that if $\H^m(P_V(E))>0$ for positively many $V$, then $E$ is not purely unrectifiable. First he used induction on $n$ to get the result for $m=1$ and all $n$. To do this he applied the induction hypothesis to $P_W(E)$ on some suitably chosen hyperplane $W$. To get to $m>1$ he used an elegant argument applying the case $m=1$ to the intersections $X\cap E$ of $E$ with appropriate affine $(n-m+1)$-planes $X$. The planes $W$ and $X$ are found by some integral geometry. This is only a rough imprecise idea.

Jones, Katz and Vargas \cite{JKV97} also gave another proof in the case $m=n-1$. They used induction on $n$ beginning with Besicovitch's result.

O'Neil proved in \cite{One96} a local version of the projection theorem: if lower and upper $m$-densities of a measure $\mu$ are positive and finite and all projections  on $m$-planes of the supports of the tangent measures are convex, then $\mu$ is $m$-rectifiable.

We also have the Crofton formula: if $E\subset\Rn$ is an $m$-rectifiable Borel set with $\H^m(E)<\infty$, then
$$\H^m(E)= c(n,m)\int_{G(n,m)}\int_V\card(E\cap P_V^{-1}\{a\})\, d\H^ma\, d\gamma_{n,m}V.$$
This can be stated in terms of the \emph{integral-geometric measure} $\mathcal I_1^m$: $\H^m(E)=\mathcal I_1^m(E)$ for rectifiable sets. There is a continuum of integral-geometric measure $\mathcal I_t^m, 1\leq t \leq\infty$, defined by
$$\mathcal I_t^m(A)=\lim_{\delta\to 0}\inf\{\sum_{i=1}^{\infty}\zeta_{n,m,t}(B_i): A\subset \cup_{i=1}^{\infty}B_i, B_i\ \text{Borel sets}, d(B_i)<\delta\},$$
where
$$\zeta_{n,m,t}(B)=c(n,m,t)\left(\int\H^m(P_V(B))^t\,d\gamma_{n,m}V\right)^{1/t}\ \text{if}\ 1\leq t<\infty,$$
$$\zeta_{n,m,\infty}(B)=\text{esssup}_V\H^m(P_V(B)).$$
For any Borel set $B\subset\Rn$,
$$\mathcal I_1^m(B)=c(n,m)\int_{G(n,m)}\int_V \card(B\cap P_V^{-1}\{a\})\, d\H^ma\, d\gamma_{n,m}V.$$
All these measures agree with $\H^m$ for $m$-rectifiable sets and they have the same null-sets: $\mathcal I_t^m(A)=0$ if and only if $A$ is contained in a Borel set $B$ such that $\H^m(P_V(B))=0$ for almost all $V\in G(n,m)$. Moreover, $\mathcal I_s^m\leq\mathcal I_t^m$ if $s\leq t$. It is also known that if $\mathcal I_{\infty}^m(A)<\infty$, then all $\mathcal I_t^m(A)$ agree, \cite[3.3.16]{Fed69}. 

In \cite{Mat86a} a compact set $F\subset\R^2$ was constructed for which $\mathcal I_1^1(F)<\mathcal I_{\infty}^1(F)=\infty$. But it is not known if, for example, all $\mathcal I_t^m, 1<t<\infty$, agree with $\mathcal I_1^m$ or with $\mathcal I_{\infty}^m$? Here are two theorems relevant for rectifiability. The first is often called structure theorem. It is essentially a restatement of what already was said above. The proof of the second requires still more work, see \cite{Fed69}, 3.3.13 and 3.3.14.

\begin{thm}
If $E\subset\Rn$ with $\H^m(E)<\infty$ and $1\leq t\leq\infty$, then $E=R\cup P$ where $R$ is $m$-rectifiable and $\mathcal I_t^m(P)=0$. In particular, $E$ is $m$-rectifiable if and only if $\H^m(E)=\mathcal I_t^m(E)$.
\end{thm}

\begin{thm}
If $E\subset\Rn$  with $\mathcal I_{\infty}^m(E)<\infty$, then $E$ is $\mathcal I_{\infty}^m$-rectifiable, that is, $\mathcal I_{\infty}^m$ almost all of $\Rn$ can be covered with countably many Lipschitz images of $\R^m$.
\end{thm}

In addition to Hausdorff and integral-geometric measures there are many other natural $m$-dimensional measures which all agree on $m$-rectifiable sets, see \cite[Theorem 3.2.26]{Fed69} and \cite[Theorem 17.11]{Mat95}.

Brothers \cite{Bro69} proved the analogue of Theorem \ref{profed} and studied integral-geometric measures on $n$-dimensional manifolds $X$ with a transitive group $G$ of diffeomorphisms. Let $Y$ be an $(n-m)$-dimensional submanifold. Then the condition that almost all projections of $E$ have $\H^m$ measure zero can be stated as $E\cap g(Y)=\emptyset$ for almost all $g\in G$. Hovila, E. and M. J\"arvenp\"a\"a and Ledrappier \cite{HJJL14} proved the projection theorem for more general transversal families of linear maps $\Rn\to\R^m$ and applied it to invariant measures of geodesic flows on surfaces. Besicovitch's three alternatives are still there in both of these approaches.

\subsection{Multiscale approximations}\label{multiscale}
In Chapters \ref{unifrect} and \ref{Rectmeas} we shall discuss more extensively multiscale approximations in terms of Jones type square functions, recall Section \ref{tsp}, but now let us just state the following result of Azzam and Tolsa from \cite{Tol15b} and \cite{AT15}. Define
\begin{equation}\label{beta2}
\beta_{E}^{m,2}(x,r)=\inf_{V\ \text{affine}\ m-\text{plane}}\left(r^{-m}\int_{E\cap B(x,r)}\left(\frac{d(y,V}{r}\right)^2\,d\H^my\right)^{1/2}.
\end{equation}

\begin{thm}\label{Azzam}
If $E\subset\Rn$ is $\H^m$ measurable and $\mathcal H^m(E)<\infty$, then $E$ is  $m$-rectifiable if and only
\begin{equation}\label{ATeq1}
\int_0^{1}\beta_{E}^{m,2}(x,r)^2r^{-1}dr < \infty\ \text{for}\  \mathcal H^m\ \text{almost all}\ x\in E.
\end{equation}
\end{thm}\

The proof is very technical and complicated. The part that \eqref{ATeq1} implies rectifibility was proved in \cite{AT15}. It uses stopping time arguments where the rough idea is the following.  The assumption \eqref{ATeq1} tells us that large part of $E$ is well approximated by $m$-planes at most scales, so we start with some finite family of generalized dyadic cubes where this happens. Then we go to smaller subcubes and stop when the approximation is not good enough. The stopping cubes contain only a small part of $E$. In the others the approximation becomes better and better. That allows us to build Lipschitz graphs which in the limit tend to a Lipschitz graph which meets $E$ in a set of positive measure. This is a vastly oversimplified sketch. In fact, this type of scheme is not typical only to the above paper. It is quite commonly used in particular in connection of uniform rectifiability, see Section \ref{Basic tools}.  Also the other direction, proved in \cite{Tol15b}, is based on stopping time arguments. There Tolsa first proves the corresponding result with $\a$ numbers, see \eqref{alphaur}. The proofs would be much easier if $E$ had positive lower density, now the stopping also takes place when the density ratios get too small. 

Edelen, Naber and Valtorta \cite{ENV16} proved a sufficient condition for rectifiability which extends that part of Theorem \ref{Azzam}, see Theorem \ref{ENV}. Naber and Valtorta proved in \cite{NV17} and \cite{NV20} closely related quantitative results and applied them to harmonic maps between manifolds and to stationary varifolds, see Chapter \ref{Singularities}.

Many results on rectifiability and uniform rectifiablity hold with the exponent 2 in the $\beta$ numbers replaced by a range of exponents $p$. Rather surprisingly Tolsa showed in \cite{Tol19} that Theorem \ref{Azzam} holds only for $p=2$. With additional density conditions other exponents $p$ work, too, see \cite{Paj97} and \cite{BS16}.

In \cite{AS18} Azzam and Schul obtained higher dimensional analogues of both theorems \ref{jones} and \ref{betathm} using $\beta$-integrals defined in terms of Hausdorff content. For related results, see \cite{Vil19a}, \cite{AV21} and \cite{Hyd21a}. Hilbert space versions were proven by Hyde in \cite{Hyd21b}.

\subsection{Reifenberg type results} \label{Reiftype}
For $0<m<n, x\in\Rn, r>0$ and $E\subset\Rn$ define the $\beta$ number
\begin{equation}\label{reifbeta}\beta_E^m(x,r)=\inf_V\sup_{y\in E\cap B(x,r)}d(y,V)/r,\end{equation}
and the bilateral $\beta$ number
\begin{equation}\label{reifbwg}b\beta_E^m(x,r)=\inf_V\left(\sup_{y\in E\cap B(x,r)}d(y,V)/r+\sup_{y\in V\cap B(x,r)}d(y,E)/r\right),\end{equation}
where the infima are taken over all $m$-planes in $\Rn$. Reifenberg proved in \cite{Rei60}
\begin{thm}\label{reifthm}
For any $0<\a<1$, there is $\delta=\delta(n,\a)>0$ such that if $E\subset B^n(0,2)$ is closed and $b\beta_E^m(x,r)<\delta$ when $x\in E$ and $B(x,r)\subset B^n(0,2)$, then there is an $\a$-bi-H\"older map from $B^m(0,1)$ onto $E\cap B^n(0,1)$. 
\end{thm}

A nice exposition of the proof and related matters is given by Naber in \cite{Nab20a}.

Suppose that $E$ satisfies this condition locally for every $\delta>0$; $b\beta_E^m(x,r)<\delta$ when $x\in E, 0<r<r(\delta)$. Then it follows that $\dim E = m$ but $E$ need not have $\sigma$-finite $\H^m$ measure, in particular, it need not be $m$-rectifiable. This is easily seen by von Koch snowflake type examples for which the angle between the segments in consecutive generations goes to zero slowly. Then the approximating lines can turn around infinitely many times.

So to get rectifiability from Reifenberg type assumptions, involving approximation by planes, we need something more. In a way many results discussed earlier are of this sort, but here I restrict 'Reifenberg type' to mean that in addition to fixing the dimension of the planes, we do not make any other dimensional assumptions, such as with densities. Jones's traveling salesman theorem \ref{jones} is of this type. But Reifenberg type could also refer to results with bijective parametrizations. 

Now we give a Reifenberg type result of Simon, see \cite[Section 4.2]{Sim96}, which he used to prove rectifiability of singularities of minimal surfaces and harmonic maps. We shall return to this in Chapter \ref{Singularities}.

The formulation of Simon's theorem is a bit complicated, so I state a simpler special case. The assumptions of the actual result allow at each scale a small exceptional set which seems to be essential in the applications.

\begin{thm}\label{simrect}
For any $0<\delta<1$, there is $\e=\e(n,\delta)>0$ such that the following holds. Let $E\subset \Rn$ be closed such that $\beta_E^m(x,r)<\e$ when $x\in E$ and $0<r<1$ and suppose that $E$ has the following property: Let $x_0\in E, 0<r_0<1$, and $V\in G(n,m)$ for which  $d(x,V+x_0)\leq\e r_0$ for all $x\in E\cap B(x_0,r_0)$. If $x\in E\cap B(x_0,r_0)$ and $0<r<r_0$ are such that $E\cap B(y,\delta s)\neq\emptyset$ for all $y\in (V+x)\cap B(x,s), r\leq s\leq r_0$, then  $d(y,V+x)\leq\e s$ for all $y\in E\cap B(x,s), r\leq s\leq r_0$. Then $E$ is $m$-rectifiable.
\end{thm}

That is, the main assumptions say something like this: if $E$ is $\e$-well approximated in some ball with some plane, then it continues to be $\e$-well approximated in smaller balls with translates of the same plane as long as there are no $\delta$-gaps, that is, as long as there is bilateral $\delta$-approximation. So if at a generic point the approximating plane can turn (wildly, as for von Koch-type examples) only because of the gaps, then the set is rectifiable (and the plane cannot turn wildly).

To see why this might be true, observe first that if there are no gaps at all, there is approximation at all scales with planes parallel to a fixed plane. This easily implies that $E$ is contained in a Lipschitz graph. Secondly, if $E\cap B(x,r)$ is contained in an $\e r$ neighbourhood of an $m$-plane through $x$ and there is a $\delta$-gap, then, since $\e$ is much smaller than $\delta$, $E\cap B(x,r)$ can be covered with balls $B_i$ for which $\sum_id(B_i)^m < \lambda (2r)^m$ with $\lambda<1$ depending on $\delta$. Thus gaps at many places and scales leads to small measure.

David and Toro found in \cite{DT12} conditions which, when added to the assumptions of Theorem \ref{reifthm}, guarantee that the map $f$ can be  chosen to be bi-Lipschitz. One such condition is that $\sum_k\beta_E^m(x,10^{-k})^2$ is bounded. They also discussed relations to uniform rectifiability. The methods are partially based on the earlier work of Toro \cite{Tor95}. In \cite{Mer17} Merhej, generalizing a result from \cite{Tor95}, established sufficient conditions for the bi-Lipschitz parametrization of a codimension one AD-regular  rectifiable set in terms of the Poincar\'e inequality and quadratic oscillation of the unit normal.

Naber and Valtorta used in \cite{NV17} a $\beta$ assumption to get a $W^{1,p}, p>n,$ parametrization. Edelen, Naber and Valtorta \cite{ENV16} gave a Reifenberg-type result for measures, see Theorem \ref{ENVreif}. Results of this type have been applied to the structure of singularities, see Chapter \ref{Singularities}. 

\subsection{Lebesgue null-sets and singular measures}\label{Lebnull}
Many of the results of this section are described in \cite{ACP05} and \cite{ACP10}, but the full proofs have not yet been published.

This begins when Preiss \cite{Pre90} discovered a set $A$ in the plane of zero Lebesgue measure such that every Lipschitz function $f:\R^2\to\R$ is differentiable at some point of $A$. After that a lot of work has been done by Alberti, Cs\"ornyei, Preiss and others on the differentiability properties of Lipschitz maps on null-sets. I don't go into that here, but see 
the surveys \cite{ACP05} and \cite{ACP10}, and \cite{AM16}. I shall concentrate on geometric properties of Lebesgue null-sets and general singular measures.

Preiss's set is simple to state: any $G_{\delta}$ null-set of the plane containing the countable set of lines with rational coordinates is fine. Then for a given Lipschitz function there are a lot of directional derivatives. But to get differentiability, one should be able to combine them to a derivative mapping. The following definition has turned out to be relevant for differentiability and other questions:

\begin{df}\label{weaktf}
A Borel mapping $\tau:E\to G(n,m)$ is a \emph{weak $m$-tangent field} of a set $E\subset\Rn$ if for every $m$-rectifiable set $F$ with $\H^m(F)<\infty$,\ $\apTan(F,x)=\tau(x)$ for $\H^m$ almost all $x\in E\cap F$.
\end{df}

Let $E$ be $\H^m$ measurable and $\H^m(E)<\infty$. If $E$ is $m$-rectifiable, then $E$ has an $\H^m$  unique weak $m$-tangent field, while if $E$ is purely $m$-unrectifiable, every $\tau:E\to G(n,m)$ is a weak $m$-tangent field of $E$. In general, a weak $m$-tangent field is unique up to purely $m$-unrectifiable sets.

Alberti, Cs\"ornyei and Preiss have proven

\begin{thm}\label{ACP}
Any set $E\subset\Rn$ with $\mathcal L^n(E)=0$ admits a weak $(n-1)$-tangent field.
\end{thm}

In the words of the authors of \cite{ACP10}: "This result can be understood as saying the rather mysterious fact that one can
prescribe in which direction an $(n-1)$-surface meets a null set $E$, without knowing the surface itself."

Integral representations with rectifiable measures play an important role in the investigations of singular measures. They were used by Alberti in \cite{Alb93} for the proof of his rank one theorem \ref{Alberti} for BV-functions, where he established deep analytic and geometric properties  of singular measures. Nowadays they are called Alberti representations by many authors.

Alberti, Cs\"ornyei and Preiss presented the following general definitions and results in \cite{ACP10}:
\begin{df}\label{rectrepr}
A measure $\nu\in\mathcal M(\Rn)$ is called \emph{$m$-rectifiably representable} if it can be written as $\nu=\int\mu_t\,dPt$ where  $\mu_t\ll \H^m\restrict E_t$ for some $m$-rectifiable set $E_t$ with $\H^m(E_t)<\infty$ and $P$ is some probability measure. For such a $\nu$ we say that $\tau:\Rn\to G(n,m)$ is an \emph{$m$-tangent field} of $\nu$ if $\apTan(E_t,x)=\tau(x)$ for $\H^m$ almost all $x\in E_t$ and $P$ almost all $t$.
\end{df}

\begin{thm}\label{ACP1}
Let $\mu\in\mathcal M(\Rn)$. Then 

$\mu$ is $m$-rectifiably representable if and only if $\mu(E)=0$ for every purely $m$-unrectifiable set $E\subset\Rn$.

If $\mu$ is  $(n-1)$-rectifiably representable, then it admits an $(n-1)$-tangent field if and only if it is singular.

$\mu$ has a unique decomposition as $\mu=\mu_n+\mu_{n-1}\dots+\mu_0$, where each $\mu_m$ is $m$-rectifiably representable and it lives on a  purely $(m+1)$-unrectifiable set.
\end{thm}

We shall return to tangent fields and Alberti representations in Chapter \ref{metricspaces}.

Cs\"ornyei, Preiss and Tiser \cite{CPT05} and Maleva and Preiss \cite{MP19} introduced large subclasses of purely 1-unrectifiable sets such that for any set in these classes some Lipschitz function is non-differentiable at every point of it. They used them to describe many other detailed (non-)differentiability properties of Lipschitz functions, too.

Alberti and Marchese \cite{AM16} introduced the decomposition bundle $V(\mu,x), x\in\Rn,$ of any measure $\mu\in\mathcal M(\Rn)$: $V(\mu,x)$ is the smallest linear subspace of $\Rn$ with the following property: if $\nu=\int\mu_t\,dPt$ is 1--rectifiably representable, as in Definition \ref{rectrepr}, and $\nu\ll\mu$, then $\apTan(E_t,x)\subset V(\mu,x)$ for $\H^1$ almost all $x\in E_t$ and $P$ almost all $t$. They used it to obtain interesting Lipschitz differentiability results and they also applied it to normal currents. 

We have by a fairly easy result of \cite[Proposition 2.9]{AM16}:

\begin{pr}\label{AMpr}
$\mu$ is purely 1-unrectifable if and only if $V(\mu,x)=\{0\}$ for $\mu$ almost all $x\in\Rn$.
\end{pr}

Del Nin and Merlo \cite{DM21} found a nice application of the decomposition bundles. They proved a dichotomy in Fourier restriction:
\begin{thm}\label{DM21}
If $0<s\leq n$, $\mu\in\mathcal M(\Rn)$ with $0<\Theta^{\ast s}(\mu,x)<\infty$ for $\mu$ almost all $x\in\Rn$ and the restriction inequality
\begin{equation}\label{Fourier}
\|\hat f\|_{L^q(\mu)}\lesssim \|f\|_{L^p(\Rn)}\end{equation}
holds when $q=sp'/n$, with $p'=p/(p-1)$, then either $q=p'$, that is, $s=n$, and so $\mu\ll\mathcal L^n$, or $\mu$ is purely 1-unrectifiable.
\end{thm}
The value $q=sp'/n$ is the end-point, the largest value for which this restriction could hold under the density assumption.

To prove theorem \ref{DM21} one first checks by direct computation that \eqref{Fourier} is preserved for tangent measures. At $\mu$ almost all points $x$ the tangent measures are shown to be of the form $\nu\otimes\H^{k(x)}\restrict V(\mu,x), \nu\in\mathcal M(V(\mu,x)^{\perp})$, where $k(x)=\dim V(\mu,x)$. If $k(x)>0$, the restriction estimate for 
$\H^{k(x)}\restrict V$ follows, which, by trivial scaling, is only possible if $q=p'$. Hence the theorem follows from Proposition \ref{AMpr}.

Another application of of the decomposition bundles was found by Marchese and Merlo in \cite{MM22}. They showed that the Lusin type approximation property \eqref{Lusin} holds for a Radon measure $\mu$ on $\Rn$ in place of the Lebesgue measure if and only $\mu$ is a sum of absolutely continuous (with respect to the corresponding Hausdorff measures) rectifiable measures of various dimensions.

\subsection{Minkowski content and discrete energies}

For $0<m<n$ the \emph{lower and upper Minkowski contents} are defined for $A\subset\Rn$ by
$$\mathcal M^m_{\ast}(A)=\liminf_{\delta\to0}\a(n-m)^{-1}\delta^{m-n}\mathcal L^n(\{x\in A: d(x,A)<\delta\}),$$
$$\mathcal M^{\ast m}(A)=\limsup_{\delta\to0}\a(n-m)^{-1}\delta^{m-n}\mathcal L^n(\{x\in A: d(x,A)<\delta\}).$$
If they agree, their common value $\mathcal M^{m}(A)$ is the Minkowski content of $A$.

It is easy to see that $\H^m(A)\lesssim \mathcal M^m_{\ast}(A)$ and that there are compact sets $A$ for which $\H^m(A)=0$ and $\mathcal M^m_{\ast}(A)=\infty$. In particular, when $m$ is an integer, such an $A$ is $m$-rectifiable. Federer proved in \cite[Theorem 3.2.39]{Fed69}

\begin{thm}\label{Mink}
If $0<m<n$ are integers and a closed set $F\subset\Rn$ is a Lipschitz image of a bounded subset of $\R^m$, then $\mathcal M^{m}(F)=\mathcal H^{m}(F)$.
\end{thm}

Borodachov,  Hardin and  Saff used in \cite{BHS08} and \cite{BHS14} this connection to prove a very general result on the asymptotics of discrete $s$-energies.  For $s>0$ and a compact subset $F$ of $\Rn$ define
$$\mathcal E_s(F,N)=\min\{\sum_{i\neq j}|x_i-x_j|^{-s}:x_1,\dots,x_N\subset F\}.$$
For a given $F$ the question how does $\mathcal E_s(F,N)$ behave when $N\to\infty$ has been extensively studied. The most classical cases are when $F$ is a sphere and $s=n-2$, when $n\geq 3$, and $|x|^{-s}$ is replaced by $-\log|x|$, when $n=2$. See the book  \cite{BHS19} for a huge amount of generalizations and related questions, connections and applications to a wide variety of topics. Here I only discuss one result related to rectifiability.

Let $m=\dim F$ be the Hausdorff dimension of $F$. The cases $s<m$ and $s\geq m$ are quite different. In the first case by classical potential theory, see \cite[Theorem 4.2.2]{BHS19},
$$\lim_{N\to\infty}N^{-2}\mathcal E_s(F,N)=\iint|x-y|^{-s}\,d\mu_{m,F}x\,d\mu_{m,F}y$$
where $\mu_{s,F}\in\mathcal M(F)$ is an equilibrium probability measure which minimizes the energies $I_s(\mu):=\iint|x-y|^{-s}\,d\mu x\,d\mu y$ among all probability measures $\mu\in\mathcal M(F)$. The second case is much more delicate, because $I_m(\mu)=\infty$ for all $\mu\in\mathcal M(F)$, at least if $\H^m(F)>0$. We have by \cite{BHS08} and \cite{BHS14}, see also \cite[Theorem 8.5.2]{BHS19}
\begin{thm}\label{BHS}
Let $F\subset\Rn$ be compact and $s>m=\dim F$. If $F$ is $m$-rectifiable and $\mathcal M^{m}(F)=\mathcal H^{m}(F)$, then
$$\lim_{N\to\infty}N^{-1-s/m}\mathcal E_s(F,N)=C_{s,m}\H^m(F)^{-s/m}.$$
Moreover, if $\H^m(F)>0$ and $\{x_{N,1},\dots,x_{N,N}\}$ is a minimizing configuration for $\mathcal E_s(F,N)$, then
\begin{equation}\label{BHSeq}
\lim_{N\to\infty}\frac{1}{N}\sum_{i=1}^N\delta_{x_{N,i}}\to \frac{1}{\H^s(F)}\H^s\restrict F\ \text{weakly as}\ N\to\infty.\end{equation}
\end{thm}

Since the constant $C_{s,m}$ does not depend on $F$ one can derive a lot of information about it looking, for example, at the case where $F$ is an $m$-sphere or a cube in $\R^m$.

When $s=m$ similar results hold for compact subsets of $C^1$ submanifolds of $\Rn$, and a little more generally, see \cite[Theorem 9.5.4]{BHS19}. Then $N^{1+s/m}$ is replaced by $N^2\log N$. It seems to be unknown whether they are valid under rectifiability assumptions as in Theorem \ref{BHS}.

In the case $s<m$ the cluster points of the extremal measures as in \eqref{BHSeq} are equilibrium measures $\mu_{s,F}$, which in general are not a multiples of Hausdorff measures. 
For example, when $F$ is a ball in $\R^m$ and $m-2<s<m$, the equilibrium measure is absolutely continuous with density going to infinity at the boundary.

For recent related results, see \cite{HSV22}.

\section{Uniform rectifiability}\label{unifrect}
The various rectifiability characterizations we have seen are qualitative. This is unavoidable if we start with qualitative size conditions of  measure and densities. But if we start with quantitative size conditions like AD-regularity, then quantitative equivalent uniform rectifiability conditions are in the core of the theory developed by David and Semmes in \cite{DS91} and \cite{DS93}. Most of the material of this chapter is from those books.

\subsection{One-dimensional sets}

Let us first have a quick look at one-dimensional sets. We say that $E\subset\R^n$ is \emph{uniformly 1-rectifiable} if it is closed,  AD-1-regular and it is contained in an AD-regular curve, that is, there is a curve $\Gamma$ and a constant $C$ such that $E\subset\Gamma$ and 
$$\H^1(\Gamma\cap B(x,r))\leq Cr\ \text{for all}\ x\in\Rn, r>0.$$

The key to the related results for approximation by lines is in the traveling salesman type conditions discussed in Section \ref{tsp}. Recall the definition of $\beta_E(x,r)$ from \eqref{beta}.

The following is a local version of Jones's traveling salesman theorem \ref{jones}:

\begin{thm}\label{beta1thm}
If $E\subset\Rn$ is closed and AD-1-regular, then $E$ is uniformly 1-rectifiable if and only if
\begin{equation*}
\int_0^{R}\int_{E\cap B(x,R)}\beta_E(y,r)^2\,d\H^1y\,r^{-1}dr \lesssim R\ \text{for all}\ x\in E, R>0.
\end{equation*}
\end{thm}

We also have a characterization by Menger curvature, recall its definition from \eqref{menger}:

\begin{thm}\label{mengerthm}
If $E\subset\R^2$ is closed and AD-1-regular, then $E$ is uniformly 1-rectifiable if and only if
$$\int_{E\cap B(a,R)}\int_{E\cap B(a,R)} \int_{E\cap B(a,R)} c(x,y,z)^2\,d\H^1x\,d\H^1y\,d\H^1z \lesssim R$$
for all $a\in E, R>0.$
\end{thm}

That uniformly rectifiable sets satisfy this curvature condition was proved by Melnikov and Verdera in \cite{MV96}. The other direction was proved in \cite{MMV96}.

\subsection{Lipschitz maps and approximation by planes}\label{Lipmap-appr}

Now $0<m<n$ will be integers. In higher dimensions there is no as simple definition of uniform rectifiability, but there are many interesting characterizations. Let us again take Lipschitz maps as the basis of the definition:

\begin{df}
We say that $E\subset\R^n$ is \emph{uniformly $m$-rectifiable} if it is closed,  AD-$m$-regular and there are positive numbers $M$ and $\theta$ such that for all $x\in E$ and $0<R<d(E)$ there is a Lipschitz map $f:B^m(0,R)\to\Rn$ with $\Lip(f)\leq M$ and
$$\H^m(E\cap B(x,R)\cap f(B^m(0,R)))\geq \theta R^m.$$
\end{df}

David and Semmes call this property as $E$ having \emph{big pieces of Lipschitz images} of $\R^m$. There is a similar characterization by bi-Lipschitz images but \emph{big pieces of Lipschitz graphs} (graphs over $m$-planes) is strictly stronger by an unpublished Venetian blind construction of Hrycak, see \cite{Azz21} or \cite{Orp21}. However, iterating this, big pieces of big pieces of Lipschitz graphs is equivalent to uniform rectifiability, see \cite{AS12}. This with iteration seems to be a general phenomenon, see Section \ref{Parabolic} and \cite{BHHLN20}. 

It is clear that uniformly rectifiable sets are rectifiable, but not vice versa. Lipschitz graphs are basic examples of uniformly rectifiable sets.

Define the $L^p$ versions of $\beta$s for any $1\leq p<\infty$ by

\begin{equation}\label{betap}
\beta_{E}^{m,p}(x,r)=\inf_{V\ \text{affine}\ m-\text{plane}}\left(r^{-m}\int_{E\cap B(x,r)}\left(\frac{d(y,V)}{r}\right)^p\,d\H^my\right)^{1/p}.
\end{equation}
Set also $\beta_{E}^{m,\infty}(x,r)=\beta_E^m(x,r)$, recall \eqref{reifbeta}.

According to \cite{DS91} Jones and Fang have produced 3-dimensional Lipschitz graphs which show that Theorem \ref{beta1thm} does not hold for $\beta_{E}^{3,\infty}$. But we have the following:

\begin{thm}\label{betamthm}
Let $1\leq p \leq \infty$, if $m=1$, and $1\leq p < 2m/(m-2)$, if $m>1$. If $E\subset\Rn$ is closed and AD-$m$-regular, then $E$ is uniformly $m$-rectifiable if and only
\begin{equation}\label{betap2}
\int_0^{R}\int_{E\cap B(x,R)}\beta_E^{m,p}(y,r)^2\,d\H^my\,r^{-1}dr \lesssim R^m\ \text{for all}\ x\in E, R>0.
\end{equation}
\end{thm}

In addition to being natural as the $L^p$ versions of the uniform approximation of sets, the $\beta$s have counterparts in $L^p$ differentiation of functions due to a result of Dorronsoro in \cite{Dor85}. 




The validity of \eqref{betap2} is often stated as $E$ satisfying the \emph{geometric lemma}. One also says that $\beta_E^{m,p}(x,r)^2\,d\H^mx\,r^{-1}dr$ is a Carleson measure on $E\times (0,\infty)$:

\begin{df}
Let $E\subset\Rn$. A Borel measure $\lambda$ on $E\times (0,\infty)$ is a \emph{Carleson measure} if
$$\lambda(B(x,r)\times (0,R))\lesssim R^m\ \text{for all}\ x\in E, R>0.$$
A set $A\subset E\times (0,\infty)$ is a \emph{Carleson set} if $\chi_A(x,r)\,d\H^mx\,r^{-1}dr$ is a Carleson measure.
\end{df}

This terminology comes from Carleson's solution of the complex analysis corona problem in 1962 and the methods he introduced. Several other Carleson measure and set conditions characterizing uniform rectifiability can be found in \cite{DS93}. Below we shall see some of them and some more recent ones.

Condition \eqref{betap2} guarantees that $E$ is well approximable by planes at most scales, so it is a relative of the existence of approximate tangent planes of rectifiable sets. Next we shall give different conditions in that spirit.

\begin{df}\label{wgl}
Let $E\subset\Rn$ be AD-$m$-regular. We say that $E$ satisfies \emph{weak geometric lemma} if for every $\e>0$,
$$\{(x,r)\in E\times (0,\infty): \beta_{E}^{m}(x,r)>\e\}$$
is a Carleson set, that is,
$$\int_0^{R}\H^m(\{y\in E\cap B(x,R):\beta_{E}^{m}(y,r)>\e\})r^{-1}dr \lesssim R^m\ \text{for all}\ x\in E, R>0.$$ 
\end{df}

Weak geometric lemma does not imply even ordinary rectifiability, recall Section \ref{Reiftype}. But it is useful in combination with some other conditions. As in the case of rectifiability and tangent measures, the corresponding biliteral approximation does the job. Recall the bilateral $\beta$ from \eqref{reifbwg}. The \emph{bilateral weak geometric lemma} characterizes uniform rectifiability:

\begin{thm}\label{bwgl}
If $E\subset\Rn$ is closed and AD-$m$-regular, then $E$ is uniformly $m$-rectifiable if and only for every $\e>0$,
$$\{(x,r)\in E\times (0,\infty): b\beta_{E}^m(x,r)>\e\}$$
is a Carleson set. 
\end{thm}

A weaker condition where approximation is allowed with unions of planes already implies uniform rectifiability. To state this define $ub\beta_E(x,r)$ as in  \eqref{reifbwg} but $V$ replaced by a union of $m$-planes. Then, see \cite{DS93}, Proposition II.3.18,

\begin{thm}\label{ubwgl}
If $E\subset\Rn$ is closed and AD-$m$-regular, then $E$ is uniformly $m$-rectifiable if and only for every $\e>0$,
$$\{(x,r)\in E\times (0,\infty): ub\beta_{E}(x,r)>\e\}$$
is a Carleson set. 
\end{thm}

I have stated this a bit more exotic characterization since it will be useful in Chapter \ref{singular integrals}.

There is also a local symmetry characterization. Recall a similar condition in Lemma \ref{refl} and its role in the proof of Theorem \ref{densthm1a}.

\begin{thm}\label{LSthm}
If $E\subset\Rn$ is closed and AD-$m$-regular, then $E$ is uniformly $m$-rectifiable if and only for every $\e>0$,
$$\{(x,r)\in E\times (0,\infty): \exists y,z\in E\cap B(x,r)\ \text{such that}\ d(2y-z,E)>\e r\}$$
is a Carleson set. 
\end{thm}

The local symmetry condition follows immediately from the bilateral approximation in Theorem \ref{bwgl}, the converse requires more work.

Instead of approximating sets by planes one can approximate the Hausdorff measure on $E$ with Lebesgue measures on planes. It is more natural to state this for measures. So we define an AD-$m$-regular measure $\mu$ to be uniformly $m$-rectifiable if its support is uniformly $m$-rectifiable. For a ball $B(x,r)$ and $\mu, \nu\in \mathcal M(B(x,r))$ define
\begin{equation}\label{measmetr}
F_{x,r}(\mu,\nu)=\sup\{|\int f\,d\mu - \int f\,d\nu|: \spt f\subset B(x,r), \Lip(f)\leq 1\}.\end{equation}
Then $F_{x,r}$ is a metric which metrizes weak convergence, see, e.g., \cite[p. 195]{Mat95}. Tolsa  introduced the following $\a$ coefficients in \cite{Tol09}  for a measure $\mu$:
\begin{equation}\label{alphaur}
\a_{\mu}^{m}(x,r)=r^{-m-1}\inf\{F_{x,r}(\mu,c\H^m\restrict V): c\geq 0, V\ \text{an affine}\ m-\text{plane}\}.\end{equation}
He proved  the following: 

\begin{thm}\label{ATTthm}
If $\mu\in\mathcal M(\Rn)$ is AD-$m$-regular, then $\mu$ is uniformly $m$-rectifiable if and only if
\begin{equation}\label{betap1}
\int_0^{R}\int_{B(x,R)}\a_{\mu}^{m}(x,r)^2\,d\mu x\,r^{-1}dr \lesssim R^m\ \text{for all}\ x\in \Rn, R>0.
\end{equation}
\end{thm}

It is easy to show that the $\a^{m}$ numbers dominate the $\beta^{m,1}$ numbers, so one direction follows from Theorem \ref{betamthm}. The other direction uses corona decompositions, see Section \ref{Basic tools}.

In \cite{Tol12} Tolsa proved an $L^p, 1\leq p\leq 2,$ version of this, the case $p=1$ is essentially the above. Then the $\a$ coefficients are defined using the mass transport, Wassertein, distance
$$W_p(\mu,\nu)=\inf_{\pi}\left(\int|x-y|^p\,d\pi(x,y)\right)^{1/p},$$
where the infimum is taken over all probability measures $\pi$ on $\Rn\times\Rn$ whose marginals are $\mu$ and $\nu$. Dabrowski \cite{Dab20a}, \cite{Dab21b} proved the similar characterization of (non-uniform) rectifiability.

In \cite{AD20} Azzam and Dabrowski gave a characterization of the $L^p$ norms $\|f\|_{L^p(\mu)}$ for uniformly rectifiable measures $\mu$ in terms of the $\a$s.

\subsection{Density ratios} 
There is also an analogue of Preiss's theorem \ref{preiss} 'existence of density is equivalent to rectifiability':

\begin{thm}\label{urdens}
If $\mu\in\mathcal M(\Rn)$ is AD-$m$-regular, then $\mu$ is uniformly $m$-rectifiable if and only if the complement in $\spt\mu\times (0,\infty)$ of the set 
\begin{align*}
&\{(x,r)\in \spt\mu\times (0,\infty):\text{there exists}\ \delta(x,r)>0\ \text{such that}\\ &
|\mu(B(y,t))-\delta(x,r)t^m|<\e r^m\ \text{for all}\ y\in \spt\mu\cap B(x,r), 0<t\leq r\}\end{align*}
is a Carleson set for all $\e>0$.
\end{thm}

The 'if' part is the more difficult one. The proof for that is based on uniform measures. It can be shown that the condition of the theorem implies that $\mu$ can be approximated at most locations and scales by $m$-uniform measures. For $m=1,2,n-1$ the $m$-uniform measures are either flat or conical as was discussed after Preiss's theorem \ref{preiss}. Using this David and Semmes completed the proof of Theorem \ref{urdens} for $m=1,2,n-1$ in \cite{DS93}, they had proven the 'only if' direction already \cite{DS91}. For the remaining dimensions no such concrete information about uniform measures is available. However, Preiss's paper contains enough information (uniform measures are 'flat at infinity') so that Tolsa could finish the proof in \cite{Tol15a}.


Of course, Theorem \ref{urdens} implies that uniform measures are uniformly rectifiable.  Tolsa showed more in the same paper: the uniform measures have the 'big pieces of Lipschitz graphs' property.

Chousionis, Garnett, Le and Tolsa \cite{CGLT16} characterized uniform rectifiability with the differences of density ratios, as in Theorem \ref{TT}:

\begin{thm}\label{CGLT}
Let $\mu\in\mathcal M(\Rn)$ be AD-$m$-regular. Then $\mu$ is uniformly $m$-rectifiable if and only if for all $x\in\spt\mu, R>0$,

$$\int_0^R\int_{B(x,r)}\left|\frac{\mu(B(y,r))}{r^m}-\frac{\mu(B(y,2r))}{(2r)^m}\right|^2\,d\mu y\,r^{-1}dr \lesssim R^m.$$ 
\end{thm}

Azzam and Hyde \cite{AH20} proved a sufficient condition for uniform rectifiability in terms of density ratios involving Hausdorff content.

\subsection{Projections}

There is no satisfactory characterization of uniform rectifiability in terms of projections. Tao \cite{Tao09} proved a quantitative multiscale version of Besicovitch's projection theorem \ref{proj1} quantizing Besicovitch's proof. However, this does not seem to relate to David-Semmes uniform rectifiability. But we have the following theorem due to Orponen \cite{Orp21}:

\begin{thm}\label{orpproj}
If $E\subset\Rn$ is closed and AD-$m$-regular, then $E$ has big pieces of Lipschitz graphs if and only there is $\theta>0$ such that for every $x\in E$ and $0<r<d(E)$ there is $V\in G(n,m)$ for which 
\begin{equation}\label{BP}
\H^m(P_W(E\cap B(x,r))\geq \theta r^m\ \text{for}\ W\in B(V,\theta).
\end{equation}
\end{thm}

It is easy to see that big pieces of Lipschitz graphs implies \eqref{BP}: if $G$ is a Lipschitz graph over $V$, then it is also a Lipschitz graph over $W$ when $W$ is sufficiently close to $V$.  David and Semmes showed earlier that one projection satisfying \eqref{BP} is enough if $E$ in addition satisfies the weak geometric lemma, Definition \ref{wgl}. Obviously one projection alone does not imply even rectifiability, think about the four corners Cantor set of Example \ref{Cantor}. Orponen showed that \eqref{BP} implies the weak geometric lemma. His complicated argument involves, among many other things, ingredients from the classical Besicovitch-Federer argument for the proof of Theorem \ref{profed}. The example of Hrycak mentioned in Section \ref{Lipmap-appr} shows that uniformly rectifiable sets need not satisfy \eqref{BP}.

For a related result of Martikainen and Orponen, see \cite{MO18b}. Dabrowski and Villa \cite{DV22} proved an analyst’s traveling salesman theorem for sets satisfying \eqref{BP}.





\subsection{Basic tools}\label{Basic tools}
Many of the proofs use generalized dyadic cubes constructed by David in \cite{Dav91}, another construction was given by Christ in \cite{Chr90b}. The standard dyadic cubes usually are not good enough and they are  replaced by a family $\Delta$ of Borel subsets of an AD-$m$-regular set $E$. In case $E$ is unbounded, which is not essential, $\Delta$ splits into subfamilies $\Delta_j, j\in\Z$. Each $\Delta_j$ is a disjoint partition of $E$, $\H^m(Q)\sim 2^{jm}$ for $Q\in \Delta_j$, and if $Q\in \Delta_j, Q'\in\Delta_k, j\leq k,$ then either $Q\subset Q'$ or $Q\cap Q'=\emptyset$. This is very much as for the standard dyadic cubes but there is a fourth more special 'small boundary' property: there is $\a>0$ such that for $Q\in \Delta_j$ and for all $0<\tau<1$,
$$\H^m(\{x\in Q:d(x,E\setminus Q)\leq \tau 2^j\})+\H^m(\{x\in E\setminus Q:d(x,Q)\leq \tau 2^j\})\lesssim\tau^{\a}2^{jm}.$$
This is more rarely used, but it is useful in particular in connection of singular integrals. 

Often the proofs use stopping time arguments in the spirit we discussed in Section \ref{multiscale}. This leads to the general concept of corona decompositions.  

In the \emph{corona decomposition} the family $\Delta$ as above is decomposed into a good subfamily $\mathcal G$ and a bad subfamily $\mathcal B$. Further, $\mathcal G$ is decomposed into stopping time families $\mathcal S$. Each of them has a unique maximal top cube $Q(\mathcal S)$ which contains all other cubes of $\mathcal S$. In addition, if $Q\in\mathcal S$ and $Q\subset Q'\subset Q(\mathcal S)$, then $Q'\in\mathcal S$, and either all children of $Q$ belong to $\mathcal S$ or none of them do. Then we say that $E$ admits a corona decomposition if for all positive numbers $\eta$ and $\theta$ such corona decomposition can be found with the following two properties:

The bad cubes and the maximal top cubes satisfy the Carleson packing condition for every $Q\in\Delta$:
$$\sum_{Q'\subset Q, Q'\in \mathcal B}\H^m(Q') + \sum_{\mathcal S:Q(\mathcal S)\subset Q}\H^m(Q(\mathcal S))\lesssim \H^m(Q).$$

For every stopping time family $\mathcal S$ there exists a Lipschitz graph $\Gamma(\mathcal S)$ such that $\Lip(\Gamma(\mathcal S))\leq \eta$ and $d(x,\Gamma(\mathcal S))\leq\theta d(Q)$ whenever $x\in 2Q$ and $Q\in \mathcal S$.

David and Semmes proved in \cite{DS91}

\begin{thm}\label{corona}
If $E\subset\Rn$ is closed and AD-$m$-regular, then $E$ is uniformly rectifiable if and only if it admits a corona decomposition.
\end{thm}

This often is a useful link for going from one characterizing proper to another. For instance, in \cite{DS91} David and Semmes proved that the geometric lemma implies corona decomposition which implies big pieces of Lipschitz images.

\subsection{Parabolic rectifiability}\label{Parabolic}

Parabolic uniform rectifiability was introduced by Hofmann, Lewis and Nystr\"om in \cite{HLN03} and \cite{HLN04}. Later there have been several related papers, but now I only comment the recent work of Bortz, Hoffman, Hofmann, Luna Garcia and Nystr\"om in \cite{BHHLN20} and \cite{BHHLN21a}, see the discussions and references there for other developments. Although there are deep results, this theory is not yet very far developed and mainly the codimension one case has been considered.  

In the following balls, Hausdorff measures and AD-regularity are defined with the parabolic metric $d$ in $\R^{n+1}=\Rn\times\R$:
$$d((x,s),(y,t))=|x-y|+\sqrt{|s-t|},\ (x,s), (y,t)\in \Rn\times\R.$$
Notice that the Hausdorff dimension of $\R^{n+1}$ is then $n+2$ and the codimension one Hausdorff measure is $\H^{n+1}$. Define the parabolic $\beta$ by
\begin{equation}\label{betapar}
p\beta_{E}(x,r)=\inf_{V}\left(r^{-n-1}\int_{E\cap B(x,r)}\left(\frac{d(y,V)}{r}\right)^2\,d\H^{n+1}y\right)^{1/2},
\end{equation}
where the infimum is taken over all vertical hyperplanes $V$, that is,  the planes $V=W\times\R$ where $W$ is an affine $(n-1)$-plane in $\Rn$. Using only vertical hyperplanes is natural since they and the 'horizontal' plane $H=\{(x,t):t=0\}$ are the only linear hyperplanes invariant under the dilations $\delta_r(x,t)=(rx,r^2t)$. Notice that $d(\delta_r(x,t),\delta_r(x',t'))=rd((x,t),(x',t'))$.

\begin{df}\label{paraur}
A closed set $E\subset\R^{n+1}$ is \emph{parabolic uniformly rectifiable} if it is AD-$(n+1)$-regular and 
\begin{equation}\label{paraur1}
\int_0^{R}\int_{E\cap B(x,R)}p\beta_{E}(x,r)^2\,d\H^{n+1} x\,r^{-1}dr \lesssim R^{n+1}\ \text{for all}\ x\in E, R>0.
\end{equation}
\end{df}

For David-Semmes uniformly rectifiable sets Lipschitz graphs are basic examples. But not here. There are examples of parabolic Lipschitz graphs which are not uniformly rectifiable. Basic parabolic uniformly rectifiable sets are \emph{regular Lipschitz graphs}, that is, graphs over vertical planes of regular Lipschitz functions $g$. Here  $g:\R^{k}\times\R\to\R$ is a regular Lipschitz function if it is Lipschitz in the parabolic metric $d$ and the half derivative with respect to $t$ belongs to BMO, with BMO defined using parabolic balls. The half derivative  can be defined by
$$D^t_{1/2}g(x,t)=\int_{\R}\frac{g(x,s)-g(x,t)}{|s-t|^{3/2}}\,ds.$$
See \cite{DDH18} for several characterizations of what it means that $D^t_{1/2}g\in BMO$. 

In \cite{BHHLN20} the authors introduced a very general setting for corona decompositions in metric spaces. They proved for any family $\mathcal E$ of AD-$m$-regular sets with bounded constants that if an AD-$m$-regular set $E$ admits a corona decomposition with $\mathcal E$ (that is, $\mathcal E$ replaces Lipschitz graphs in the  definition of Section \ref{Basic tools}), then $E$ has big pieces of big pieces of $\mathcal E$. Combining this with \cite{BHHLN21a} they have

\begin{thm}\label{paraurthm}
If $E\subset\Rn\times \R$ is closed and AD-$(n+1)$-regular, then the following are equivalent:
\begin{itemize}
\item[(1)] $E$ is parabolic uniformly rectifiable.
\item[(2)] $E$ admits a corona decomposition with regular Lipschitz graphs.
\item[(3)] $E$ has big pieces of big pieces of regular Lipschitz graphs.
\end{itemize}
\end{thm}


The study of parabolic rectifiability is motivated by desire to develop the boundary behaviour theory for the heat equation
$\Delta_xu(x,t) = \partial_tu(x,t)$ 
in the same vein as has been done for the Laplace equation, see Chapter \ref{Harmonic, elliptic measures}. In particular, the preference of  regular Lipschitz graphs over ordinary Lipschitz graphs in the parabolic metric comes from heat equation examples.

What should parabolic (non-uniform) rectifiability mean? In the light of the above a natural definition would seem to be the one based on covering with regular Lipschitz graphs. This was suggested in \cite{MPT20} but it has not yet been developed. It would probably be the right notion from analysis (heat equation) point of view. But we could also ask for a notion that would correspond to almost everywhere existence of approximate tangent planes and to tangent measures and this turns out to be different. This question has been studied in \cite{Mat21}. Recall the remarks after Definitions \ref{m-apprtan} and \ref{tanmeasdf} according to which approximate tangent planes and tangent measures can be defined in the parabolic setting as there.  

Let now $m$ be any integer with $1\leq m\leq n+1$. The linear subspaces of $\R^{n+1}$ of the parabolic Hausdorff dimension $m$ that are invariant under the parabolic dilations $\delta_r$ are again the vertical hyperplanes when $m=n+1$, the lines through 0 in $H=\{(x,t):t=0\}$ when $m=1$, and the linear $m$-dimensional subspaces of $H$ and linear vertical $(m-1)$-dimensional subspaces of $\Rn\times\R$ when $2\leq m\leq n$. Let us denote by $P(n,m)$ the family of such linear subspaces. The following theorem was proved in \cite{Mat21}. All the concepts are parabolic.

\begin{thm}\label{pararect}
Let $E\subset \Rn\times \R$ be $\H^m$ measurable and $\mathcal H^m(E)<\infty$. Then the following are equivalent:
\begin{itemize}
\item[(1)] For every $\e>0$ there are Lipschitz graphs $G_i$ over $V_i\in P(n,m)$ with Lipschitz constants less than $\e$ such that $\H^m(E\setminus \cup_iG_i)=0$.
\item[(2)] There are $C^1$ graphs $G_i$ over $V_i\in P(n,m)$ such that $\H^m(E\setminus \cup_iG_i)=0$.
\item[(3)] $E$ has an approximate tangent $m$-plane $V\in P(n,m)$ at $\H^m$ almost all of its points.
\item[(4)] For $\H^m$ almost all $a\in E$\ there is an $m$-flat measure $\lambda_a=\H^m\restrict V_a, V_a\in P(n,m),$ such that $\tanm(\H^m\restrict E,a)=\{c\lambda_a: 0<c<\infty\}$.
\item[(5)] For $\H^m$ almost all $a\in E$\ $\H^m\restrict E$ has a unique tangent measure at $a$.
\end{itemize}
\end{thm}

The proof has many similar ingredients as in the Euclidean case. That (5) implies (4) follows from Theorem \ref{M2}. Here $C^1$ is defined in a parabolic sense. One candidate for parabolic rectifiability is given by these conditions.

For $2\leq m\leq n$ the planes in $P(n,m)$ have linear dimension $m$ or $m-1$, but they can be treated simultaneously.

We need Lipschitz graphs with arbitrarily small Lipschitz constants because Lipschitz graphs themselves need not satisfy the other conditions, as shown in \cite{Mat21}. There one also finds an example which satisfies the conditions of Theorem \ref{pararect} but intersects every regular Lipschitz graph in measure zero. So these two possible classes of rectifiable sets are different.

\section{Rectifiability of measures}\label{Rectmeas}
\subsection{Some basic facts and examples}

Badger \cite{Bad19} has a nice survey covering parts of this topic. Recall that we have defined a measure $\mu\in\mathcal M(\Rn)$ to be $m$-rectifiable if there are Lipschitz maps $f_i:\R^m\to\Rn$ such that 
$$\mu(\Rn\setminus\bigcup_{i=1}^{\infty}f_i(\R^m))=0,$$
without requiring absolute continuity with respect to $\H^m$. 

Studying sets $E$ with $0<\H^m(E)<\infty$ is, in many respects, the same as  studying measures $\mu$ with almost everywhere positive and finite upper $m$-density, because $\mu$ and $\H^m\restrict\{x:0<\Theta^{\ast m}(\mu,x)<\infty\}$ are mutually absolutely continuous by Theorem \ref{densmeas}. Also $\mu\ll\H^m$ if and only if $\Theta^{\ast m}(\mu,x)<\infty$ for $\mu$ almost all $x\in\Rn$. If the upper density is infinite, we can have completely different rectifiable measures, in particular many lower dimensional measures are such. One can show fairly easily that all AD-$s$-regular measures with $0<s<m$ are $m$-rectifiable, see \cite[Theorem 4.1]{MM88}. For example, take $0<s<1$ and let $\mu_s$ be the product with itself of a standard $s/2$-dimensional Cantor measure in $\R$. Then $\mu_s$ is 1-rectifiable, but all attempts to approximate with lines clearly must fail. This example also shows that for general measures using Lipschitz maps is a quite different thing than using $C^1$ maps. The measure $\mu_s$ is doubling in the sense that $\mu(B(x,2r))\lesssim\mu(B(x,r)$ for $x\in\spt\mu$ and $r>0$. A more dramatic example was given by Garnett, Killip and  Schul  in \cite{GKS10}. They constructed a 1-rectifiable doubling measure $\mu$ whose support is the whole space $\Rn$. Then $\mu$ is purely unrectifiable with respect to Lipschitz graphs. In these examples the lower density is infinite. One of the few general things one can say about rectifiable measures is that the lower density is positive, see \cite{BS15}, Lemma 2.7.

\begin{thm}\label{BSlowerdens}
If $\mu\in\mathcal M(\Rn)$ is $m$-rectifiable, then $\Theta^m_{\ast}(\mu,x)>0$ for $\mu$ almost all $x\in\Rn$.
\end{thm}

There has been a lot of recent interest in finding criteria for the rectifiability of measures in terms of variants of Jones type square functions. Recall the definition of $\beta_E(x,r)$ from \eqref{beta}. As already mentioned in Section \ref{subtanplanes} $\beta_E(x,r)\to 0$ does not imply rectifiability even for AD-1-regular sets. For measures with finite upper density bilateral approximation together with positive lower density implies rectifiability by Theorem \ref{tanmthm1}, but not without positive lower density because of the example of Preiss in \cite[5.9]{Pre87} of a purely unrectifiable measure with flat tangent measures. But multiscale $\beta$ sums and integrals in the spirit of Theorems \ref{jones} and \ref{betathm} have lead to interesting results for measures.

\subsection{Square functions in general dimensions}\label{sqgendim}

Azzam and Tolsa proved their Theorem \ref{Azzam} for more general measures. In analogy to \eqref{beta2}, we define
\begin{equation}\label{betapmeas}
\beta_{\mu}^{m,2}(x,r)^2=\inf_{V\ \text{affine}\ m-\text{plane}}r^{-m}\int_{B(x,r)}\left(\frac{d(y,V)}{r}\right)^2\,d\mu y.
\end{equation}

\begin{thm}\label{AT1}
Let $\mu\in\mathcal M(\Rn)$. If $0<\Theta^{\ast m}(\mu,x)<\infty$ for  $\mu$ almost all $x\in \Rn$, then $\mu$ is  $m$-rectifiable if and only if  
\begin{equation}\label{ATeq}
\int_0^{1}\beta_{\mu}^{m,2}(x,r)^2r^{-1}dr < \infty
\end{equation}
for  $\mu$ almost all $x\in \Rn$.
\end{thm}

Observe that \eqref{ATeq} alone does not imply rectifiability; it is satisfied by the Lebesgue measure.

Edelen, Naber, and Valtorta improved in \cite[Theorem 2.19]{ENV16}, the sufficient condition for rectifiability:

\begin{thm}\label{ENV}
Let $\mu\in\mathcal M(\Rn)$. If $\Theta^{\ast m}(\mu,x)>0$ and 
$$\int_0^{1}\beta_{\mu}^{m,2}(x,r)^2r^{-1}dr < \infty$$
for  $\mu$ almost all $x\in \Rn$,  then $\mu(\Rn\setminus E)=0$ for some $m$-rectifiable set $E$. If also $\Theta_{\ast}^m(\mu,x)<\infty$ for  $\mu$ almost all $x\in \Rn$,  then $\mu\ll\H^m$, so $\mu$ is $m$-rectifiable.
\end{thm}

This paper contains much more related material with quantitative estimates. As in the proof of Theorem \ref{Azzam} several delicate stopping time  arguments are key tools in the proofs. Tolsa gave in \cite{Tol19} a different proof for the second statement based on \cite{AT15}.

Badger and Naples \cite{BN20} characterized measures which live on countably many Lipschitz graphs in terms of Jones type sums where the cubes are restricted to cones. Dabrowski characterized rectifiability and big pieces of Lipschitz graphs in terms of conical square functions, see \cite{Dab20b}, \cite{Dab21a}. 

But what if we don't make  any density assumptions? Edelen, Naber and Valtorta  proved in \cite[Theorem 2.17]{ENV16} the following Reifenberg-type theorem, which is a very special case of their results:

\begin{thm}\label{ENVreif}
Let $\mu\in\mathcal M(\Rn)$ with $\spt\mu\subset B(0,1)$ and $\eta>0$. Suppose that
$$\int_0^{2}\int\beta_{\mu}^{m,2}(x,r)^2\,d\mu x\,r^{-1}dr \leq \eta^2.$$
Then $\mu=\mu_1+\mu_2$ where $\mu_1(B(0,1)\setminus E)=0$ for some $m$-rectifiable set $E$,\ $\mu_2(\Rn)\leq C(n)\eta$ and $\Theta^m(\mu_2,x)=0$ for $\mu_2$ almost all $x\in\Rn$ . 
\end{thm}

Here $\mu_2$ could, for example, be $\mathcal L^n\restrict A$ for some $A$ with $\mathcal L^n(A)>0$. For Hilbert and Banach space versions, see \cite{ENV19}. The proof is based on corona type decompositions. Naber gives a proof also in \cite{Nab20a}. He formulates it in terms of neck decompositions. Roughly speaking this means that most of $B(0,1)$ is covered with balls $B$ whose neck regions have small measure. The neck region of $B$  is a complement in $B$ of two sets. For $x$ in the first there is a good approximation by planes at all scales. That part is rectifiable. The second set is a union of balls $B(x,r_x)$ in which there is a good approximation by planes at the scales bigger than $r_x$. Of course this is very vague and the interested reader should consult the references above and in the next sentence. Variants of the neck decomposition have been used in several places, see \cite{NV19}, \cite{JN21} and \cite{CJN21}.


\subsection{Square functions and one-dimensional measures}

Badger and Schul introduced new Jones type square functions to study the one-dimensional rectifiability of measures. Define a variant of the quadratic $\beta$ for $\mu\in\mathcal M(\Rn)$ and for cubes $Q$ by

\begin{equation}\label{BSbeta2}
\beta^1_{\mu}(Q)^2=\inf_{L\ \text{a line}}\mu(Q)^{-1}\int_{Q}\left(\frac{d(y,L)}{d(Q)}\right)^{2}\,d\mu y.
\end{equation}

All the cubes $Q$ in this section will be dyadic cubes of side-length at most 1. Then

$$J_{\mu}(x)=\sum_Q\beta^1_{\mu}(3Q)^2\chi_Q(x), x\in\Rn,$$
is essentially the original Jones function for measures. Badger and Schul proved in \cite{BS15} that $J_{\mu}(x)<\infty$ for $\mu$ almost all $x\in\Rn$ if $\mu$ is 1-rectifiable and $\mu\ll\H^1$. The absolute continuity is needed: for the rectifiable measures $\mu_s, 0<s<1,$ (recall the beginning of this chapter) we have $\beta^1_{\mu_s}(3Q)\sim 1$ for the squares $Q$ (with $d(Q)\leq 1$) meeting $\spt\mu_s$, so $J_{\mu}=\infty$ on $\spt\mu_s$. The next definition avoids this situation: 

\begin{equation}\label{BSJ}
\tilde{J}_{\mu}(x)=\sum_Q\beta^1_{\mu}(3Q)^2\frac{d(Q)}{\mu(Q)}\chi_Q(x), x\in\Rn.
\end{equation}

\begin{thm}\label{BSthm1}
Let $\mu\in\mathcal M(\Rn)$. If $\mu$ is 1-rectifiable, then $\tilde J_{\mu}(x)<\infty$ for $\mu$ almost all $x\in\Rn$. The converse holds if $\limsup_{r\to 0}\frac{\mu(B(x,2r))}{\mu(B(x,r))}<\infty$ for $\mu$ almost all $x\in\Rn$. 
\end{thm}

The proof of the first part, in \cite{BS15}, is based on Jones's traveling salesman theorem. The proof of the second part is in \cite{BS17}. 


Martikainen and Orponen showed in \cite{MO18a} that the second part cannot be extended to general measures. They constructed an example where $\tilde J_{\mu}$ is bounded and $\mu$ has zero lower density, so it is not rectifiable.

Consider the modified $\beta$ numbers

\begin{equation}\label{BSbeta1}
\beta_{\mu}^{1\ast}(Q)^2=\inf_{L\ \text{a line}}\max_{R}\min\left(\frac{1}{{d(3R)}},\frac{1}{{\mu(3R)}}\right)\int_{3R}\left(\frac{d(y,L)}{d(3R)}\right)^2\,d\mu y.
\end{equation}
The maximum is taken over the dyadic cubes $R$ for which $d(Q)\leq d(R)\leq 2d(Q)$ and $3R\subset 1600\sqrt{n}Q$. Define
$$J^{\ast}_{\mu}(x)=\sum_Q\beta_{\mu}^{1\ast}(Q)^2\frac{d(Q)}{\mu(Q)}\chi_Q(x), x\in\Rn.$$
Badger and Schul proved in \cite{BS17} the following characterization:

\begin{thm}\label{BSthm2}
Let $\mu\in\mathcal M(\Rn)$ be such that $\Theta^m_{\ast}(\mu,x)>0$ for $\mu$ almost all $x\in\Rn$.  Then $\mu$ is 1-rectifiable if and only if  $J^{\ast}_{\mu}(x)<\infty$ for $\mu$ almost all $x\in\Rn$.
\end{thm}

Combining with Theorem \ref{BSlowerdens} we have 

\begin{thm}\label{BSthm3}
Let $\mu\in\mathcal M(\Rn)$. The $1$-rectifiable and purely $1$-unrectifiable parts in the decomposition $\mu=\mu_r+\mu_u$ are given by
$$\mu_r=\mu\restrict\{x:\Theta^m_{\ast}(\mu,x)>0\ \text{and}\ J^{\ast}_{\mu}(x)<\infty\},$$
$$\mu_u=\mu\restrict\{x:\Theta^m_{\ast}(\mu,x)=0\ \text{or}\ J^{\ast}_{\mu}(x)=\infty\}.$$
\end{thm}



Badger, Li and Zimmerman \cite{BLZ21} proved analogous results in Carnot groups.

In \cite{Ler03} Lerman used different modified Jones functions to get sufficient conditions for 1-rectifiability of measures with quantitative estimates. He did not make any a priori density or absolute continuity assumptions. Naples \cite{Nap20} proved extensions to Hilbert spaces.

\subsection{Square functions and distance of measures}\label{alphas}

Often it is more natural to approximate with Lebesgue measures on planes than with planes. Recall the metric $F_{x,r}$ and the $\a$ coefficients from \eqref{measmetr} and \eqref{alphaur}. Azzam, Tolsa and Toro \cite{ATT20} used a slightly different $\a$s:
$$\tilde{\a}^m_{\mu}(x,r)=\frac{1}{r\mu(B(x,r))}\inf\{F_{x,r}(\mu,c\H^m\restrict V): c\geq 0, V\ \text{an affine}\ m-\text{plane}\}.$$
They proved 

\begin{thm}\label{ATTthm1}
Let $\mu\in\mathcal M(\Rn)$ and suppose that $\limsup_{r\to 0}\frac{\mu(B(x,2r))}{\mu(B(x,r))}<\infty$ for $\mu$ almost all $x\in\Rn$. Then $\mu$ is $m$-rectifiable and $\mu\ll\H^m$ if and only if
$$\int_0^{1}\tilde{\a}_{\mu}^m(x,r)^2r^{-1}dr < \infty$$
for  $\mu$ almost all $x\in \Rn$.
\end{thm}

The 'only if' direction was proved in \cite{Tol15b}. The difference between the $\beta$s and $\a$s is something like what we had before; small $\beta$ tells us locally that most of the measure lives close to a plane but small $\a$ tells us that most of the plane also is close to the support of the measure. So in a way the $\beta$s are smaller than the $\a$s but not in a precise sense. 

Azzam, David and Toro \cite{ADT16} proved rectifiability of doubling measures under different $\a$ assumptions. They do not specify the dimension but define an $\a_{\mu}(x,r)$ minimizing first the distance to normalized $m$-flat measures and then minimizing over $m=0,1,\dots,n$. The finiteness of $\int_0^1\a_{\mu}(x,r)/r\,dr$ for $\mu$ almost all $x\in\Rn$ implies that the doubling measure $\mu$ is a sum over $m=0,1,\dots,n$ of $m$-rectifiable measures. They also proved related quantitative results. In \cite{ADT17} they defined $\a$ numbers measuring self-similarity properties of a measure. If this $\a_{\mu}(x,r)$ is small, then the blow-up at the scale $r$ is close to blow-ups at certain smaller scales, possibly after a rotation.  They proved that then $\mu$ has unique flat tangent measures. From that they got rectifiability, and more, as above.

\section{Rectifiable sets in metric spaces}\label{metricspaces} 

\subsection{Definition and norm}

Let $(X,d)$ be a metric space. There is no problem with the definition:
 
\begin{df}\label{rect-metr}
A set $E\subset X$ is \emph{$m$-rectifiable} if there are Lipschitz maps $f_i:A_i\to X, A_i\subset\R^m, i=1,2,\dots,$ such that
$$\H^m(E\setminus\bigcup_{i=1}^{\infty}f_i(A_i))=0.$$
A set $E\subset X$ is \emph{purely $m$-unrectifiable} if $\H^m(E\cap F)=0$ for every $m$-rectifiable set $F\subset X$.
\end{df}

But everything else is problematic. What properties can we prove? There are no linear subspaces, so can we talk about tangent planes or projections? Any metric space $X$ can be embedded isometrically into a Banach space, and if $X$ is separable, which rectifiable sets are, into $l^{\infty}$. 
Thus we may consider $X$ as a metric subspace of a Banach space $Y$, that is, we can assume that the metric $d$ is given by a norm $\|\cdot\|$. This is often convenient and it gives us linear subspaces. But it does not solve everything. Anyway, Lipschitz maps from subsets of $\R^m$ to Banach spaces can be extended to all of $\R^m$, so in the definition we could consider $f_i:\R^m\to Y$. 

Let us begin with densities which anyway are defined as before.

\subsection{Densities when $m=1$}\label{dens2}
Using \cite[Example 6.4]{AK00a} we first observe that there is no hope, even in Hilbert spaces, to get Besicovitch-Preiss theorem 'existence of density implies rectifiability': Let $X=(0,1)$ with the metric $d(x,y)=\sqrt{|x-y|}$. Then $\H^2_d=(\pi/4)\mathcal L^1$, so $\Theta^2(X,x)=1/2$ for all $x\in X$. However, $X$ is purely 2-unrectifiable by Theorem \ref{densmetr} below. To see the same in Hilbert spaces, consider  $E=\{\chi_{[0,t]}:0<t<1\}\subset L^2([0,1])$. Nevertheless, maybe the weaker 'density 1 implies rectifiability' could be true?

There is no counter-example to Besicovitch's $1/2$-conjecture \ref{1/2} even in general metric spaces. The best  result known is the following theorem of Preiss and Tiser \cite{PT92} improving and extending Besicovitch's $3/4$ Theorem \ref{3/4-dens}:

\begin{thm}\label{densPT}
If $E\subset X$ is $\H^1$ measurable, $\mathcal H^1(E)<\infty$ and 
$$\Theta^1_{\ast}(E,x)> \frac{2+\sqrt{46}}{12}$$ 
for $\H^1$ almost all $x\in E$, then $E$ is 1-rectifiable.
\end{thm}

Notice that $\frac{2+\sqrt{46}}{12}$ is between $\frac{58}{80}$ and $\frac{59}{80}$, so it is less than but  close to $\frac{3}{4}$.

\begin{cor}\label{corPT}
If $E\subset X$ is $\H^1$ measurable and $\mathcal H^1(E)<\infty$, then $E$ is 1-rectifiable if and only if 
$\Theta^1(E,x)=1$ for $\H^1$ almost all $x\in E$.
\end{cor}

That rectifiable sets have density 1 also in metric spaces follows from Kirchheim's theorem which we shall soon discuss.

The proof of  Theorem \ref{densPT} has similar basic ingredients as that of Theorem \ref{3/4-dens}; Besicovitch circle pairs, in a generalized form, are used to find a continuum $C$ with finite measure which intersects $E$ in a set of positive measure. By Theorem \ref{cont} such a $C$ is rectifiable even in metric spaces. This result of Eilenberg and Harrold is perhaps the first result on rectifiability in metric spaces.

To say a bit more, recall from the discussion on the proof of Theorem \ref{3/4-dens} that for any $\a > 0$ we had for some compact subset $F$ of $E$ with $\H^1(F)>0$: if $\sigma=\tfrac{3}{4}+\a$,
\begin{equation}\label{pt1}
\H^1(E\cap U(x,r)) > \sigma 2r\ \text{for}\ x\in F, 0<r<r_0,
\end{equation}
and  
\begin{equation}\label{pt2}
\H^1(E\cap B) \leq (1+\a)d(B)\ \text{whenever}\ E\cap B\not=\emptyset\ \text{and}\ d(B)<r_0,
\end{equation}
then 
\begin{equation}\label{pt3}\H^1(E\cap U(x,y)) \geq \a|x-y|\ \text{for}\ x,y\in F\ \text{with}\ d(x,y)<r_0/3,\end{equation}
where $U(x,y)=U(x,r)\cap U(y,r), r=|x-y|$. For this one  cannot push $\sigma$ below $3/4$. But Preiss and Tiser showed that it is possible to use other sets in place of $U(x,y)$  to reduce $\sigma$ to $\frac{2+\sqrt{46}}{12}$. More precisely (but not quite precisely), they showed that the following condition holds with any $\sigma>\frac{2+\sqrt{46}}{12}$ and some $\tau>0$: whenever \eqref{pt2} holds with small $\a>0$ and $E_1, E_2$ are Borel subsets of $E$ with $d(E_1,E_2)>0$ small satisfying \eqref{pt1} in place of $F$, then there is $U\subset X$ meeting both $E_1$ and $E_2$ such that $\H^1(E\cap U\setminus (E_1\cup E_2)) > \tau d(U).$ Then they showed that this condition together with $\Theta^1_{\ast}(E,x)>\sigma$ (with any $\sigma>0$) implies that some continuum with finite measure intersects $E$ in a set of positive measure. To relate the Preiss-Tiser condition to \eqref{pt3}, observe that $E\cap U(x_1,x_2)\setminus (E_1\cup E_2)=E\cap U(x_1,x_2)$ if $d(E_1,E_2)=|x_1-x_2|, x_i\in E_1$.

\subsection{Densities and area formula for general $m$}

We still consider $X$ as a metric subspace of a Banach space $Y$ and we would like to use Rademacher's theorem. However Lipschitz maps from $\R^m$ to $Y$ need not be differentiable in the standard sense. But Kirchheim found a useful substitute  in \cite{Kir94}. This paper has been very influential in analysis and geometric measure theory in metric spaces.

Kirchheim's idea was to introduce \emph{metric differentials} $MD(f,x)$, which are seminorms on $\R^m$, by
$$MD(f,x)(v)=\lim_{r\to 0}\|f(x+rv)-f(x)\|/r,\ x,v\in \R^m,$$
whenever the limit exists. He showed that if $f$ is Lipschitz, it does exist for almost all $x\in\R^m$, and then for all $y,z\in\R^m$, 
\begin{equation}\label{md}
\|f(z)-f(y)\| -MD(f,x)(z-y) = o(|z-x|+|y-x|).
\end{equation}
So we have something like \eqref{diff} and a substitute for Rademacher's theorem. Then Kirchheim proceeded to prove an area formula. For this define for any seminorm $s$ on $\R^m$ the "Jacobian"
$$J(s) = \a(m)m\left(\int_{S^{m-1}}s(x)^{-m}d\H^{m-1}x\right)^{-1}.$$
Then
\begin{thm}\label{areametr} If $f:\R^m\to X$ is Lipschitz and $A\subset\R^m$ Lebesgue measurable, then
$$\int \card A\cap f^{-1}\{y\}\,d\mathcal H^m y = \int_AJ(MD(f,x))\,d\mathcal L^mx.$$
\end{thm}

The proof is based on the following lemma, \cite[Lemma 4]{Kir94}, going back to Federer's proof of the Euclidean area theorem and \cite[Lemma 3.2.2]{Fed69}:

\begin{lm}\label{lipdecomp} Let $f:\R^m\to X$ be Lipschitz and let $B$ be the set of $x\in\R^m$ for which $MD(f,x)$ exists and is a norm. Then for any $\lambda >1$ there are norms $\|\cdot\|_i$ on $\R^m$ and a Borel partition $(B_i)$ of $B$ such that
$$\|x-y\|_i/\lambda\leq d(f(x),f(y))\leq \lambda\|x-y\|_i\ \text{for}\ x,y\in B_i, i=1,2,\dots.$$
\end{lm}

After this one gets

\begin{thm}\label{densmetr} If  $E\subset X$ is $\H^m$ measurable, $m$-rectifiable and $\H^m(E)<\infty$, then
$\Theta^m(E,x)=1$ for $\H^m$ almost all $x\in E$.
\end{thm}

As we saw, for $m=1$ this is a characterization of rectifiability. Whether it is a characterization when $m>1$ is open; perhaps surprisingly no counter-example is known. But if we replace the Hausdorff measure $\H^m$ by the spherical Hausdorff measure $\mathcal S^m$, then Heisenberg groups give us purely $m$-unrectifiable metric spaces $X$ with $\mathcal S^m(X)<\infty$ and $\Theta^m(\mathcal S^m \restrict X,x)=1$ for $\mathcal S^m$ almost all $x\in X$. We shall come back to this in Chapter \ref{heisenberg}.

\subsection{Tangent planes}

Since we are considering $X$ as a subspace of a Banach space $Y$ we have linear subspaces and we can hope for a tangent plane characterization of rectifiability. Now $Y$ is $l^{\infty}$, or more generally the dual of a separable Banach space. In addition to Kirchheim's metric differentiability result Ambrosio and Kirchheim \cite{AK00a} proved the weak star differentiability for a Lipschitz map $f:\R^m \to Y$: for almost all $x\in\R^m$ there exists a $w^{\ast}$-differential of $f$ at $x$, that is, a linear map $wdf_x:\R^m\to Y$ such that 
\begin{equation}\label{weakstar}
w^{\ast} - \lim_{y\to x}(f(y)-f(x) -wdf_x(y-x))/|y-x|) = 0.
\end{equation}

If  $E\subset X$ is $\H^m$ measurable and $m$-rectifiable, then by Lemma \ref{lipdecomp} we can decompose almost all of it into sets $f_i(B_i)$ where $B_i\subset\R^m$ is a Borel set and $f_i$ is bi-Lipschitz on $B_i$. Then for $\H^m$ almost all $a\in E$ we can define the \emph{weak approximate tangent plane} of $E$ at $a=f_i(x)\in f_i(B_i)$ as 
$$\tanm^{(m)}(E,a) = wdf_{ix}(\R^m).$$

It further follows that if $\pi_a:Y\to \tanm^{(m)}(E,a)$ is a weak star continuous projection ($\pi_a(x)=x$ for $x\in \tanm^{(m)}(E,a)$), then $\|\pi_a(x)-a\|/\|x-a\|\to 1$ as $x\to a, x\in E_a$, where $\Theta^m(E\setminus E_a,a)=0$. 

Ambrosio and Kirchheim also had a converse in \cite[Theorem 6.3]{AK00a}:

\begin{thm}\label{akrect}
Let $E\subset Y$ be $\H^m$ measurable with $\H^m(E)<\infty$. 
Then $E$ is $m$-rectifiable if and only if for $\H^m$ almost all $a\in E$ there is a weak star continuous linear map
$\pi_a:Y\to Y, \dim \pi_a(Y)=m,$ such that for some $s>0$,
$$\lim_{r\to 0}r^{-m}\H^m(\{x\in E\cap B(a,r):\|\pi_a(x-a)\|<s\|x-a\|\})=0.$$
\end{thm}

Ambrosio and Kirchheim had this with positive lower density condition, but using the argument from the Euclidean case one easily sees that this is not needed, see \cite[Lemma 15.14]{Mat95}.


Ambrosio and Kirchheim also proved in \cite{AK00a} an area formula for Lipschitz maps between rectifiable sets and a coarea formula for $\R^k$ valued Lipschitz maps on rectifiable subsets of metric spaces.

\subsection{Cheeger's differentiability spaces and Alberti representations}\label{Cheeger}

As discussed above Kirchheim showed that Lipschitz maps from $\R^m$ to $X$ are differentiable in a sense. In \cite{Che99} Cheeger introduced conditions under which the differentiability of  Lipschitz maps from $X$ to $\R^m$ also makes sense and is true. This was further developed and generalized by Keith \cite{Kei04}.  Suppose that $X$ is equipped with a Borel measure $\mu$. We say 
that $(X,d,\mu)$ is an $m$-dimensional \emph{Lipschitz differentiability space}, LDS, if there is a countable collection of local Lipschitz charts $\phi_i: U_i \to \R^m, U_i\subset X, X=\bigcup_iU_i$, with respect to which every Lipschitz function $f:X\to\R$ is differentiable at $\mu$ almost all $x\in U_i$ in the sense that there is a unique linear function $Df(x) : \R^m\to\R$ such that 
$$f(y) - f(x) = Df(x)(\phi_i(y) - \phi_i(x)) + o(d(x,y)), y\in U_i.$$

It is clear that Euclidean spaces and Riemannian manifolds with the Lebesgue measure are LDS. Not all, even compact, metric measure spaces are LDS, but there are a lot of non-Euclidean ones that are, for example Heisenberg groups. 

Recall the role of Alberti representations in Euclidean spaces from Section \ref{Lebnull}. \emph{Alberti representation} of a measure $\mu$ on $X$ is a Fubini-type decomposition of $\mu$ into $1$-rectifiable measures $\mu_{\gamma}$:
\begin{equation}\label{Alberti3}
\mu(B)=\int\mu_{\gamma}(B)\,dP\gamma,\ B\subset X\ \text{Borel}.\end{equation}
Here $P$ is a probability measure on the space of curve fragments, that is, bi-Lipschitz mappings $\gamma:C_{\gamma}\to X, C_{\gamma}\subset\R$ compact, and the measures $\mu_{\gamma}$ are absolutely continuous with respect to $\H^1\restrict \gamma(C_{\gamma})$. 

For $\phi: X \to \R^m$ the Alberti representations $\gamma_i, i=1,\dots,m$, are said to be  $\phi$-independent if the derivatives of $\phi\circ\gamma_i, i=1,\dots,m$,  are linearly independent and belong to disjoint cones.

Bate  characterised in \cite{Bat15} the Lipschitz differentiability spaces via Alberti representations. The following theorem tells us a great deal of relations between Lipschitz differentiability spaces, Alberti representations and rectifiablity:

\begin{thm}\label{batli}
Suppose that $\H^m(X)<\infty$ and $\Theta_{\ast}^m(X,x)>0$ for $\H^m$ almost all $x\in X$. Then the following conditions are equivalent.
\begin{itemize}
\item[(1)] $X$ is $m$-rectifiable.
\item[(2)] There are Borel sets $U_i\subset X, i=1,2,\dots,$ with $\H^m(X\setminus\bigcup_iU_i)=0$ such that each $(U_i,d,\H^m\restrict U_i)$  is an $m$-dimensional LDS.
\item[(3)] There are Borel sets $U_i\subset X, i=1,2,\dots,$ with $\H^m(X\setminus\bigcup_iU_i)=0$ and Lipschitz functions $\phi_i:U_i\to\R^m$ such that each $\H^m\restrict U_i$ has $m$\ $\phi_i$-independent Alberti representations.
\end{itemize}
\end{thm}

That (1) implies (2) follows from Kirchheim's work \cite{Kir94}. He showed  that in the definition of rectifiable sets Lipschitz maps can be replaced by bi-Lipschitz maps, recall Lemma \ref{lipdecomp}. Then these can be used to find the required charts. The implication from (2) to (3) follows from \cite{Bat15}. Rectifiability from  independent Alberti representations was proved by Bate and Li in \cite{BL17}. They used some of the ideas of David from \cite{Dav15}, who proved related weaker results. A key feature in this technically very complicated argument is to show that the maps $\phi_i$ are bi-Lipschitz on large subsets.

Notice that $m$ in the assumptions and in (2) is the same. This leaves out many interesting LDSs. For example, the Hausdorff dimension of the Heisenberg group $\h^1$ is 4 but it is a 2-dimensional LDS.

David and Kleiner \cite{DK19}  proved a general result for a measure $\mu$ in a metric space without any density conditions: if at $\mu$ almost all points $\mu$ is pointwise doubling, has two independent Alberti representations, and pointed Gromov-Hausdorff limits (see Section \ref{metrictan}) are homeomorphic to $\R^2$, then $\mu$ is 2- rectifiable.

\subsection{Projections as Lipschitz images}\label{Batepro}

First some bad news. The Besicovitch-Federer projection theorem fails in every infinite-dimensional separable Banach space $Y$: Bate, Cs\"ornyei and Wilson \cite{BCW17} constructed a purely unrectifiable set $E\subset Y$ with $\H^1(E)<\infty$ for which $\mathcal L^1(L(E))>0$ for every non-zero continuous linear function $L:Y\to\R$. An earlier weaker result was given by De Pauw in \cite{Dep17}.

However, instead of linear maps Bate \cite{Bat20} considered Lipschitz maps from the metric space $X$ to $\R^n$ and obtained a rectifiability characterization in the spirit of the Besicovitch-Federer projection theorem. There is no natural measure on the space of Lipschitz maps (at least for this purpose), but Baire category gives a notion of typical maps. Related results were proven by Pugh \cite{Pug16} and Galeski \cite{Gal18}.

For the rest of this section we shall assume that $X$ is complete. To state Bate's result let us equip the space of bounded Lipschitz functions $f:X\to\R^n$ such that $\Lip(f)\leq 1$ with the supremum norm to have the complete metric space $\Lip_1(X,n)$. A subset of $\Lip_1(X,n)$ is \emph{residual} if it contains a countable intersection of dense open sets. These are the complements of the 'small' first category sets, countable unions of nowhere dense sets. Bate proved the following

\begin{thm}\label{bate}
Let $0<m\leq n$. If $E\subset X$ is purely $m$-unrectifiable with $\H^m(E)<\infty$ and $\Theta^m_{\ast}(E,x)>0$ for $\H^m$ almost all $x\in E$, then the set of all $f\in \Lip_1(X,n)$ with $\mathcal H^m(f(E))=0$ is residual.

Conversely, if $F\subset X$ is $m$-rectifiable with $\H^m(F)>0$, then the set of all $f\in \Lip_1(X,n)$ with $\mathcal H^m(f(F))>0$ is residual.
\end{thm}

Bate also showed that the positive lower density assumption is not needed if $X=\Rn$ or $m=1$. It may be that it is never needed, but this depends on an announced but unpublished result of Cs\"ornyei and Jones.

The first statement of the theorem is the more essential one. To prove it one needs to approximate well any $f\in \Lip_1(X,n)$ with $g\in \Lip_1(X,n)$ for which $\mathcal L^m(g(E))=0$, or at least arbitrarily small. There are several very interesting ingredients in this argument.

The main tools consist of Alberti representations and weak tangent fields. Recall from Section \ref{Lebnull} that both play essential roles in the investigation of differentiability of Lipschitz maps in \cite{ACP05} and \cite{ACP10}. According to Definition \ref{weaktf} a set $E\subset\Rn$ has a weak $(m-1)$-dimensional tangent field if, generically, the tangents of its 1-rectifiable subsets  span at most $(m-1)$-dimensional subspaces. In the metric space $X$ this spanning condition is interpreted via Lipschitz maps $f:E\to\R^n, 0<m\leq n$. A little more precisely, a Borel function $\tau:E\to G(n,m-1)$ is a weak tangent field with respect to $f$ if for 1-rectifiable sets $\gamma\subset E$  at almost all points $x\in\gamma$ the tangent of $f(\gamma)$ at $f(x)$ lies in $\tau(x)$. 

If $E$ is purely $m$-unrectifiable, then by Theorem \ref{batli} $E$ can have at most $m-1$ independent Alberti representations, which roughly means that one can  move along $E$ from the points of $E$ in at most to $m-1$ directions. This leads to the existence of a weak tangent field. The next step in the proof of Theorem \ref{bate} is to perturb $f$ slightly to a Lipschitz function $g$ with $\mathcal H^m(g(E))$ small. This can be done contracting along directions orthogonal to the planes $\tau(x)$. As these planes are $(m-1)$-dimensional, the small measure is achieved.

As an application of Bate's result David and Le Donne proved in \cite{DL20} that if a subset of a compact metric space has finite $\H^m$ measure, positive lower density and topological dimension $m$, then it is not purely $m$-unrectifiable. In Euclidean spaces this follows, without positive lower density assumption, from the Besicovitch-Federer projection theorem, see \cite{Fed47b}.

There exist easy examples of sets in Euclidean spaces with positive Hausdorff $m$-measure which project to measure zero in \emph{all} $m$-planes. This is true in the Lipschitz setting, too. Vitushkin, Ivanov and Melnikov \cite{VIM63} constructed a purely $1$-unrectifiable subset $E$ of the plane with $\H^1(E)=1$ such that $\mathcal L^1(f(E))=0$ for every Lipschitz function $f:E\to\R$, see \cite{Kel95} for a simplification. 
Modifying an idea of Konyagin, Ambrosio and Kirchheim \cite{AK00a} showed that for any $m>0$, not necessarily an integer, there exists a compact metric space $X$ such that $\H^m(X)=1$ and $\H^m(f(X))=0$ for every Lipschitz map $f$ from $X$ into a Euclidean space. 

Kun, Maleva and  M\'ath\'e \cite{KMM05} proved that an analytic subset $A\subset\Rn$ is purely 1-unrectifiable if and only if for any compact subset $F$ of $A$,\ $\mathcal L^1(f(F))=0$ for every local Lipschitz quotient map $f:F\to\R$. A Lipschitz function $f:F\to\R$ is a local Lipschitz quotient if there is $c>0$ such that $f(F)\cap B(f(x),cr)\subset f(B(x,r))$ for all $x\in F, r>0.$

\subsection{Metric tangents}\label{metrictan}

In Theorem \ref{akrect} we had a characterization of rectifiability in terms of approximate tangent planes, which are linear subspaces of the bigger Banach space. Bate \cite{Bat21} proved corresponding results approximating in the Gromov-Hausdorff distance by $m$-dimensional Banach spaces. Even more, he proved an analogue of Theorem \ref{tanmthm1} in metric spaces. 

First we have to define the relevant concepts. The \emph{Gromov–Hausdorff distance} $d_{GH}(X_1,X_2)$ between metric
spaces $(X_i, d_i), i = 1, 2,$ is the infimum of those $\e > 0$ for which there
exists a metric space $Z$ and isometric embedddings $Z_1$ and $Z_2$ of $X_1$ and $X_2$ into $Z$ such that the Hausdorff distance $d_H(Z_1,Z_2)<\e$. Then $d_{GH}(X_1,X_2)=0$ if and only 
if the completions of $(X_1, d_1)$ and $(X_2, d_2)$ are isometric. If $(X_i,d_i), i=0,1,2,\dots$, are metric spaces and $x_i\in X_i$, we say that the sequence $(X_i,d_i,x_i)$ of pointed metric spaces converges to $(X_0,d_0,x_0)$ if  for any $r>0$ there is a sequence $r_i\to r$ such that $d_{GH}(B_{d_i}(x_i,r_i),B_{d_0}(x_0,r))\to 0$ as 
$i\to\infty$. A \emph{pointed metric measure space} $(X, d, \mu, x)$ consists of a metric space $(X, d)$, a locally finite (bounded sets have finite measure) Borel measure $\mu$ on $X$ and a
distinguished point $x\in\spt\mu$.

Let $(X,d)$ be a complete metric space.  For $K \geq 1$ let $\biLip(K)$ be the set of metric spaces $Y=(\Rn,\rho)$ such that $\rho$ is $K$-bi-Lipschitz equivalent to the Euclidean norm. Denote by $\biLip(K)^{\ast}$ the set of pointed metric measure spaces $(X, d, \mu, x)$ with $(X,d)\in \biLip(K)$.  First a geometric version from \cite{Bat21}:

\begin{thm}\label{tanmthmbate}
Let $E\subset $ be $\H^m$ measurable with $\mathcal H^m(E)<\infty$. Suppose that for $\H^m$ almost all $x\in E,\ \Theta_{\ast}^m(E,x)>0$  and there exists  $K_x \geq 1$ such that for each $r>0$ there exist $Y_r\in$ $\biLip(K_x)$ and $E_r\subset E\cap B(x,r)$ for which 
$$\lim_{r\to 0}r^{-m}\H^m(E\cap B(x,r)\setminus E_r)=0$$
and
$$\lim_{r\to 0}r^{-1}d_{GH}(Y_r\cap B(0,r),E_r)=0.$$
Then $E$ is $m$-rectifiable.
\end{thm}

This result is new also in Euclidean spaces, because to conclude rectifiability it allows much more general approximating sets than planes.

Next we define the tangent measure spaces. For a pointed metric measure space $(X,d,\mu,x)$ let $\tanm(X,d,\mu,x)$ be the set of all pointed metric measure spaces $(Y,\rho,\nu,y)$ such that there exist $r_i>0$ for which $r_i\to 0$ and 
\begin{equation}
(X,r_i^{-1}d,\mu(B(x,r_i))^{-1}\mu,x)\to (Y,\rho,\nu,y).
\end{equation}
Finding a metric which is suitable to define this convergence is a bit of a delicate matter, since for the standard metric measure space convergence there may not be any tangent measure spaces. Standard meaning that one uses the Hausdorff distance for the spaces in $Z$, as in the definition of $d_{GH}$, and a distance metrizing the weak convergence for the push-forwards of the measures. Instead, Bate defined a metric $d_{\ast}$ that only considers the
distance between the measures  and disregards the Hausdorff distance between the embedded metric spaces.

\begin{thm}\label{tanmthmbate1}
Let $E\subset X$ be $\H^m$ measurable with $\mathcal H^m(E)<\infty$ and $\Theta_{\ast}^m(E,x)>0$ for $\H^m$ almost all $x\in E$.  Then the following are equivalent:
\begin{itemize}
\item[(1)] $E$ is $m$-rectifiable.
\item[(2)] For $\H^m$ almost all $x\in E$ there exists an $m$-dimensional Banach space $(\R^m,\|\cdot\|_x)$ such that 
$$\tanm(X,d,\H^m\restrict E,x)=\{(\R^{m},\|\cdot\|_x,\mathcal L^m,0)\}.$$
\item[(3)] For $\H^m$ almost all $x\in E$ there exists  $K_x \geq 1$ such that 
$$\tanm(X,d,\H^m\restrict E,x)\subset \biLip(K_x)^{\ast}.$$
\end{itemize}
\end{thm}

That (1) implies (2) uses Kirchheim's results  in \cite{Kir94} discussed above. The main and hardest arguments involve the proofs of Theorem \ref{tanmthmbate} and the implication (3) $\Rightarrow$ (1) in Theorem \ref{tanmthmbate1}. They use Bate's Lipschitz projection Theorem \ref{bate}. Recall that the proof of the Euclidean Theorem \ref{tanmthm1} also  used projections. Another basic tool needed  consists of approximation of $E$ with continuous images of cubes in $\R^m$. Its formulation is rather complicated, so I only give Bate's own informal description from his introduction: Roughly speaking, it shows that, under the
hypotheses of Theorem \ref{tanmthmbate}, for $\H^m$ almost all  $x\in E$, the following is true: for any $\e>0$
and any sufficiently small $r > 0$, there exists a metric space $\tilde E$ containing $E$ and a continuous
(in fact H\"older) map $\iota:[0,r]^m \to \tilde E$ such that $\iota|\partial [0,r]^m$ is close to having Lipschitz
inverse and $\H^m_{\infty}(\iota ([0,r]^m)\setminus E)<\e r^m$. Here $\H^m_{\infty}$ denotes the $m$-dimensional Hausdorff content. 

To finish the proofs Bate showed that as a consequence of Theorem \ref{bate} $E$ is rectifiable if all of its subsets posses this approximation property. To get some (not quite correct) idea how Theorem \ref{bate} is applied imagine that for some $x\in E, r>0,$ with $\H^m(E\cap B(x,r))\sim r^m$ we would have $\iota([0,r]^m)=E\cap B(x,r)$ for $\iota$ as above. Extending, if possible, $(\iota|\partial [0,r]^m)^{-1}$ to a Lipschitz map $f$ from 
$E\cap B(x,r)$ onto $[0,r]^m$ we would then have $\mathcal{L}^m(f(E\cap B(x,r)))\sim r^m$. By some topology the same would be true for small perturbations of $f$. If $E$ were purely unrectifiable this would contradict Theorem \ref{bate}.

In addition, for the proof of Theorem \ref{tanmthmbate1} Bate extensively developed the properties of the metric $d_{\ast}$ and the metric tangent measure spaces.

\subsection{Menger curvature}
The curvature of a smooth plane curve $C$ at a point $p$ can be defined as the limit of the reciprocals of the largest radii of the circles approaching and only touching $C$ at $p$. The obstruction for not being able to make the radii bigger comes from triples of points on the curve near $p$. So the curvature at $p$ is also the limit of $c(x,y,z)$, where $x,y,z\in C$  go to $p$ and $c(x,y,z)=1/R$ with $R$ equal to the radius of the circle passing through $x,y$ and $z$. If we take three points in a metric space $X$ there is an isometric triple in the plane. Thus we can define the \emph{Menger curvature}  $c(x,y,z)$ for any $x,y,z\in X$. It is not necessary to pass through the plane, because there is a formula which gives  $c(x,y,z)$ in terms of the three distances, see \cite{Men30},\cite{Hah05} or \cite{Hah08}. 

Menger introduced his curvature in \cite{Men30} as a part of extensive studies of geometry of metric spaces. In particular, he characterized in terms of this curvature the metric arcs which are isometric to Euclidean line segments. An interested reader might also wish to have a look at the book \cite{BK70}.

We already discussed this topic in Section \ref{tsp}. In \cite{Hah05}  Hahlomaa proved a variant of Jones's theorem \ref{jones} in general metric spaces. Now $\beta(F)$ is defined in terms of the Menger curvature. In \cite{Hah08} he generalized the David-L\'eger theorem \ref{davleg} to metric spaces:

\begin{thm}\label{Hah}
If $\mathcal H^1(X)<\infty$ and 
$$\int_X\int_X \int_X c(x,y,z)^2\,d\H^1x\,d\H^1y\,d\H^1z < \infty,$$
then $X$ is 1-rectifiable
\end{thm}

The proof follows similar lines as David's original proof in the plane, but much also has to be changed. 


\section{Heisenberg and Carnot groups}\label{heisenberg}
For an introduction to Heisenberg and Carnot groups, see for example \cite{CDPT07} or \cite{BLU07}. Nice surveys related to rectifiability are given by Serapioni \cite{Ser08} and by Serra Cassano \cite{Ser16}. Both Euclidean spaces and Heisenberg groups are special cases of Carnot groups. Except for the Euclidean case they are non-Abelian and have a similar structure as Heisenberg groups but instead of two levels there can be any finite number of levels with different dilation exponents. I do not discuss them explicitly, but I make some comments on them along the way.

\subsection{The Heisenberg group $\h^n$}

Heisenberg group $\hn$ is $\R^{2n+1}$ as a set but with a different metric and non-Abelian group structure. 
We denote the points of $\h^n$ by $p=(z,t)=(x,y,t),z=(x,y)\in\Rn\times\Rn,t\in\R,$ and define the non-Abelian group operation by 
$$p\cdot p'=(z+z',t+t'+\omega(p,p')),$$
where
$$\omega(p,p')=\omega(z,z')=-2\sum_{i=1}^n(x_iy'_{i}-y_{i}x'_i).$$
We shall use the Koranyi metric $d$ given by
\begin{equation}\label{heismetr}
d(p,p')=(|z-z'|^4+(t-t'-\omega(z,z'))^2)^{\frac{1}{4}}.
\end{equation}
Then the ball of radius $r$ at the origin 
$$B(0,r)=\{p:(|z|^4+t^2)^{\frac{1}{4}}\leq r\}$$
is like a cylinder of width $2r$ and height $2r^2$, so $\mathcal L^{2n+1}(B(0,r))\sim r^{2n+2}$. The ball $B(p,r)$ is the image of $B(0,r)$ under the left translation $\tau_p;
\tau_p(q)=p\cdot q.$ The metric is left invariant, that is, $\tau_p$ is an isometry. Moreover, $\mathcal L^{2n+1}=c(n)\mathcal H_d^{2n+2}$ is a Haar measure of the group, so $\mathcal L^{2n+1}(B(p,r))=c'(n)r^{2n+2}$. In particular, we see that the Hausdorff dimension of $\h^n$ is $2n+2$. We also define the dilations $\delta_r, r>0$, by
$$\delta_r(z,t)=(rz,r^2t).$$
Then $d(\delta_r(p),\delta_r(q))=rd(p,q)$. 

In $\h^n$ the metric on  all \emph{horizontal lines} $L_e:=\{te,0):t\in\R\}, e\in S^{2n-1}$, and their left translates is Euclidean, so their subsets are $1$-rectifiable. There are many other 1-rectifiable sets; any two points of $\h^n$ can be joined with a rectifiable curve. This fact also leads to the geodesic metric which is equivalent to the one we have chosen. 

When $n>1$ there are many $m$-rectifiable sets for $m=1,\dots,n$, in the same way. Define the horizontal plane $\h=\{t=0\}$ identified with $\R^{2n}$ and define the space of \emph{horizontal subgroups} 
$$G_h(2n,m)=\{V\in G(2n,m): \omega(p,q)=0\ \text{for all}\ p,q\in V\}.$$
The metric on each $V\in G_h(2n,m)$ is Euclidean, so they are nice $m$-rectifiable sets.
The unitary transformations act transitively on $G_h(2n,m)$ and lead to an invariant measure $\mu_{n,m}$. Not all linear subspaces of $\h$ of dimension $1<m\leq n$ are subgroups, and none of dimensions bigger than $n$.

The vertical subgroups of linear dimension $m-1, 1\leq m-1\leq 2n$, and Hausdorff dimension $m$ are the vertical planes $V\times\R, V\in G(2n,m-2)$. They together with the horizontal subgroups are the only non-trivial homogeneous subgroups of $\hn$, that is, they are closed and invariant under the dilations. But the narrower collection of complementary subgroups will be more relevant for us.

The subgroups $V\in G_h(2n,m)$ and $W=V^{\perp}$  are complementary: $V\cap W = \{0\}$ and they span $\h^n$ in the sense that $V\cdot W = \h^n$. In particular, any $p\in\hn$ has a unique decomposition
\begin{equation}\label{heisproj}
p=p_V\cdot p_W, p_V\in V, p_W\in W.
\end{equation}

Define $\mathcal G(\h^n,m)$ as the set of $V\in G_h(2n,m)$, when $1\leq m\leq n$, and $V\times\R, V\in G(2n,m-2), V^{\perp}\in G_h(2n,2n+2-m)$, when $n+2\leq m\leq 2n+1$. They are the homogeneous subgroups of Hausdorff dimension $m$ which admit a complement in the sense of \eqref{heisproj}.

Also the vertical subgroups in  $\mathcal G(\h^n,m)$ will be $m$-rectifiable, but we have to change the definition, as we shall soon do.  

\subsection{Some analytic tools in Heisenberg and Carnot groups}\label{Heistools}

The analytic structure of $\hn$ is generated by the following vector fields:

$$X_i=\partial_{x_i} + 2y_i\partial_{t}, Y_i= \partial_{y_i} - 2x_i\partial_{t}, i=1,\dots n, \nabla_H = (X_1,\dots,X_n,Y_1,\dots,Y_n).$$

These vector fields at the origin span the horizontal plane $\h$ and at $p$ the plane $\tau_p(\h)$. The tangent vectors of rectifiable curves are spanned by them. We say that a continuous function $u$ on an open subset $U$ of $\h^n$ belongs to $C^1_H$ if $\nabla_Hu$ is continuous. 

There is a very general Rademacher theorem due to Pansu \cite{Pan89}. It says that Lipschitz maps $f$ between Carnot groups are almost everywhere  differentiable. Now the derivative $d_Hf(p)$ is a homogeneous (it commutes with the dilations) homomorphism between the groups such that
$$\lim_{q\to p}\frac{d(f(p)^{-1}\cdot f(q),d_Hf(p)(p^{-1}\cdot q))}{d(p,q)}=0.$$

Franchi, Serapioni and Serra Cassano proved an implicit function theorem for real-valued functions on subsets of $\hn$ and a Whitney extension theorem for Euclidean valued (or values in the horizontal subbundle) functions on subsets of $\hn$, see \cite{FSS01} and \cite{FSS03}. But a Whitney extension theorem for $\hn$ valued functions is missing, and hence also approximation of Lipschitz maps with differentiable maps.

There also are several versions of the area formula, see \cite{JNV20} and \cite{Vit20} and the references given there.

\subsection{Definitions of rectifiability}\label{Heisrectdef}

There are several natural ways to define rectifiable sets in Heisenberg and Carnot groups. Some of them are known to be equivalent in some cases, for some the relations are unknown. We begin with the definition we have already used in metric spaces:
\begin{df}\label{Lrect}
Let $1\leq m\leq n$. A set $E\subset\h^n$ is $m$-rectifiable if there are Lipschitz maps $f_i:A_i\to\h^n, A_i\subset\R^m, i=1,2,\dots,$ such that $\H_d^m(E\setminus\bigcup_{i=1}^{\infty}f_i(A_i))=0$.
\end{df}

The corresponding notion for measures is defined as in \ref{m-rectmeas}. I only gave this definition for low dimensional sets for a good reason: there are no non-trivial $m$-rectifiable subsets of $\h^n$ for $m>n$ by the following result. In $\h^1$ it was proved by Ambrosio and Kirchheim in \cite{AK00a}, a very general statement in Carnot groups is due to Magnani \cite{Mag04}.

\begin{thm}
If $m=n+1,\dots,2n+2$ and $f:A\to\h^n, A\subset\R^m,$ is Lipschitz, then $\H_d^m(f(A))=0$. In particular, $\h^n$ is purely $(2n+2)$-unrectifiable.
\end{thm}

To get an idea why this is true, suppose that $n=1$ and $f=(f_1,f_2)$ maps $\R^3$ into the vertical plane $W=\{y=0\}$. The metric in $W$ is given by the 'norm' $(|x|^4+t^2)^{1/4}$. Then $|f_2(u)-f_2(v)|\lesssim|u-v|^2$, so $f_2$ is constant and $\H_d^3(f(\R^3))=0$ follows. Of course $f$ need not map into a vertical plane, but its Pansu differentials $d_Hf(u)$ must, because they are group homomorphisms and thus $d_Hf(u)(\R^3)$ is an Abelian subgroup. From this one can argue similarly.

Are there any non-trivial $m$-rectifiable subsets when $m>n$? No with our present definition, as we saw above. But there are many with an alternate definition. Recall that Euclidean rectifiable sets can be defined using level sets of regular functions. Based on this we first define \emph{regular surfaces} $S$. If $1\leq m\leq n$ this means that $S$ is locally the image of an open subset of $\R^m$ under an injective continuously differentiable (in Pansu's sense) map with injective derivative. If $n+2\leq m\leq 2n+1$ we say that $S$ is regular if for every $p\in S$ there are an open set $U$ with $p\in U$ and a function $u:U\to \R^{2n+2-m},$ whose coordinate functions belong to $C^1_H$ such that $S\cap U=\{q\in U:u(q)=0\}$ and for $q\in U$ the Pansu derivative $d_Hu(q)$ is surjective.
These surfaces are also regular in the sense that they have tangent subgroups, in the first case via Pansu derivative and in the second via the kernel of $d_Hf(p)$. Notice that when $m\geq n+2$, the topological dimension of $S$ is $m-1$. 

\begin{df}\label{Hrect}
Let $m=1,\dots,n,n+2,\dots,2n+1$. A set $E\subset\h^n$ is $(m,\h)$-rectifiable if
there are $m$-regular surfaces $S_i, i=1,2,\dots,$ such that $\H_d^m(E\setminus\bigcup_{i=1}^{\infty}S_i)=0$.
\end{df} This definition in codimension 1 is due to Franchi, Serapioni and Serra Cassano \cite{FSS01}. They introduced this concept and used it to develop De Giorgi's theory of sets of finite perimeter in Heisenberg groups. We shall return to this in Section \ref{perheis}. For general dimensions and Carnot groups, see \cite{Mag06} and \cite{FSS07}.

As an example, consider the vertical plane $W=\{y=0\}\subset\h^1$.  Then $u(x,y,t)=y$ belongs to $C^1_{H}$ with  $Y_1u\not=0$ on $W$. Hence $W$ is a regular surface and a $(3,\h)$-rectifiable set. As an another example the horizontal plane $\h$ in $\hn$ is not a regular surface, because it has a singularity, a characteristic point, at $0$. But it is regular outside $0$, and hence a $(2n+1,\h)$-rectifiable set. 

In fact, all $C^1$ smooth Euclidean $m$-dimensional, $n+1\leq m\leq 2n$, surfaces are $(m+1,\h)$-rectifiable. They have positive and locally finite 
$\H^{m+1}_d$ measure, they are regular outside the set of characteristic points, which has $\H_d^{m+1}$ measure zero by results of Balogh \cite{Bal03} for $m=2n$ and Magnani \cite{Mag06} for general $m$. A point $p$ is a characteristic point of a hypersurface $S$ if the tangent space (in Heisenberg sense) at $p$ is spanned by the horizontal vector fields $X_i,Y_j$. More generally, the Euclidean $m$-rectifiable sets, $n+1\leq m\leq 2n$, are $(m+1,\h)$-rectifiable since  $\H^{m}(A)=0$ implies $\H_d^{m+1}(A)=0$ for $A\subset\hn$. The converse is false, because $(2n+1,\h)$-rectifiable sets can have Euclidean Hausdorff dimension bigger than $2n$, see \cite{KS04} for an example in $\h^1$ of Hausdorff dimension is $2.5$. General comparisions of Euclidean and Carnot Hausdorff measures can be found in \cite{BTW09}.

When $1\leq m\leq n$, clearly $(m,\h)$-rectifiable sets are $m$-rectifiable, but it is not known if the converse holds. Although Lipschitz maps are almost everywhere differentiable by Pansu's theorem, one would need something like Whitney's extension theorem to go to $C_H^1$ from Lipschitz. This is not known for Heisenberg valued maps.

Intrinsic differentiable graphs and intrinsic Lipschitz graphs have recently been investigated intensively. They provide another definition for rectifiability. In Euclidean spaces cones were used to characterize rectifiability in terms of the approximate tangent planes, and they are directly connected to Lipschitz maps as in the argument preceding Theorem \ref{tanthm1}. In Heisenberg groups the situation is more complicated but we can define a class of Lipschitz maps geometrically in terms of cones. This was done by Franchi, Serapioni and Serra Cassano in \cite{FSS06} in Heisenberg groups and in \cite{FS16} in general Carnot groups. 

For a homogeneous subgroup $G$ define the cone
\begin{equation}\label{heiscone}
\begin{split}
X(p,G,s)&=\{q\in\h^n:d(p^{-1}\cdot q,G)<sd(p,q)\}\\
&=p\{q\in\h^n:d(q,G)<sd(q,0)\}.
\end{split}
\end{equation}
Geometrically these cones look rather different from the Euclidean cones.

Let $V$ and $W$  be complementary subgroups of   $\h^n$; $V\cap W=\{0\}$ and $V\cdot W = \h^n$, with $V$ horizontal and $W$ vertical. For much that follows they could also be complementary homogeneous subgroups of a general Carnot group. We say that $S\subset\hn$ is a (vertical) \emph{intrinsic Lipschitz graph} if there is $s>0$ such that for all $p\in S$,
$$S\cap X(p,V,s)=\{p\}.$$
We say that a function $f:A\to V, A\subset W,$ is (vertical) \emph{intrinsic Lipschitz} if $gr(f):=\{p\cdot f(p):p\in A\}$  is an intrinsic Lipschitz graph. In Euclidean spaces this just means that $f$ is Lipschitz. But now the intrinsic Lipschitz functions need not be Lipschitz in the metric sense, and vice versa, see \cite[Example 3.21]{AS09}.

Changing the roles of $V$ and $W$ we get horizontal intrinsic Lipschitz graphs and functions. Arena and Serapioni proved in \cite[Proposition 3.20]{AS09} that a horizontal function $f$ is intrinsic Lipschitz if and only if $p\mapsto p\cdot f(p)$ is metric Lipschitz.

Intrinsic Lipschitz graphs over the whole planes have positive and locally finite Hausdorff measure (in the appropriate dimension), they are even AD-regular. This holds in general Carnot groups, see \cite[Theorem 3.9]{FS16}.

\begin{df}\label{Graphrect}
Let $m=1,\dots,n,n+2,\dots,2n+1$. A set $E\subset\h^n$ is $(m,\H_{intL})$-rectifiable if
there are $m$-dimensional (in terms of Hausdorff dimension) intrinsic Lipschitz graphs $S_i, i=1,2,\dots,$ such that $\H_d^m(E\setminus\bigcup_{i=1}^{\infty}S_i)=0$.
\end{df} 


There is a slightly more complicated notion of intrinsic differentiable functions and graphs, see \cite{FSS06}, \cite{AS09}. Arena and Serapioni proved in \cite[Theorem 4.2]{AS09} that regular surfaces and intrinsic differentiable graphs are the same locally. So the rectifiable sets defined in terms of intrinsic differentiable graphs are the same as $(m,\h)$-rectifiable sets. To get to intrinsic Lipschitz one would need a Rademacher type theorem for intrinsic Lipschitz functions. Such a theorem was proved by Franchi, Serapioni and Serra Cassano \cite{FSS11} in the codimension 1 case $(m=2n+1)$ in Heisenberg groups. For this they used their results on sets of finite perimeter. Franchi, Marchi and Serapioni \cite{FMS14} extended this to some Carnot groups. In \cite{Vit20} Vittone proved that intrinsic Lipschitz functions $f:A\to V, A\subset W, W$ vertical, are almost everywhere intrinsic differentiable. From this using Whitney type arguments he further showed that intrinsic Lipschitz graphs can be approximated in measure by regular surfaces. This leads to, see \cite[Corollary 7.4]{Vit20}, 

\begin{thm}\label{Vit}
Let $n+2\leq m\leq2n+1$ and let $E \subset\hn$ be $\H_d^m$ measurable with $\H_d^m(E)<\infty$. Then  $E$ is $(m,\h)$-rectifiable if and only if
 $E$ is $(m,\H_{intL})$-rectifiable.
\end{thm}

There are several ingredients of independent interest in Vittone's proof. He introduced an alternative equivalent definition. According to it $S$ is a vertical intrinsic Lipschitz $W$-graph if and only if there exist $\delta>0$ and a Lipschitz map $g:\hn\to W^{\perp}$ such that $g(x)=0$ on $S$ and $g$ satisfies the ellipticity type condition $(g(p\cdot v)-g(p))\cdot v\geq \delta|v|^2$ for $v\in W^{\perp}, p\in\hn$ (the second $\cdot$ is the inner product). In the main part of the argument he used currents (in Heisenberg sense) to show that the blow-ups of $S$ converge to a vertical plane locally uniformly, which is the content of Rademacher's theorem in this setting. 

Counter-examples to Rademacher's theorem for intrinsic Lipschitz graphs in some Carnot groups were given by Julia, Nicolussi Golo and Vittone in \cite{JNV21a}. In \cite{JNV21b} they proved the almost everywhere tangential differentiability of Euclidean valued functions on $C^1_{\h}$ submanifolds of $\hn$ yielding on $\h$-rectifiable sets a Lusin type approximation of Lipschitz functions and a coarea formula.

As said before, there are no Lipschitz maps from Euclidean spaces to parametrize higher dimensional non-trivial sets. But perhaps one could consider Lipschitz maps from other spaces. For instance, in $\h^1$ could 3-dimensional rectifiable sets
 be defined using Lipschitz maps from a fixed vertical plane? Although an intrinsic Lipschitz map need not be metric Lipschitz, maybe some other Lipschitz map could be used to parametrize an intrinsic Lipschitz graph? Such a definition of rectifiability and its consequences, in the generality of Carnot groups, was studied by Pauls \cite{Pau04} and by Cole and Pauls \cite{CP06}. Let $(G_1,d_1)$ and $(G_2,d_2)$ be Carnot groups and $G$ a subgroup of $G_1$ with Hausdorff dimension $m$. Let us say that $E\subset G_2$ is $G$-rectifiable if up to $\H^m_{d_2}$ measure zero it can be covered with countably many Lipschitz images of subsets of $G$.  Except for the cases where $G_1$ is Euclidean, not much is known about the relation of $G$-rectifiability to other concepts of rectifiability, but some partial information exists. 
 
Fix a vertical subgroup $W_n\in\mathcal G(\h^n,2n+1)$. It does not matter which, they all are isometric. Recall that $C^1$ Euclidean hypersurfaces in $\hn$ are $(2n+1,\h)$-rectifiable. Cole and Pauls proved in \cite{CP06} that the $C^1$ hypersurfaces in $\h^1$ are $W_1$-rectifiable, too. Di Donato, F\"assler  and Orponen \cite{DFO19} proved  that the $C^{1,\a}, \a>0,$ hypersurfaces in $\h^n$ are $W_n$-rectifiable, moreover, they have big pieces of bi-Lipschitz images of $W_n$. Earlier the rectifiability was proved by Antonelli and Le Donne \cite{AL20} for $C^{\infty}$ surfaces.

Antonelli and Le Donne showed in \cite{AL20} that there exists a Carnot group containing a $C^{\infty}$ hypersurface without characteristic points which is not  $G$-rectifiable for any Carnot group $G$, but it is rectifiable in the sense of Franchi, Serapioni and Serra Cassano.

The papers \cite{DFO19} and \cite{AL20} show more than rectifiability, instead of Lipschitz maps they use bi-Lipschitz maps. Antonelli and Le Donne discuss more generally definitions of rectifiability based on bi-Lipschitz maps. Bigolin and Vittone \cite{BV10} give a counter-example concerning bi-Lipschitz parametrizations.  Orponen \cite{Orp20} shows that  bi-Lipschitz images of $W_1$ in $\h^1$ admit a corona decomposition by intrinsic bi-Lipschitz graphs. 

In \cite{LY19} Le Donne and Young  proved that a sub-Riemannian manifold which has a Carnot group $G$ as a constant Gromov-Hausdorff limit, see Section \ref{metrictan}, is $G$-rectifiable with bi-Lipschitz parametrizations. The converse also holds by Pansu's Rademacher theorem.

A possibly weaker notion of rectifiability with cones was proposed by Don, Le Donne, Moisala and Vittone in \cite{DLMV19}. They applied it to finite perimeter sets in general Carnot groups.

\subsection{Rectifiable sets and tangent subgroups} 


Next we shall give a characterization of rectifiability in terms of approximate tangent subgroups and tangent measures. They are defined as in the Euclidean case but using the Heisenberg structure. Recall the definition of the cone $X(p,G,s)$ from \eqref{heiscone}.

We say that $G\in\mathcal G(\h^n,m)$  is  an \emph{approximate tangent subgroup} of $E\subset\h^n$ at a point $p\in\hn$ if $\Theta^{\ast m}(E,p)>0$ and for all $0<s<1$,
$$\lim_{r\to 0}r^{-m}\H_d^m(E\cap B(p,r)\setminus X(p,G,s))=0.$$

To define the tangent measures we now set
$$T_{a,r}(p)=\delta_{1/r}(a^{-1}\cdot p),\ p,a\in\hn, r>0.$$
Then, as before, if $\mu$ is a Radon measure on $\hn$, a non-zero Radon measure $\nu$ is called a \emph{tangent measure} of $\mu$ at $a\in\hn$ if there are sequences $(c_i)$ and $(r_i)$ of positive numbers such that $r_i\to 0$ and $c_iT_{a,r_i\#}\mu\to\nu$ weakly. We denote the set of tangent measures of $\mu$ at $a$ by $\tanm(\mu,a)$. The following theorems were proven in \cite{MSS10}:

\begin{thm}\label{MSS}
Let $1\leq m\leq n$ and let $E \subset\hn$ be $\H_d^m$ measurable with $\H_d^m(E)<\infty$. Then the following are equivalent:
\begin{itemize}
\item[(1)] $E$ is $m$-rectifiable.
\item[(2)] $E$ has an approximate tangent subgroup $G_p\in\mathcal G(\h^n,m)$ at $\H_d^m$ almost all $p\in E$.
\item[(3)] For $\H_d^m$ almost all $p\in E$ there is $G_p\in\mathcal G(\h^n,m)$ such that 
$$\tanm(\H_d^m\restrict E,p)=\{c\H_d^m\restrict G_p:0<c<\infty\}.$$
\item[(4)] For $\H_d^m$ almost all $p\in E$\ $\H_d^m\restrict E$ has a unique tangent measure at $p$.
\end{itemize}
\end{thm}

\begin{thm}\label{MSS1}
Let $n+2\leq m\leq 2n+1$ and let $E \subset\hn$ be $\H_d^m$ measurable with $\H_d^m(E)<\infty$. Suppose also that $\Theta_{\ast}^m(E,p)>0$ for $\H_d^m$ almost all $p\in E$. Then the following are equivalent:
\begin{itemize}
\item[(1)] $E$ is $(m,\h)$-rectifiable.
\item[(2)] $E$ has an approximate tangent subgroup $G_p\in\mathcal G(\h^n,m)$ at $\H_d^m$ almost all $p\in E$.
\item[(3)] For $\H_d^m$ almost all $p\in E$ there is $G_p\in\mathcal G(\h^n,m)$ such that 
$$\tanm(\H_d^m\restrict E,p)=\{c\H_d^m\restrict G_p:0<c<\infty\}.$$
\item[(4)] For $\H_d^m$ almost all $p\in E$\ $\H_d^m\restrict E$ has a unique tangent measure at $p$.
\end{itemize}
\end{thm}

That (4) implies (3) follows from Theorem \ref{M2}. There the uniqueness means uniqueness up to multiplication by positive constants. 

It is not known if the positive lower density assumption in Theorem \ref{MSS1} is needed. Also it is not known if we can replace $m$-rectifiable by $(m,\h)$-rectifiable in Theorem \ref{MSS}(1). As mentioned before, the problem here is the lack of Whitney's extension theorem for $\hn$ valued functions.

A few words about the proofs. The case $m\leq n$ is proved by arguments similar to those used in the Euclidean case. But the cones now are harder to deal with if the lower density is 0. This is overcome by proving first that the existence of tangent subgroups implies positive lower density. When $m\geq n+2$ the key for proving that approximate tangents imply rectifiability is the Whitney extension theorem for $\R^k$ valued maps of Franchi, Serapioni and Serra Cassano, \cite{FSS01}. 

Theorem \ref{MSS} was generalized by Idu, Magnani and Maiale \cite{IMM20} to general homogeneous groups, which need not be Carnot groups.

\subsection{Densities}

For low dimensional sets ($m\leq n$) we can apply Kirchheim's theorem \ref{densmetr} to conclude that $m$-rectifiable subsets of $\hn$ with finite measure have density 1 almost everywhere. The converse is not known, except for $m=1$, see Theorem \ref{antmer2}. 

Recently there has been a remarkable progress on densities by Merlo in the two papers \cite{Mer19} and \cite{Mer20}. He proved Preiss's density characterization of rectifiability for codimension 1 subsets of $\hn$:

\begin{thm}\label{mer1}
Let $\mu \in\mathcal M(\hn)$ be such that the positive and finite limit\\ 
$\lim_{r\to 0}r^{-2n-1}\mu(B(p,r))$ exists for $\mu$ almost all $p\in\hn$. Then $\mu$ is  $(2n+1,\h)$-rectifiable.
\end{thm}

As a corollary we have

\begin{thm}\label{mercor}
Let $E \subset\hn$ be $\H_d^{2n+1}$ measurable with $\H_d^{2n+1}(E)<\infty$. Then $E$ is $(2n+1,\h)$-rectifiable if and only if
$\Theta^{2n+1}(E,p)$ exists for $\H_d^{2n+1}$ almost all $p\in E$.
\end{thm}

The fact that the density exists for rectifiable sets follows from \cite[Corollary 3.6]{JNV20}. Julia, Nicolussi Golo and Vittone proved this for sets of general dimensions in arbitrary Carnot groups. They derive the result from an area formula. See also, \cite{AM22}.

Often the density of high dimensional Heisenberg rectifiable sets is strictly less than $1$, but for many metrics the density of the spherical Hausdorff measure is 1 almost everywhere, while for some others it is strictly less than $1$, see \cite{Mag17}.
 
As in the Euclidean case, the proof of Theorem \ref{mer1} splits into two parts. In the first part \cite{Mer19} Merlo proves that the tangent measures are flat:

\begin{thm}\label{mer2}
Let $\mu \in\mathcal M(\hn)$ be such that the positive and finite limit\\ 
$\lim_{r\to 0}r^{-2n-1}\mu(B(p,r))$ exists for $\mu$ almost all $p\in\hn$. Then for $\mu$ almost all $p\in\hn$, 
$$\tanm(\mu,p)\subset \{c\H^{2n+1}\restrict V: V\in \mathcal G(\h^n,2n+1), 0<c<\infty\}.$$ 
\end{thm}

Again an easy argument shows that the tangent measures $\nu$ are uniform; 
$$\nu(B(p,r))=cr^{2n+1}\ \text{for}\ p\in\spt\nu,\ r>0.$$ 
Using moments in the spirit of \cite{Pre87} and \cite{KP87} Merlo proves that their supports are contained in quadratic conical surfaces. Then  disconnectedness is verified saying that vertical flat uniform measures are separated from the others. As in the Euclidean case this implies that typically only vertical flat uniform tangent measures exist. Although many of the arguments run as in the Euclidean case, many essential changes have to be made, in particular, replacing algebraic computations by geometric arguments. I find it surprising that Preiss's proof can be followed at all, since it used heavily the fact that the metric is given by an inner product. Merlo overcomes this by cleverly applying the polarization 
$$V(p,q)= (\|p\|^4+\|q\|^4-\|p^{-1}\cdot q\|^4)/2.$$
It can be written, with a rather complicated formula, in terms of the inner product in $\R^{2n}$ involving for $p=(w,s), q=(z,t)$ all the coordinates $w,z,s,t$ in various combinations.

In the second paper \cite{Mer20} Merlo proves that the flatness of the tangent measures implies rectifiability. This is the case $m=2n+1$ in the following theorem:

\begin{thm}\label{mer3}
Let $m=1,\dots,n$ or $m=2n+1$ and let $\mu \in\mathcal M(\hn)$ be such that for $\mu$ almost all $p\in\hn$,\ $0<\Theta_{\ast}^{m}(\mu,p)\leq\Theta^{\ast m}(\mu,p)<\infty$ and  
$$\tanm(\mu,p)\subset \{c\H^{m}\restrict V: V\in \mathcal G(\hn,m), 0<c<\infty\}.$$
Then  $\mu$ is $(m,\h)$-rectifiable.
\end{thm}

The proof for $m=2n+1$ is very complicated. It is given in general Carnot groups. As for Theorem \ref{tanmthm1} one of the ideas is to use big projections implied, but not easily, by the assumptions. For this some ideas of David and Semmes \cite{DS93} in uniform rectifiability help. But now projections mean the splitting projections $p\mapsto p_W$, recall \eqref{heisproj}, which are much more complicated to handle than the Euclidean projections. They are not even Lipschitz. Big projection on $W$ improves the information given by the assumption: the approximating planes have to be rather close to $W$. Then one can proceed to show that some intrinsic Lipschitz graph intersects $\spt\mu$ in positive measure. 

The proof for low dimensions was given by Antonelli and Merlo in \cite{AM22}. It follows similar patterns, but now the projections are Lipschitz, which helps. On the other hand, serious difficulties are caused by the fact that they prove this in quite general Carnot settings; it is assumed that the tangent subgroups admit at least one normal complementary subgroup.

There is also interesting information about the uniform measures in $\h^1$: they are just the flat measures on horizontal lines, $t$-axis and vertical planes. The first two cases are due to Chousionis, Magnani and Tyson \cite{CMT20}, the third to Merlo \cite{Mer19}, using also some results from \cite{CMT20}. It seems to be open whether there can be non-flat uniform measures in higher dimensional Heisenberg groups.

Combining the case $m=1$ with Theorem \ref{mer3} we get the Besicovitch-Preiss theorem for one-dimensional measures in $\h^1$:

\begin{thm}\label{antmer2}
Let $\mu \in\mathcal M(\h^1)$ be such that the positive and finite limit\\ $\lim_{r\to 0}r^{-1}\mu(B(p,r))$ exists for $\mu$ almost all $p\in\hn$. Then $\mu$ is $1$-rectifiable.
\end{thm}

In Theorems \ref{mer1} and \ref{antmer2} it is important that we use the Koranyi metric $d$. We have been using it all the time, but everywhere else any equivalent metric would be fine. For other, essentially different, metrics the validity of this theorem is unknown. Recall the discussion about the Euclidean $l^{\infty}$ norm  in Section \ref{Densities}.

In \cite{AM20}, \cite{AM21} and \cite{AM22} Antonelli and Merlo base in arbitrary Carnot groups the definition of the rectifiability of a measure on the condition that almost everywhere the tangent measures are flat or that there is a unique flat tangent measure. Flat means a constant multiple of the Hausdorff measure on a homogeneous subgroup that has a complementary subgroup. 
Compare with Theorems \ref{MSS}, \ref{MSS1}, \ref{mer2} and \ref{mer3}. Antonelli and Merlo make a systematic study of such measures with many interesting results generalizing a lot of earlier rectifiability theory. In addition to the analogue of Theorem \ref{mer3}, they prove a characterization with intrinsic Lipschitz and  differentiable graphs, a very general area formula, and rectifiability of level sets of Lipschitz maps.

Let us still have a look at the full group $\hn$ from the point of view of general metric spaces. We know that it has positive and locally finite $\H^{2n+2}$ measure and it is purely unrectifiable. The density of the spherical measure $\Theta^{2n+2}(\mathcal S^{2n+2},p)=1$ but $\Theta^{2n+2}(\H^{2n+2},p)=c(n)<1$ for every $p\in\hn$, see \cite{Mag17}. This difference comes from the fact that balls are not extremals for the isodiametric inequality in $\hn$, see \cite{Rig11}. Recall that it is open in general metric spaces if the Hausdorff density one implies rectifiability. So the Heisenberg group does not give a counter-example. 

\subsection{Projections}
Recall from \eqref{heisproj} that if the subgroups $V\in G_h(2n,m)$ and $W=V^{\perp}$  are complementary, then any $p\in\hn$ has a unique decomposition
$p=p_V\cdot p_W, p_V\in V, p_W\in W$. So we have the projections $P_V(p)=p_V$ and $Q_W(p)=p_W$. Here $V$ is horizontal and $P_V$ is just the standard orthogonal projection onto $V$ but $Q_W$ is more awkward, in particular, it is not Lipschitz. The effect of these projections on the Hausdorff dimension (Marstrand type theorems) was studied in \cite{BFMT12}.

Very little is known about the relations between projections and rectifiability.  There is a bit of that in the proof of Theorem \ref{mer3}, as briefly explained above. As mentioned in Section \ref{projections} Hovila, E. and M. J\"arvenp\"a\"a and Ledrappier \cite{HJJL14} proved the Besicovitch-Federer projection theorem for transversal families of linear maps $\Rn\to\R^m$. In \cite{Hov14} Hovila proved that the family $\{P_V:V\in G_h(2n,m)\}$ is transversal. Thus the Besicovitch-Federer projection theorem holds in $\R^{2n}$ for $G_h(2n,m)$, which is strictly lower dimensional than the whole $G(2n,m)$ when $m>1$. This immediately leads to the following result in $\hn$. There $\pi(z,t)=z$. 

\begin{thm}
Let $1\leq m\leq n$ and let $E \subset\hn$ be a Borel set with $\H^m(\pi(E))<\infty$. Then 
$\H^m(P_V(E))=0$ for $\mu_{n,m}$ almost all $V \in G_h(2n, m)$ if and only if $\pi(E)$ is purely $m$-unrectifiable.
\end{thm}

\subsection{Uniform rectifiability}\label{UR}
The above results on rectifiability could give a starting point for uniform rectifiability. So far it is not nearly fully developed, but there are many interesting results, however, as far as I know, only in the cases of dimension 1 and codimension 1. Some of the terminology here is analogous to that in Chapter \ref{unifrect}.

For one-dimensional sets the question again is about traveling salesman type results. Define $\beta_E(x,r)$ as in \eqref{beta} but the infimum is taken only over the horizontal lines, that is, the left translates of the of the lines through 0 in the horizontal plane. 
\begin{thm}\label{tsthmheis}
Let $G$ be a step two Carnot group of Hausdorff dimension $Q$. If $\Gamma\subset G$ is a rectifiable curve, then
$$\int_{G}\int_0^{\infty}r^{-Q}\beta_{\Gamma}(x,r)^4\,dr\,dx\lesssim \H^1(\Gamma).$$
\end{thm}

\begin{thm}\label{tsthmheis1}
If $p<4$ and $E\subset\h^1$ is compact and such that 
$$\beta_p(E):=\int_{\h^1}\int_0^{\infty}r^{-4}\beta_E(x,r)^p\,dr\,dx<\infty,$$
then $E$ is contained in a rectifiable curve $\Gamma$ for which $\H^1(\Gamma)\lesssim d(E)+ \beta_p(E).$
\end{thm}

Theorem \ref{tsthmheis} was first proved in $\h^1$ by Li and Schul \cite{LS16a} and then extended by  Chousionis, Li and Zimmerman \cite{CLZ19a}, with a modified satement, to arbitrary Carnot groups. Theorem \ref{tsthmheis1} was proved by Li and Schul in \cite{LS16b}. The exponent $4$ of $\beta$ comes from the geometry of $\h^1$. 

It is not known if Theorem \ref{tsthmheis1} holds with $p=4$. However, Li \cite{Li19} produced a Carnot group where there is a gap for the exponents. He then defined modified $\beta$ numbers, based on the stratification of the group, and proved that they give a necessary and sufficient condition in arbitrary Carnot groups; the analogues of Theorems \ref{tsthmheis} and \ref{tsthmheis1} hold with the exponent $2s$, where $s$ is the step of the group.

As we have seen, intrinsic Lipschitz graphs play an important role in low codimensional rectifiability, so they probably should be basic examples of uniformly rectifiable sets. Here are some results in this direction for codimension 1 subsets of $\hn$. The $\beta$ numbers are defined as before but restricting the approximation to vertical planes.  

Chousionis, F\"assler and Orponen proved in \cite{CFO19a} that an AD-3-regular set in $\h^1$ has big pieces of intrinsic Lipschitz graphs if and only if it  satisfies the weak geometric lemma and has big vertical projections. The projections here are the group projections $Q_W$ onto vertical planes defined above. 


Let $E\subset\hn$ be Lebesgue measurable such that $E$ and $\hn\setminus E$ are AD-$(2n+2)$-regular and $\partial E$ is AD-$(2n+1)$-regular. Then Naor and Young \cite{NY18} established a corona decomposition with intrinsic Lipschitz graphs for $\partial E$. This lead to an isoperimetric-type inequality, which was used to solve a fundamental combinatorial problem. See also \cite{NY20}. By \cite[Section 8]{FOR20} the set $\partial E$ as above is a particular case of a Semmes surface, that is, a closed set $F$ such that for every $p\in F$ and $0<r<r_0$ there are balls of radius $cr$ in different components of $B(p,r)\setminus F$. Partially using techniques of \cite{NY18} F\"assler, Orponen and Rigot proved in \cite{FOR20} that Semmes surfaces have big pieces of intrinsic Lipschitz graphs. David \cite{Dav88} had earlier proved the corresponding Euclidean result for which this paper gives a new proof. 

Di Donato, F\"assler  and Orponen proved in \cite{DFO19}  that $C^{1,\a}, \a>0$, codimension 1 intrinsic graphs have big pieces of bi-Lipschitz images of vertical hyperplanes in $\hn$. 

Chousionis, Li and Young proved in \cite{CLY20} that intrinsic Lipschitz graphs in $\hn$ satisfy the geometric lemma when $n\geq 2$. They used a Dorronsoro inequality of F\"assler and Orponen \cite{FO20}  and reduction to lower dimensional groups by slicing. This slicing technique does not work when $n=1$ and the problem is open.

We shall return to Heisenberg groups in connection of singular integrals.

\section{Bounded analytic functions and the Cauchy transform}\label{Analytic}

General references for this chapter are \cite{Tol14}, \cite{Paj02}, \cite{Chr90a} \cite{Dud10}, \cite{Mat95} and \cite{Gar72}. Verdera has a recent survey in \cite{Ver21}.  Many of the topics are continued in the higher dimensional case in the next chapter.

\subsection{Removable sets and Menger curvature}
In 1888 Painlev\'e \cite{Pai1888} proved that any compact subset $E$ of the plane with one-dimensional Hausdorff measure zero is removable for bounded analytic functions.  This was before the existence of Hausdorff measures, but the condition simply means that $E$ can be covered with finitely many discs the sum of whose  diameters is arbitrarily small. 

The removability means that any bounded complex analytic function in $U\setminus E$, where $U$ is an open set containing $E$, has an analytic extension to $U$.  It is easy to see by the Cauchy integral formula that it is enough to consider $U=\C$ and then by Liouville's theorem $E$ is removable if and only if every bounded analytic function in the complement of $E$ must be constant. The problem of finding a geometric characterization of removability is called Painlev\'e's problem.

In 1947 Ahlfors \cite{Ahl47} characterized removable sets in terms of the \emph{analytic capacity} $\gamma$:
$$\gamma(E)=\sup\{\lim_{|z|\to\infty}|z(f(z)-f(\infty))|: |f|\leq 1, f\ \text{analytic in}\ \C\setminus E\}.$$
Then $\gamma(E)=0$ if and only $E$ is removable. However, this is still a complex analytic characterization.

Here is an easy proof of Painlev\'e's theorem: Let $f:\C\setminus E\to\C$ be analytic with $|f|\leq 1, f(\infty)=0$, and let $z\in\C\setminus E$. Let $\e>0$ and cover $E$ with discs $B_j, j=1,\dots,k,$ such that $\sum_{j=1}^kd(B_j)<\e$ and $z$ is outside them. Let $\Gamma$ be the boundary of their union. Then by the Cauchy integral formula
$$f(z)= - \frac{1}{2\pi i}\int_{\Gamma}\frac{f(\zeta)}{\zeta-z}\,d\H^1\zeta.$$
For small $\e$ this is bounded by $C(z)\e$, so $f(z)=0$.

What can we say if we only have $\H^1(E)<\infty$? From the above argument we still get that 
$$f(z)= - \int_{\Gamma_j}\frac{f(\zeta)}{\zeta-z}\,d\H^1\zeta,$$
where the $\Gamma_j$ are surrounding $E$ with $\H^1(\Gamma_j)<\pi\H^1(E)+1$ and $d(\Gamma_j,E)\to 0$. Some subsequence of $f\H^1\restrict \Gamma_j$ converges weakly to a complex measure on $E$. Looking at this a bit more closely one finds that this measure is absolutely continuous with respect to $\H^1\restrict E$. This leads to the representation
$$f(z)= C_E(\varphi)(z):=\int_E\frac{\varphi(\zeta)}{\zeta-z}\,d\H^1\zeta,\ z\in \C\setminus E,$$
where $\varphi$ is a bounded complex valued Borel function on $E$ and $C_E(\varphi)$ is the \emph{Cauchy transform} of $\varphi$.

So we can rephrase Painlev\'e's problem, at least for sets with finite length, in terms of the Cauchy transform: when can we put a bounded non-trivial Borel function on $E$ whose Cauchy transform is bounded? How does this relate to rectifiability? For $\H^1$ almost all points $z\in E$ with $\varphi(z)\neq 0$ the integral $\int_E|\zeta-z|^{-1}|\varphi(z)|\,d\H^1\zeta = \infty$, so $C_E(\varphi)(z)$ is not defined on $E$. Thus $C_E(\varphi)(z)$ can be bounded when $z$ is near $E$ only due to cancellation and the cancellation comes from the symmetries of $E$: for many points $z\in E$, if $\zeta\in E$ the symmetric point $2z-\zeta$ should be close to $E$. Heuristically, rectifiable sets have such symmetries and purely unrectifiable fail to have them. We have a theorem of David \cite{Dav98}:

\begin{thm}\label{David}
Let $E\subset \C$ be compact with $\H^1(E)<\infty$. Then $E$ is removable for bounded analytic functions if and only if $E$ is purely 1-unrectifiable.
\end{thm}

Although the above discussion may give some indication why this could be true, completing the proof has been a long and difficult journey. It was finished when David in \cite{Dav98} showed that purely unrectifiable sets are removable. But let us first look at the other direction. So we should show that if $E$ is rectifiable with $0<\H^1(E)<\infty$ then there is a bounded non-constant analytic function in its complement. We may assume that $E$ is a subset of a $C^1$ graph $\Gamma$, even with a small Lipschitz constant. If $E$ is a line segment, finding $\varphi$ is easy: if $\varphi$ is a smooth function on $E$ vanishing at the end-points, then $C_E(\varphi)$ is bounded by a direct computation. If $E$ is a subset of the graph $\Gamma$ of a $C^{1,\a}$ function for some $\a>0$, then one can still construct $\varphi$, see \cite[Theorem I.7.1]{Gar72}. But nobody knows how to construct it if $\Gamma$ is only $C^1$. However, it is known to exist by duality methods involving the Hahn-Banach theorem or some of its equivalent forms, see \cite{Chr90a}, \cite{Paj02} or \cite{Tol14}. To do this one considers $C_{\Gamma}$ as a singular integral operator on $\Gamma$. We shall discuss this and other singular integrals more in the next chapter. The final step needed was Calder\'on's theorem, \cite{Cal77}, saying $C_{\Gamma}$ is bounded in $L^2(\Gamma)$. More precisely,
$$\int_{\Gamma}|C_{\Gamma,\e}(g)|^2\,d\H^1\lesssim \int_{\Gamma}|g|^2\,d\H^1\ \text{for}\ g\in L^2(\Gamma),$$
with constant independent of $\e$, where
\begin{equation}\label{Cauchymax}
C_{\Gamma,\e}(g)(z)=\int_{\{\zeta\in\Gamma:|\zeta-z|>\e\}}\frac{g(\zeta)}{\zeta-z}\,d\H^1\zeta,\ z\in \Gamma.
\end{equation}
This only gives us non-zero functions $g$ such that $C_{\Gamma,\e}(g)\in L^2(\Gamma)$ uniformly in $\e$, but then duality methods based on the Hahn-Banach theorem can be used to produce $\varphi$ for which $C_{\Gamma}(\varphi)$ is bounded in $\C\setminus E$.

Calder\'on \cite{Cal77} proved the $L^2$ boundedness of the Cauchy transform on Lipschitz graphs with small Lipschitz constant in 1977. In 1982 Coifman, McIntosh and Meyer \cite{CMM82} proved this for general Lipschitz graphs. Later many people have given different proofs. David \cite{Dav84} proved that the Cauchy transform is bounded on a rectifiable curve $\Gamma$ if and only if  $\Gamma$ is AD-1-regular, that is, uniformly 1-rectifiable.

To prove the converse statement of Theorem \ref{David} we need to show that if $E$ is not removable, then it contains a rectifiable subset of positive measure. We can again start with a bounded complex valued Borel function $\varphi$ on $E$ such that $C_E(\varphi)$ is bounded in $\C\setminus E$. There are three main steps in the proof:

(1) Modify $\H^1\restrict E$ and $\varphi$ to a finite Borel measure $\mu$ and a bounded Borel function $g$ such that $\mu\sim\H^1\restrict E$ on some $F\subset E$ with $\mu(F)>0, \mu(B(z,r))\leq r$ for $z\in \C, r>0$, the real part of $g \geq\delta >0$ on $\C$ and $C_{\mu}g\in BMO(\mu)$.

This was done in \cite{DM00} except that we only derived $L^2$ estimates for $C_{\mu}g$, the BMO estimates were then (later, although the papers appeared in different order) proved in \cite{Dav98}. For the modifications we needed to construct generalized dyadic cubes for non-doubling measures similar to those mentioned in Section \ref{Basic tools}. Then $BMO(\mu)$ refers to BMO defined with these cubes. The construction of $\mu$ and $g$ is done by stopping time arguments similar to those used by Christ in \cite{Chr90b} in the AD-regular case.   

Before explaining the step (2) let us define that $C_{\mu}$ is bounded in $L^2(\mu)$ if the truncated \emph{Cauchy transforms} $C_{\mu,\e}$, 
$$C_{\mu,\e}g(z)=\int_{|\zeta-z|>\e}\frac{g(\zeta)}{\zeta-z}\,d\mu\zeta, \e>0,$$ 
are uniformly bounded in $L^2(\mu)$. Furthermore, we say that $C_{E}$ is bounded in $L^2(E)$ if this holds for $\mu=\H^1\restrict E$.

(2) Use a $T(b)$-theorem to prove that the singular integral operator $C_{\mu}$ is bounded in $L^2(\mu)$. Again, a suitable $T(b)$-theorem in the AD-regular case was proved in \cite{Chr90b}. But in the present general non-doubling situation no such result was known and David proved it in \cite{Dav98}. Roughly speaking it says that the existence of $b=g$ as in step (1) implies the $L^2$-boundedness of $C_{\mu}$.

(3) Recall from \eqref{menger} the Menger curvature $c(z_1,z_2,z_3)=1/R$, where $R$ is the radius of the circle passing through $z_1,z_2$ and $z_3$. The final step (which was completed first) is based on its relation to the Cauchy kernel $1/z$, which is provided by the following remarkable formula of Melnikov \cite{Mel95}:

\begin{equation}\label{Cauchymenger}
\sum_{\sigma}\frac{1}{(z_{\sigma(1)}-z_{\sigma(3)})\overline{(z_{\sigma(2)}-z_{\sigma(3)})}}=c(z_1,z_2,z_3)^2
\end{equation}
for distinct points $z_1,z_2,z_3\in\C$, where $\sigma$ runs through the six permutations of $\{1,2,3\}$. See \cite[Section 3.2]{Tol14}, for the easy elementary proof. Define the curvature of $\mu$ by
$$c^2(\mu)=\iiint c(z_1,z_2,z_3)^2\,d\mu z_1\,d\mu z_2\,d\mu z_3.$$
Melnikov and Verdera \cite{MV96} integrated \eqref{Cauchymenger} with respect to $\mu$, wrote $|C_{\mu,\e}1(z)|^2$ as a double integral and used Fubini's theorem six times to get
\begin{equation}\label{Melnikovcurv}
\int |C_{\mu,\e}1|^2\,d\mu=\frac{1}{6}c^2_{\e}(\mu) + O(\mu(\C)),
\end{equation}
where in $c^2_{\e}(\mu)$ the integration is performed over the triples with mutual distances bigger than $\e$. Since the left hand side is bounded by the step (2), so is the right hand side, whence $c^2(\mu)<\infty$, and the proof is completed by the David-Leg\'er theorem \ref{davleg}.

The formula \eqref{Cauchymenger} is remarkable firstly because it relates to the Cauchy kernel a non-negative quantity which vanishes on lines, and only on lines. This alone would be very useful. Secondly this quantity has a concrete geometric meaning.

The above argument gives

\begin{thm}\label{CauchyL2}
Let  $E\subset\C$ be $\H^1$ measurable with $\H^1(E)<\infty$. If $C_{E}$ is bounded in $L^2(E)$, then $E$ is $1$-rectifiable. 
\end{thm} 

Nazarov, Treil and Volberg, see \cite{NV99} or \cite{NTV02}, proved a little later than David a $T(b)$ theorem which also completes the proof of Theorem \ref{David}. Their method was quite different. They used random translations of the standard dyadic squares finding them in a good position with big probability. This method has turned out to be very useful in many later developments. 

The formula \eqref{Melnikovcurv} was used by Melnikov and Verdera in \cite{MV96} to give yet another proof for the $L^2$-boundedness of the Cauchy transform on Lipschitz graphs.

If in the above sketch of the proof we can start with a real-valued $\varphi$, then the modified function $g$ will be positive and we can immediately go from step (1) to step (3). This was done in \cite{DM00} and it resulted to the analogue of Theorem \ref{David} for Lipschitz harmonic functions. We shall discuss them in higher dimensions in the next chapter.

It is easy to see that sets of Hausdorff dimension bigger than 1 are not removable. Hence after David's theorem only the case with dimension 1 and infinite measure remained open. The complete solution of Painlev\'e's problem was given by Tolsa in \cite{Tol03}:

\begin{thm}\label{Tolsa}
Let  $E\subset\C$ be compact. Then $E$ is not removable if and only there is $\mu\in\mathcal M(E)$ such that $\mu(B(z,r))\leq r$ for $z\in \C$ and $r>0$ and $c^2(\mu)<\infty$. 
\end{thm}

Notice that when $\H^1(E)<\infty$ Tolsa's criterion is equivalent to rectifiablity by Theorems \ref{davleg} and \ref{mengerthm}, the latter applied to Lipschitz graphs.

By results of Melnikov \cite{Mel95} and Tolsa \cite{Tol05} we also have a quantitative version of Theorem \ref{Tolsa} in terms of the analytic capacity, see also Theorems 4.14 and 6.1 in \cite{Tol14}.

\begin{thm}\label{Tol}
Let  $E\subset\C$ be compact. Then 
$$\gamma(E)\sim \sup\{\mu(E): \mu\in\mathcal M(E), c^2(\mu)\leq\mu(E), \mu(B(z,r))\leq r\ \text{for}\ z\in\C, r>0\}.$$  
\end{thm}

For AD-regular sets we have

\begin{thm}\label{MMV}
If  $E\subset\R^2$ is closed and AD-$1$-regular, then $E$ is uniformly $1$-rectifiable if and only $C_{E}$ is bounded in $L^2(E)$.
\end{thm}

As mentioned above the boundedness on uniformly rectifiable sets is due to David \cite{Dav84}. The converse was proved in \cite{MMV96} using \eqref{Melnikovcurv}, the above mentioned results of Christ and Theorem \ref{betamthm}, with $p=2$, of David and Semmes. It also gave Theorem \ref{David} for AD-regular sets. 

Jaye and Nazarov gave in \cite{JN14} a different proof for Theorem \ref{MMV} using reflectionless measures, see some comments on them at the end of Section \ref{L2UR}. Their proof does not use Menger curvature but relies on other special properties of the Cauchy kernel.

For general measures $\mu$ on $\C$ with linear growth; $\mu(B(z,r))\leq r$ for $z\in \C, r>0$, the boundedness of $C_{\mu}$ in $L^2(\mu)$ is equivalent to $c^2(\mu\restrict B)\lesssim\mu(2B)$ for all discs $B$, see \cite[Theorem 3.5]{Tol14}. Hence, recalling Theorem \ref{beta1thm}, the following result of Azzam and Tolsa \cite{AT15} extends Theorem \ref{MMV}:

\begin{thm}\label{AT15}
Let $\mu\in\mathcal M(\C)$ with $\mu(B(z,r))\leq r$ for $z\in \C, r>0$. Then  
\begin{equation}\label{AT15eq}
c^2(\mu)+\mu(\C)\sim\int_0^{\infty}\int\beta^{1,2}_{\mu}(x,r)^2\frac{\mu(B(x,r))}{r}\,d\mu x\frac{1}{r}\,dr+\mu(\C).
\end{equation}
\end{thm}

Recall their closely related rectifiability characterization in Theorem \ref{AT1}. For the corresponding result for the Riesz kernels $R_{n-1}$ in $\Rn$, see Theorem \ref{DT21}.

Combining Theorem \ref{AT15} with Theorems \ref{Tolsa} and \ref{Tol} we have a $\beta$ characterization of removable singularities of bounded analytic functions and of analytic capacity.

\subsection{Projections}
The validity of Theorem \ref{David} was conjectured by Vitushkin in \cite{Vit67}. Actually he formulated his question for general compact sets $E$: is $E$ removable if and only if $\H^1(P_L(E))=0$ for almost all lines $L$ through the origin? By the Besicovitch projection theorem  \ref{proj1} this is equivalent to pure unrectifiability when $E$ has finite measure. In general the answer is no. I showed in \cite{Mat86b} that the projection condition is not conformally invariant whereas the removability obviously is. This did not tell us which of the implications is false, but soon afterwards Jones and Murai constructed in \cite{JM88} a non-removable set with zero projections. An easier construction with Menger curvature was done by Joyce and M\"orters in \cite{JM00}. The other direction is still open. Because of Tolsa's Theorem \ref{Tolsa} this can be stated purely as a geometric measure theory problem: if $\H^1(P_L(E))>0$ for positively many lines $L$, is it then possible to construct $\mu\in\mathcal M(E)$ such that $\mu(B(z,r))\leq r$ for $z\in \C$ and $r>0$ and $c^2(\mu)<\infty$?  Chang and Tolsa proved in \cite{CT20} a partial result in this direction.

Dabrowski and Villa \cite{DV22} proved the quantitative estimate $\gamma(E)\sim d(E)$ provided the compact set $E\subset\C$ satisfies the projection condition \eqref{BP}. For this they used Theorems \ref{Tol} and \ref{AT15}.

\subsection{Principal values}\label{PV}
From Theorem \ref{CauchyL2} we get a characterization of rectifiability in terms of the boundedness of the Cauchy transform. Now we give a characterization in terms of the convergence properties of the Cauchy transform.

\begin{thm}\label{Cauchypv}
Let $E\subset\C$ be $\H^1$ measurable with $\H^1(E)<\infty$. Then $E$ is 1-rectifiable if and only the finite limit
$$\lim_{\e\to 0}\int_{\{\zeta\in E:|\zeta-z|>\e\}}\frac{1}{\zeta-z}\,d\H^1\zeta$$
exists for $\H^1$ almost all $z\in E$.
\end{thm}

The convergence for rectifiable sets was proved in \cite{MM94}, in \cite{Ver92} Verdera gave a different Hahn-Banach proof. The other, more difficult, direction was proved by Tolsa in \cite{Tol00}. Under the assumption of positive lower density it was proved in \cite{Mat95a}. Both papers contain more general results for measures. We say that $\mu\in\mathcal M(\C)$ has a \emph{principal value} at $z$ if the finite limit
$$C\mu(z):=\lim_{\e\to 0}\int_{\{\zeta:|\zeta-z|>\e\}}\frac{1}{\zeta-z}\,d\mu\zeta$$
exists. Define the maximal transform
$$C^{\ast}\mu(z):=\sup_{\e>0}\left|\int_{\{\zeta:|\zeta-z|>\e\}}\frac{1}{\zeta-z}\,d\mu\zeta\right|.$$

Proving and using a variant of a $T(b)$ theorem of Nazarov, Treil and Volberg, see \cite{NV99} or \cite{NTV02}, Tolsa proved first that if $\Theta^{\ast 1}(\mu,z)<\infty$ and $C^{\ast}\mu(z)<\infty$ for $\mu$ almost all $z\in\C$, then $C_{\mu}$ is $L^2$ bounded on a set of large $\mu$ measure. Then he concluded that by \eqref{Melnikovcurv} there are compact sets $E_i$ such that $\mu(\C\setminus\cup_{j=1}^{\infty}E_j)=0$ and $c^2(\mu\restrict E_j)<\infty.$ Applying this and Theorem \ref{davleg} to $\mu=\H^1\restrict E$, Theorem \ref{Cauchypv} follows, and even with the existence of principal values replaced by the finiteness of the maximal function. 

In \cite{Mat95a} I proved

\begin{thm}\label{pvmu}
Let $\mu\in\mathcal M(\C)$. If $\Theta^1_{\ast}(\mu,z)>0$ and $C\mu(z)$ exists for $\mu$ almost all $z\in\C$, then $\mu$ is 1-rectifiable.
\end{thm}

The proof uses tangent measures, recall Section \ref{tanmeas}. The assumptions imply that for $\mu$ almost all $z\in\C$ every $\nu\in\tanm(\mu,z)$ is \emph{symmetric} which means that
$$\int_{B(z,r)}(\zeta-z)\,d\nu\zeta = 0\ \text{for}\ z\in\spt\nu, r>0.$$
Then up to discrete measures the symmetric measures were characterized. Every non-discrete symmetric measure is either a constant multiple of the two-dimensional Lebesgue measure, or a constant multiple of the one-dimensional Lebesgue measure on a line, or a countable sum of constant multiples of the one-dimensional Lebesgue measures on parallel lines. Theorem \ref{pvmu} was deduced from this.

Tolsa extended also Theorem \ref{pvmu} in \cite{Tol07}:

\begin{thm}\label{Tolsapv}
Let $\mu\in\mathcal M(\C)$. If $\Theta^{\ast 1}(\mu,z)>0$ and $C^{\ast}\mu(z)<\infty$ for $\mu$ almost all $z\in\C$, then $\mu$ is 1-rectifiable.
\end{thm}

The proof involves similar ingredients as \cite{Tol00} but considerable extra difficulties are caused by the weaker density assumptions. The paper contains more detailed information about measures with finite Cauchy maximal transform.

For $\mu\in\mathcal M(\Rn), x\in\Rn$ and $r>0$, define 
$$C^m_{\mu}(x,r)=r^{-m-1}\int_{B(x,r)}(y-x)\,d\mu y.$$ 
It vanishes on $\spt\mu$ for all $r>0$ if and only if $\mu$ is symmetric. The finiteness in the following theorem is a kind of approximate symmetry condition.

\begin{thm}\label{Villa}
Let $\mu\in\mathcal M(\Rn)$ with $0<\Theta^{\ast m}(\mu,x)<\infty$ for $\mu$ almost all $x\in\Rn$. Then $\mu$ is $m$-rectifiable if and only if $\int_0^{\infty}|C^m_{\mu}(x,r)|^2/r\,dr<\infty$ for $\mu$ almost all $x\in\Rn$.
\end{thm}

The sufficiency of this condition for rectifiability follows from the work of Mayboroda and Volberg \cite{MV09b} while Villa proved the necessity in \cite{Vil19b}. Villa also proved an analogous result for uniform rectifiability and considered more general kernels.

\subsection{Square functions}\label{square functions}
In \cite{DS93} David and Semmes proved the following 

\begin{thm}\label{squarecauchy}
Let $E\subset\C$ be AD-1-regular. For $f\in L^2(E)$ and $z\in\C\setminus E$ define
$$F(z)=\int_E\frac{f(w)}{z-w}\,d\H^1w.$$
Then $E$ is uniformly rectifiable if and only if 
$$\int_{\C}|F'(z)|^2d(z,E)\,dz \lesssim \int_E|f|^2\,d\H^1\ \text{for all}\ f\in L^2(E).$$
\end{thm}

See Theorem I.2.41 in \cite{DS93}, where also an $(n-1)$-dimensional version in $\Rn$ is given with the Cauchy transform replaced by the Riesz transform $R_{n-1}$. 

Here is a related result with harmonic functions

\begin{thm}\label{squareharm}
Let $\Omega\subset\Rn$ be a corkscrew domain (see Section \ref{Harmmeas}) with AD-$(n-1)$-regular boundary. 
Then $\partial\Omega$ is uniformly rectifiable if and only if 
$$\int_{\Omega\cap B(a,r)}|\nabla u(x)|^2d(x,\partial\Omega)\,dx \lesssim r^{n-1}|\sup\{u(x):x\in\Omega\}|^2$$
for all $a\in \partial\Omega, r>0$ and for every bounded harmonic function $u$ in $\Omega$.
\end{thm}

The 'only if' direction was proved by Hofmann, Martell and Mayboroda \cite{HMM16} and the 'if' direction by Garnett, Mourgoglou and Tolsa \cite{GMT18}. These papers show that this estimate is also equivalent to an approximation property of harmonic functions. Related results were proven by Hofmann and Tapiola \cite{HT20} and Bortz and Tapiola \cite{BT19}.

\subsection{Other related kernels}\label{otherkernels}

The Menger curvature trick \eqref{Cauchymenger} is particular to the Cauchy kernel $1/z$. We shall now discuss some positive results and counterexamples with other kernels. 

For $\Omega:S^1\to\C$, define $$k_{\Omega}(z) = \Omega(z/|z|)/|z|,\ z\neq 0.$$ 
Huovinen proved in \cite{Huo97} that Theorem \ref{pvmu} remains valid for $k_{\Omega_k}$ with $\Omega_k(z)=z^k/|z|^k$, where $k$ is an odd positive integer. Assuming additionally that the lower density is finite, it holds for finite linear combinations of such kernels. Now the tangent measures satisfy the $\Omega$-symmetricity $\int_{B(z,r)}\Omega(\zeta-z)\,d\nu\zeta = 0\ \text{for}\ z\in\spt\nu, r>0.$  Huovinen  characterized the supports of such non-discrete  measures. They are either unions of lines or the whole plane $\C$.  Observe that the cancellation for $z^k$ does not only come from  $-z$, but from $k-1$ other points, too. Thus, when $k\geq 3$, in addition to flat tangent  measures also certain finite sums of up to $k$ flat measures on lines through the origin, are possible tangent measures. These are called spike measures by Jaye and Merch\'an and they cause new problems as compared to the case $k=1$. Anyway, Huovinen was able to show that under positive lower density and existence of principal values only flat measures occur as tangent measures almost everywhere.

Jaye and Merch\'an \cite{JM21} proved Huovinen's $\Omega_k$ result assuming positive and finite upper density almost everywhere. New methods are required;  for $k\geq 3$ the Menger curvature is not available and tangent measures do not seem to work either. To deal with this, Jaye and Merch\'an  combined Tolsa's and Huovinen's methods. They introduced $\a$ numbers in the spirit of \eqref{alphaur}, but now minimizing the distance to spike measures. They showed that positive and finite upper density and existence of principal values imply that these $\a$ numbers tend to zero. This alone does not imply rectifiability, but combined with $L^2$-boundedness it does. The needed $L^2$-boundedness again follows by $T(b)$ theorems. The proof is also based on the earlier work \cite{JM20a} and \cite{JM20b} of these authors.



In \cite{Vil21} Villa proved the rectifiability result for odd bi-Lipschitz functions $\Omega:S^1\to S^1$ with constant close to 1 under positive lower density and finite upper density. Then the AD-1-regular $\Omega$-symmetric measures are 1-flat.

In \cite{Huo01} Huovinen considered kernels satisfying standard Calder\'on-Zygmund conditions and some additional cancellation conditions. Let  $K_t(z)=Re(z)/|z|^2-tRe(z)^3/|z|^4, t\in\R$. When $t=1$, he proved that there exist compact purely 1-unrectifiable AD-regular sets such that the principal values exist almost everywhere and the operator is $L^2$ bounded on some subset of positive measure. This phenomenon is caused by the cancellation coming from the coordinate axis. Jaye and Nazarov \cite{JN18} found for $t=3/4$ a compact purely 1-unrectifiable set with positive and finite 1-measure for which the operator is $L^2$ bounded. This set is not AD-regular but it has the interesting property that, despite the $L^2$ boundedness, the principal values do not exist. In fact, their kernel was much simpler, $\bar{z}/z^2$, but its real part is $4K_t$. Mateu and Prat \cite{MP21} gave an example in higher dimensions.

Chousionis, Mateu, Prat and Tolsa \cite{CMPT12} considered the kernels $Re(z)^k/|z|^{k+1}$ for positive odd integers $k$. They proved for them the analogues of Theorems \ref{CauchyL2} and \ref{MMV}. They did it by relating to these kernels a sum of permutations as in \eqref{Cauchymenger} and showed that it is non-negative and has properties similar to the Menger curvature. This allowed them to prove the analogue of the David-Leg\'er theorem \ref{davleg}. Chunaev \cite{Chu17} did the same for a larger class of kernels including $K_t$ as above for certain parameters $t$ for which the permutation sum is non-negative. But for some $t$ this sum takes both positive and negative values. Even for a range of such $t$ Chunaev, Mateu, and Tolsa \cite{CMT19} managed to prove analogous results.


David \cite{Dav01} and Chousionis \cite{Cho09} considered some self-similar fractals and constructed Calder\'on-Zygmund kernels for which the operators are bounded on these fractals. The kernels are defined to fit the self-similarities of the sets.

\section{Singular integrals}\label{singular integrals}

\subsection{A few words in general}

The classical Calder\'on-Zygmund theory of singular integrals deals with the Lebesgue measure on $\R^m$ and operators $T_K$;
\begin{equation}\label{TK}
T_Kf(x)=\int K(x-y)f(y)\,dy.\end{equation}
Here $K$ usually is nice except that it has a singularity of order $|x|^{-m}$ at the origin. Because of this the integral in \eqref{TK} often does not exist. But $K$ is also assumed to possess cancellation, being odd or something less. Then the principal values
\begin{equation}\label{pv}
T_Kf(x)=\lim_{\e\to 0}\int_{|x-y|>\e} K(x-y)f(y)\,dy\end{equation}
exist if $f$ is sufficiently nice, usually Lipschitz is enough.  This is easily checked for instance for the Riesz kernel $|x|^{-m-1}x$.

The core of the theory is the $L^2$-boundedness; when is $T_K:L^2\to L^2$ bounded? Very general regularity assumptions on $K$, saying that the singularity at $0$ is not too bad, suffice. For instance the condition that $|x|^{m+j}|\nabla^jK(x)|$ is bounded for $j=0,1,2,\dots$, that will appear later, is much more than enough. The $L^2$-boundedness implies, and often is equivalent to, many other fundamental properties ($L^p, 1<p<\infty,$ BMO, weak $L^1$, $T(b)$, etc.). In particular it implies that the convergence of \eqref{pv} takes place for all $f\in L^1$ for almost all $x\in\R^m$. This is because nice functions are dense in $L^1$.

Our interest is in the case where the Lebesgue measure is replaced by in some sense $m$-dimensional measure on $\Rn$. To define the $L^2$-boundedness without having the pointwise formula \eqref{pv}, we define the truncated operators $T_{K,\mu,\e}, \e>0$:
$$T_{K,\mu,\e}g(x)=\int_{|x-y|>\e}K(x-y)g(y)\,d\mu y.$$ 
We say that $T_{K,\mu}$ is bounded in $L^2(\mu)$ if the $T_{K,\mu,\e}, \e>0$, are uniformly bounded in $L^2(\mu)$. Often this is equivalent to saying that the maximal transform $T_{K,\mu}^{\ast}$,
$$T_{K,\mu}^{\ast}g(x)=\sup_{\e>0}|\int_{|x-y|>\e}K(x-y)g(y)\,d\mu y|,$$ 
is bounded in $L^2(\mu)$. In case $\mu=\H^m\restrict E$ we say that $T_{K,E}$ is bounded in $L^2(E)$.

Coifman and Weiss extended in \cite{CW71} most of the basic theory from the Lebesgue measure in $\R^m$ to doubling measures $\mu$ in metric spaces $E$, in particular to AD-$m$-regular sets in $\Rn$. Again a basic question is: when is $T_K:L^2(\mu)\to L^2(\mu)$ bounded? But now it is not only about the kernel, but also, and usually mainly, about $E$ and $\mu$. Still $L^2$-boundedness implies a lot of other things, but not any more automatically the existence of principal values. The cancellation properties of the kernel do not help even for constant functions if the measure does not have symmetry properties. But we showed with Verdera in \cite{MV09a} that under very general conditions $L^2$-boundedness implies weak convergence of the truncated operators.

Is the doubling condition necessary for the general theory? It seems that this question was seriously considered only when it was needed in connection of the Cauchy integral and analytic capacity, recall the $T(b)$-theorems of David and of Nazarov, Treil and Volberg from the previous chapter. In addition to $T(b)$, Nazarov, Treil and Volberg, and later also Tolsa, developed the non-homogeneous (i.e., non-doubling) Calder\'on-Zygmund theory with surprising success in many papers.

\subsection{$L^2$-boundedness and uniform rectifiability}\label{L2UR}

This topic has its origins in David's work in the 1980s, see \cite{Dav84} and \cite{Dav88}. Most of the basics  was developed by David and Semmes in \cite{DS91} and \cite{DS93}. The main problem related to rectifiability is the following conjecture of David and Semmes:

\begin{conj}\label{Davidsemmes}
Let $0<m<n$ be integers and let $E\subset\Rn$ be AD-$m$-regular. Then the \emph{Riesz transform} $R_{E}^m$ is bounded in $L^2(E)$ if and only if $E$ is uniformly $m$-rectifiable.
\end{conj}

The kernel of $R_{E}^m$ is the \emph{Riesz kernel} $R_m(x)=|x|^{-m-1}x, x\in\Rn\setminus\{0\}$, so $R_{E}^m=T_{R_m,E}$ according to the above notation. For a measure $\mu$ we shall also set $R_{\mu}^m=T_{R_m,\mu}$. The $L^2$-boundedness again means that the truncated operators $R^m_{E,\e}, \e>0$,
$$R_{E,\e}^mg(x)=\int_{\{y\in E:|x-y|>\e\}}R_m(x-y)g(y)\,d\H^my,$$ 
are uniformly bounded in $L^2(E)$. 

The boundedness of much more general singular integral operators on uniformly rectifiable sets was proved in \cite{DS91}. We shall return to this soon. The problem is the converse. It is known for $m=1,n-1$, and only then:

\begin{thm}\label{NToV}
Let $m=1$ or $m=n-1$ and let $E\subset\Rn$ be AD-$m$-regular. Then the Riesz transform $R_{E}^m$ is bounded in $L^2(E)$ if and only if $E$ is uniformly $m$-rectifiable.
\end{thm}

For $m=1$ the proof of Theorem \ref{MMV} gives this, too. So it is based on symmetrization as in \eqref{Cauchymenger}. For $m>1$ this method does not work because the corresponding sum takes both positive and negative values as was shown, more generally, by Farag in \cite{Far99} and \cite{Far00b}. 

The case $m=n-1$ is due to Nazarov, Tolsa and Volberg in \cite{NTV14a}. Their proof is very long and complicated and contains many brilliant ideas. Before saying a few words about it, let us look at a much weaker and simpler result which holds for all $m$, see \cite{Mat02}

\begin{pr}\label{Mat02}
Let $0<m<n$ and let $\mu\in\mathcal M(\Rn)$ be AD-$m$-regular. If the Riesz transform $R_{\mu}^m$ is bounded in $L^2(\mu)$, then for $\mu$ almost all $a\in \Rn$\ $\mu$ has some $m$-flat tangent measures at $a$.
\end{pr}

The proof is easy. A duality argument, recall the discussion around \eqref{Cauchymax}, gives a bounded function $g$ such that the maximal function $R_{\mu}^{m \ast}g$ is bounded. Then we can find an AD-$m$-regular tangent measure $\nu$ such that
\begin{equation}\label{mat02}
|\int_{r<|x-y|<R}|x-y|^{-m-1}(x-y)\,d\nu y|\leq C\ \text{for}\ x\in\spt\nu, 0<r<R.\end{equation}
Since $\spt\nu$ cannot be the whole space we can find an open ball $U$ disjoint from $\spt\nu$ such that there is $x\in\spt\nu\cap\partial U$. Any tangent measure $\pi$ of $\nu$ at $x$ is again AD-$m$-regular and satisfies \eqref{mat02}. It has support in a half-space $H$ with $0\in\spt\pi\cap\partial H$, so due to \eqref{mat02} $\spt\pi$ must be tangential to $\partial H$ at $0$.  The next tangent measure will have support in $\partial H$, perhaps with a non-constant density, but the last tangent measure will be flat, and it is also a tangent measure of $\mu$ by Preiss's 'tangent measures to tangent measures are tangent measures' principle, see \cite[Theorem 14.16]{Mat95}. 

Proposition \ref{Mat02} gives approximation by planes at \emph{some} arbitrarily small scales, but even for ordinary rectifiability we would need it at \emph{all} arbitrarily small scales.  As such it is useless for the proof of Theorem \ref{NToV}, but in \cite{NTV14a} Nazarov, Tolsa and Volberg used tangent measures to prove a quantitative form of this: sufficiently good pointwise boundedness of the Riesz transform at a range of scales gives a good approximation at some cube in the same range. Then the counter-assumption that the set in Theorem \ref{ubwgl} is not Carleson allows to build a Cantor structure based on the generalized dyadic cubes, mentioned in Section \ref{Basic tools}, where the cubes of nearly flatness of $\spt\mu$ (or rather of $\mu$) and cubes of non-flatness are alternating. The proof ends with clever applications of an extremal problem and a maximum principle. These have their origins in \cite{ENV14}. The maximum principle  requires harmonicity of the kernel and this is the main reason (or maybe the only reason) why the proof does not extend to $m<n-1$. The end result is that the $L^2$-boundedness implies the property of approximation with unions of planes as in Theorem \ref{ubwgl} and hence uniform  rectifiability. This is not a sketch of the proof, only a few tiles from a magnificent structure. 

For sets of finite measure we have the analogue of Theorem \ref{CauchyL2}:

\begin{thm}\label{RieszL2}
Let  $E\subset\R^n$ be $\H^{n-1}$ measurable with $\H^{n-1}(E)<\infty$. If $R^{n-1}_{E}$ is bounded in $L^2(E)$, then $E$ is $(n-1)$-rectifiable. 
\end{thm}

This was proved by Nazarov, Tolsa and Volberg in \cite{NTV14b}. They used a result of Eiderman, Nazarov and Volberg from \cite{ENV14} according to which $R^{n-1}_{E}$ is unbounded in $L^2(E)$ if $E$ has zero lower density. Thus they could assume that $\Theta_{\ast}^{n-1}(E,x)>0$ for $\H^{n-1}$ almost all $x\in E$. Then they used an argument of Pajot from \cite{Paj96} to find AD-$(n-1)$-regular measures $\mu_j$ such that $R_{\mu_j}^{n-1}$ is bounded in $L^2(\mu_j)$ and $\H^{n-1}\restrict E \leq \sum_j\mu_j$. This allowed them to conclude the proof by Theorem \ref{NToV}.

As in \ref{AT15} for $m=1$, there is a generalization of Theorem \ref{NToV}:

\begin{thm}\label{DT21}
Let $\mu\in\mathcal M(\Rn)$ with $\mu(B(x,r))\leq r^{n-1}$ for $x\in \Rn, r>0$. Then  
\begin{equation}\label{DT21eq}
\|R^{n-1}_{\mu}\|_{L^2(\mu)}^2+\mu(\Rn)\sim\int_0^{\infty}\int\beta^{n-1,2}_{\mu}(x,r)^2\frac{\mu(B(x,r))}{r^{n-1}}\,d\mu x\frac{1}{r}\,dr+\mu(\Rn).
\end{equation}
\end{thm}

The estimate $'\lesssim'$ was proved by Girela-Sarri\'on \cite{Gir19} for general $m$ and general kernels. The converse was proved by Tolsa \cite{Tol21} based on his paper with Dabrowski \cite{DT21a}.

The validity of $'\gtrsim'$ in Theorem \ref{DT21} for $1<m<n-1$ is an open question, as it should be, since the David-Semmes conjecture is then open. However, Jaye, Nazarov and Tolsa proved in \cite{JNT18} that the $L^2$-boundedness for all radial-type kernels $K\in\mathcal K_m(\Rn)$, see below, implies that the right hand side of \eqref{DT21eq} is finite.
 
Prat, Puliatti and Tolsa \cite{PPT21} extended Theorems \ref{NToV} and \ref{RieszL2} to  kernels which are gradients of the fundamental solutions of more general elliptic equations; $cR_{n-1}$ is the gradient of the fundamental solution of the Laplace equation. 
Let $A(x), x\in\Rn,$ be an $n\times n$ matrix with H\"older continuous entries satisfying the ellipticity conditions 
\begin{equation}\label{elliptic}
|\xi|^2 \lesssim A(x)\xi\cdot\xi\ \text{for all}\ \xi, x\in\Rn,
\end{equation} 
\begin{equation}\label{elliptic1}
A(x)\xi\cdot\eta\lesssim |\xi||\eta|\ \text{for all}\ \xi,\eta, x\in\Rn.
\end{equation} 
The equation
\begin{equation}\label{elliptic2}
L_Au(x):=-\text{div}(A\nabla u)(x)=0,
\end{equation} 
has a fundamental solution $\Gamma_A(x,y)$. 
The kernel then is its gradient, 
$K_A(x,y)=\nabla_x\Gamma_A(x,y)$. 
Technically the situation now is more complicated than for the Riesz kernel, but the authors of \cite{PPT21} managed with several modifications to follow the same strategy to prove the analogues of Theorems \ref{NToV} and \ref{RieszL2}.

Mas and Tolsa proved in \cite{MT14} a characterization of uniform rectifiability in terms of the $L^2$-boundedness of variations of the Riesz transform. They are defined by maximizing $\sum_{m\in\Z}|R^m_{\mu,\e_{m+1}}f(x)-R^m_{\mu,\e_m}f(x)|^2$ over decreasing sequences of positive numbers $\e_m$. 


Except for the Riesz kernels and some other particular cases, recall the discussion in Section \ref{otherkernels}, the question for which kernels Theorem \ref{NToV} holds is pretty much open. However, if we consider a large class of kernels the characterization of uniform rectifiability in all dimensions was already obtained by David and Semmes in \cite{DS91}. Let us denote by $\mathcal K_m(\Rn)$ the set of smooth real valued odd functions $K$ on $\Rn\setminus\{0\}$ such that $|x|^{m+j}|\nabla^jK(x)|$ is bounded for $j=0,1,2,\dots$. 

\begin{thm}\label{DSkernels}
Let $E$ be AD-$m$-regular. Then $E$ is uniformly rectifiable if and only $T_{K,E}$ is bounded in $L^2(E)$ for all kernels $K\in\mathcal K_m(\Rn)$.
\end{thm}

Already the boundedness with the kernels $K\in\mathcal K_m(\Rn)$ which are of the form $K(x)=x_jk(x)$, where $k$ is radial, is enough for uniform rectifiability, see \cite[Theorem I.2.59]{DS93}, for rectifiability this was proved in \cite{MP95b}. Related results are in \cite{Tol09}.

Jaye and Nazarov studied \emph{reflectionless} measures in several papers.  They have some resemblance to symmetric measures. A measure $\mu$ is reflectionless with respect to a kernel $K$ if $T_{K,\mu}1$ vanishes, in a weak sense, on the support of $\mu$. 
In \cite{JN19} they showed that if $K$ is a kernel with respect to which the AD-$m$-regular reflectionless measures are flat, then the $L^2$-boundedness on any  AD-$m$-regular set $E$ implies the uniform rectifiablity of $E$. This condition is known to hold for the Cauchy kernel and $R_1$ in the plane but it is unknown for higher dimensional Riesz kernels. With this property of the Cauchy kernel Jaye and Nazarov gave in \cite{JN14} the earlier mentioned new proof of Theorem \ref{MMV}. The example of \cite{JN18} mentioned near the end of Section \ref{otherkernels} relies on the fact that the Lebesgue measure on $\C$ is reflectionless with respect to the kernel $\bar z/z^2$.

\subsection{Principal values}
We have the analogue of Theorem \ref{Cauchypv} for the Riesz kernels:

\begin{thm}\label{Rieszpv}
Let $E\subset\Rn$ be $\H^m$ measurable with $\H^m(E)<\infty$. Then $E$ is $m$-rectifiable if and only if the finite limit
$$\lim_{\e\to 0}\int_{\{y\in E:|y-x|>e\}}R_m(y-x)\,d\H^my$$
exists for $\H^m$ almost all $x\in E$.
\end{thm}

Verdera's Hahn-Banach proof \cite{Ver92} for the Cauchy transform generalizes to give that the principal values exist for rectifiable sets. The converse  was proved by Tolsa in \cite{Tol08}. His proof for Theorem \ref{Cauchypv} does not work now as he used Menger curvature for that. It is kind of replaced by $L^2$ estimates on Lipschitz graphs. That is, Tolsa showed that the $L^2$ norm of the Riesz transform on Lipschitz graphs is quantitatively bounded by the $L^2$ norm of the gradient of the function, not only from above but also from below. These estimates are then used to construct Lipschitz graphs containing a positive measure of $E$. Very roughly, the existence of principal values implies some approximation with a graph of a Lipschitz function, whose gradient cannot have too big $L^2$ norm by the bounds on the Riesz transform. Arguments similar to those of Leg\'er in the proof of Theorem \ref{davleg} are essential. 

The analogue of Theorem \ref{pvmu} was proved in \cite{MP95b} but in addition to positive lower density we had to assume that it also is finite:

\begin{thm}\label{pvmun}
Let $\mu\in\mathcal M(\Rn)$. If $0<\Theta^m_{\ast}(\mu,x)<\infty$ and the finite limit $\lim_{\e\to 0}\int_{\{|y-x|>\e\}}R_m(y-x)\,d\mu y$ 
exists for $\mu$ almost all $x\in \Rn$, then $\mu$ is $m$-rectifiable.
\end{thm}

 Again we showed that tangent measures $\nu$ are symmetric:
$$\int_{B(x,r)}(y-x)\,d\nu y = 0\ \text{for}\ x\in\spt\nu, r>0.$$
We could not say very much about them, except that the AD-$m$-regular symmetric measures are flat. But we could show that if almost everywhere $\mu$ has positive lower density and the tangent measures are symmetric, then they are flat.

Combining theorems \ref{RieszL2} and \ref{Rieszpv} we see that if $m=1$ or $m=n-1$, $E\subset\Rn$ is $\H^{m}$ measurable with $\H^{m}(E)<\infty$ and the Riesz transform is bounded in $L^2(E)$, then the principal values exist $\H^m$ almost everywhere in $E$. But we only know this going through rectifiability. It would be interesting to have a direct proof, or any proof when $1<m<n-1$. Combining with Theorem \ref{Rieszpv} this would give a proof that the $L^2$-boundedness implies rectifiability.

The more general question for which kernels $L^2$-boundedness implies almost everywhere convergence is open, even for AD-regular sets. We know from Section \ref{otherkernels} that there are reasonable kernels for which this fails, but the Jaye-Nazarov example is not AD-regular. In \cite{MV09a} we proved with Verdera under very general conditions that the $L^2$-boundedness together with zero lower density implies the existence of principal values. Jaye and Merch\'an \cite{JM20a} strengthened this by replacing zero density by the condition that modifications of Tolsa's $\a$s (recall Section \ref{alphas}) tend to zero. See also \cite{JM20b} for related results and recall the discussion in \ref{otherkernels} on \cite{JM21}. 

Let $\Omega\subset\Rn$ be a domain with compact AD-$(n-1)$-regular (and a little more) boundary. D. and M. Mitrea and Verdera \cite{MMV16} proved that then $\Omega$ is a $C^{1+\a}, 0<\a<1,$ domain if and only if the Riesz transform $R^{n-1}_{\partial\Omega}$ maps $C^{\a}(\partial\Omega)$ into itself. The second half actually means several equivalent conditions, some of them involving uniform rectifiability, corresponding to different definitions of the Riesz transform, for example principal value and distributional definitions. The proof uses Clifford algebras and an interesting formula expressing the unit normal in terms of the Riesz transform and the Cauchy-Clifford transform.

\subsection{Lipschitz harmonic functions}
The kernel $R_{n-1}$ is a constant multiple of the gradient of the fundamental solution of the Laplacian in $\Rn$, which is $c_n|x|^{2-n}$ for $n\geq 3$,\ $c_2\log|x|$ for $n=2$. Hence the codimension 1 Riesz transform is related to the removability of Lipschitz harmonic functions in the same way as the Cauchy transform is related to the removable sets of bounded analytic functions.

Let us say that a compact set $E\subset\Rn$ is \emph{removable for Lipschitz harmonic functions}, abbreviated RLH, if whenever $U$ is an open set containing $E$ every Lipschitz function $u:U\to\R$ which is harmonic in $U\setminus E$ is harmonic in $U$. Since Lipschitz functions always have Lipschitz extensions, we can start with $u$ defined in all of $U$.

\begin{thm}\label{NToV1}
Let $E\subset\Rn$ be compact with $\H^{n-1}(E)<\infty$. Then $E$ is RHL if and only if $E$ is purely $(n-1)$-unrectifiable.
\end{thm}

The non-removability of rectifiable sets was proved in \cite{MP95a} by methods similar to those used for the analytic functions and discussed in the previous chapter. The converse was proved for $n=2$ in \cite{DM00} based on Menger curvature (recall the discussion after Theorem \ref{CauchyL2}) and the David-L\'eger theorem \ref{davleg}. The case of general $n$ is due to Nazarov, Tolsa and Volberg in \cite{NTV14b}. The proof is reduced to Theorem \ref{NToV} via a $T(b)$-theorem from \cite{Vol03} in a similar manner as was argued for analytic functions.

In the plane Tolsa's characterization, Theorem \ref{Tolsa}, is valid also for Lipschitz harmonic functions. Hence

\begin{thm}\label{Tolsa1}
Let  $E\subset\C$ be compact. Then the following are equivalent:
\begin{itemize}
\item[(1)] $E$ is removable for bounded analytic functions.
\item[(2)] $E$ is removable for Lipschitz harmonic functions.
\item[(3)] If $\mu\in\mathcal M(E)$ is such that $\mu(B(z,r))\leq r$ for $z\in \C$ and $r>0$, then $c^2(\mu)=\infty$. 
\end{itemize}
\end{thm}

The equivalence of (1) and (2) is only known passing through (3), and hence by a very complicated proof. It is easy to see that (1) implies (2): if $u$ is Lipschitz harmonic, the $\partial_{\overline{z}}\partial_zu=\Delta u = 0$, so $f=\partial_z u$ is bounded analytic. The converse, to get $u$ from $f$, would require some kind of integration, which is possible in some special cases but maybe not always.

Tolsa \cite{Tol21} characterized also the removable sets for Lipschitz harmonic functions in all dimensions:

\begin{thm}\label{Tol21} Let  $E\subset\Rn$ be compact. Then  $E$ is not removable for Lipschitz harmonic functions if and only if there exists $\mu\in\mathcal M(E)$ such that $\mu(B(x,r))\leq r^{n-1}$ for $x\in \Rn$ and $r>0$ and 
$$\int_0^{\infty}\int\beta^{n-1,2}_{\mu}(x,r)^2\frac{\mu(B(x,r))}{r^{n-1}}\,d\mu x\frac{1}{r}\,dr<\infty.$$
\end{thm}

This is a consequence of Theorem \ref{DT21}.

\subsection{Parabolic singular integrals}

In \cite{BHHLN20}, \cite{BHHLN21b} it is shown that a large class of parabolic singular integral operators are $L^2$-bounded on parabolic uniformly rectifiable sets.  Recall Section \ref{Parabolic}, as mentioned there the motivation for this theory comes from the heat equation $\Delta_xu(x,t) = \partial_tu(x,t) .$ Its fundamental solution $W$ is given for $x\in\Rn, t>0$, by
$$W(x,t)=ct^{-n/2}e^{-|x|^2/(4t)}.$$ 
The kernel $K$, which now replaces the Riesz kernel $R_{n-1}$, is
$$K(x,t)=\nabla_xW(x,t)=-(c/2)t^{-n/2-1}xe^{-|x|^2/(4t)},$$ 
with $K(x,t)=0$, when $t\leq 0$. Notice that it is antisymmetric only in the $x$ variable, but this is good enough for the $L^2$-boundedness on uniformly rectifiable sets by \cite{BHHLN20}. The converse, the analogues of Theorems \ref{NToV} and \ref{RieszL2}, are unknown. These would be needed to get a removability result such as Theorem \ref{NToV1} for parabolic regular Lipschitz (again BMO in $t$ variable, recall Section \ref{Parabolic}) solutions of the heat equation. Mateu, Prat and Tolsa have done in \cite{MPT20} some preliminary work in this direction. For instance, they showed that positive measure subsets of parabolic regular Lipschitz graphs are not removable. They also constructed Cantor sets with positive measure which are removable.

We shall discuss harmonic measure, induced by the Laplace equation, in the next chapter. In the same way the heat equation leads to the caloric measure. To get something like Theorem \ref{AHM16} in the parabolic case would also seem to require information about the consequences of the $L^2$-boundedness of $T_K$.

\subsection{Heisenberg groups}
As compared to the Euclidean theories, rather little is known about the singular integrals
$$T_Kf(p)=\int K(p^{-1}\cdot q)f(q)\,d\mu q$$
on lower than full dimensional subsets of the Heisenberg groups. As was seen in Chapter \ref{heisenberg} rectifiability is much better understood than uniform rectifability. Anyway intrinsic Lipschitz graphs are good candidates for basic uniformly rectifiable sets. So it makes sense to ask: For what kernels are the singular integral operators $L^2$-bounded on intrinsic Lipschitz graphs? When does $L^2$-boundedness or existence of principal values imply rectifiability?

What kernels should replace the Riesz kernels? If we just look at the expression and the scaling property of $R_m$ in $\Rn$ a similar kernel in $\hn$ is $R_m=(R_{m,1},R_{m,2})$ where for $p=(z,t)\in\hn$ and $\|p\|=(|z|^4+t^2)^{1/4}$,
$$R_{m,1}(p)=\|p\|^{-m-1}z\ \text{and}\ R_{m,2}(p)=\|p\|^{-m-2}t.$$
In \cite{CM11} an analogue of Proposition \ref{Mat02} was proved for these kernels. However, if we want a connection to harmonicity, we should start from a Laplacian. Let $X_i,Y_i, i=1,\dots,n,$ be the vector fields as is Section \ref{Heistools}. The sub-Riemannian or Kohn Laplacian in $\hn$ is defined by
$$\Delta_H=\sum_{i=1}^n(X_i^2+Y_i^2).$$
For the potential theory related to it, see \cite{BLU07}. The fundamental solution of $\Delta_Hu=0$ is $\Gamma(p)=c\|p\|^{-2n}$. Note that $2n=\dim\hn-2$, so it has the same form as in $\Rn$. The kernel related to the Lipschitz harmonic functions is $K=\nabla_H\Gamma:\hn\to\R^{2n}$. It looks a bit complicated, since its coordinate functions are for $i=1,\dots,n$, and $z=(x,y)\in\R^{2n}$,
$$K_i(z,t)=\frac{x_i|z|^2+y_it}{\|p\|^{2n+4}},\ K_{i+n}(z,t)=\frac{y_i|z|^2-x_it}{\|p\|^{2n+4}}.$$
However, it is a reasonable Calder\'on-Zygmund kernel. It is not odd, but it is horizontally antisymmetric: $K(z,t)=-K(-z,t)$. Due to this Chousionis, F\"assler and Orponen \cite{CFO19b} were able to prove in $\h^1$ that $T_K$, and more general singular integrals, are $L^2$-bounded on $C^{1,\a}$ intrinsic graphs, see also \cite{FO18} for related results. 

Orponen \cite{Orp18} proved that for AD-3-regular subsets of $\h^1$  the $L^2$-boundedness of the 3-dimensional singular integrals with horizontally antisymmetric kernels implies local symmetry of the type in Theorem \ref{LSthm} and that this implies weak geometric lemma.

In \cite{CM14} a class of self-similar purely $(2n+1,\h)$-unrectifiable subsets of $\hn$ with positive $\H^{2n+1}$ measure was introduced on which $T_K$ is $L^2$ unbounded and which are removable for Lipschitz harmonic functions. Further results in this direction were proven by Chousionis and Urbanski in \cite{CU15}.

Chousionis and Li \cite{CL17} introduced a class of non-negative 1-homogeneous kernels in $\h^1$ which vanish on the vertical axis $\{z=0\}$. For some of them the operator is $L^2$-bounded on regular curves and for some  the $L^2$-boundedness on an AD-1-regular set implies that it is contained in a regular curve. Extension of the first statement to general Carnot groups is given by Chousionis, Li and Zimmerman \cite{CLZ19b}.

F\"assler and Orponen \cite{FO21} proved in $\h^1$ the $L^2$-boundedness of many singular integral operators on AD-regular curves and on a class of vertical graphs called Lipschitz flags. 

\section{Harmonic measure and elliptic measures}\label{Harmonic, elliptic measures}
For the harmonic measure in the plane, see \cite{GM05}. Toro's survey \cite{Tor19} discusses relations between geometric measure theory and harmonic measure both in the plane and in higher dimensions.

\subsection{Harmonic measure}\label{Harmmeas}

Let $\Omega\neq\Rn$ be an open connected subset of $\Rn, n\geq 2$. We assume that $\partial\Omega$ is not too small; it is not a polar set. This is true if  $\H^{n-1}(\partial\Omega)>0$. By classical potential theory, see, for example, \cite{AG01}, for any continuous function $f$ on $\partial\Omega$ one can solve the Dirichlet problem to find a harmonic function $u_f$ in $\Omega$ with boundary values $f$, in a generalized sense if $\partial\Omega$ is not sufficiently nice. Fixing  $p\in\Omega$ for a while the Riesz representation theorem and the maximum principle can be used to show that there is a probability measure $\omega_{\Omega}^p\in\mathcal M(\partial \Omega)$ such that 
\begin{equation}\label{hm}
u_f(p)=\int f\,d\omega_{\Omega}^p.\end{equation} 
Then $u(p)=\int u\,d\omega_{\Omega}^p$ if $u$ is continuous in $\overline{\Omega}$ and harmonic in $\Omega$. The measure $\omega_{\Omega}^p$ is called the \emph{harmonic measure} of $\Omega$ at $p$. It depends on $p$, but by the Harnack inequality for any two points $p,q\in\Omega$ the measures $\omega_{\Omega}^p$ and  $\omega_{\Omega}^q$ are comparable.

According to Kakutani's probabilistic characterization  $\omega_{\Omega}^p(A)$ is the probability that the Brownian traveller starting from $p$ hits $A$ before hitting any other part of the boundary. This helps to visualize the fact that in the case of complicated boundaries $\omega_{\Omega}^p$ lives on  parts of the boundary which are more easily accessible from $\Omega$. A more precise statement in the plane is a result of  Wolff \cite{Wol93} saying that harmonic measure lives on a set of sigma-finite $\H^1$ measure. The corresponding statement in $\R^3$ is false by his example in \cite{Wol95}. However, Bourgain \cite{Bou87} proved that in $\Rn$ it lives on a set of dimension at most $n-\e(n)$, where $\e(n)$ is a small positive constant, whose best value is unknown.

From the point of view of this survey the main question is: what are the relations between harmonic measure and the geometry of the boundary? More precisely, for $E\subset\partial\Omega$ does rectifiability of $E$ imply something on the harmonic measure on $E$, and conversely, do some properties of the harmonic measure lead to rectifiability? The first general result was proved already in 1916 by the brothers Riesz: if $\Omega\subset\R^2$ is simply connected and $\H^1(\partial\Omega)<\infty$, then $\omega_{\Omega}^p$ and $\H^1\restrict\partial\Omega$ are mutually absolutely continuous. Here the whole boundary is rectifiable but Bishop and Jones extended this in \cite{BJ90} by showing that if $\Omega\subset\R^2$ is simply connected and $E\subset\partial\Omega\cap\Gamma$, where $\Gamma$ is a rectifiable curve, then $\omega_{\Omega}^p(E)=0$ if and only if $\H^1(E)=0$. They used Jones's traveling salesman theorem \ref{jones} for this.

In the plane complex analytic tools are very effective, but in higher dimensions quite different methods are required and the situation is in many ways different. First, the analogue of the Riesz brothers theorem fails. Ziemer constructed in \cite{Zie74} a domain in $\R^3$ whose boundary is 2-rectifiable, it even has an ordinary tangent plane a each point, and it has finite $\H^2$ measure, but $\H^2$ is not absolutely continuous with respect to harmonic measure. To the other direction, Wu showed in \cite{Wu86} that there exists a domain $\Omega\subset\R^3$ and $E\subset\partial\Omega\cap V$, where $V$ is a plane, with positive harmonic measure and $\H^2(E)=0$, even $\dim E=1$.

Starting from Dahlberg \cite{Dah77} in 1977 and Lipschitz boundaries and followed by David, Jerison, Kenig and Semmes and more general boundaries in the 1980s a great number of people have produced results in the spirit that various geometric properties of the boundary imply absolute continuity of harmonic measure, often with quantitative estimates such as being an $A_{\infty}$ weight. Uniform rectifiability of the boundary alone is not sufficient to get such results by an example of Bishop and Jones in \cite{BJ90}. But starting with some natural geometric conditions of the boundary, 
uniform rectifiability often comes into play leading to, or even characterizing, quantitative properties of the harmonic measure. 
Commonly used conditions are corkscrew and Harnack chain conditions. Roughly speaking the former says that every ball centered in the boundary contains a ball of comparable size inside the domain and the latter that any two points in the domain can be joined with a chain of balls whose size is comparable to the distance to the boundary. 

The proof of the following theorem was completed by Azzam,  Hofmann,  Martell, Mourgoglou and Tolsa in \cite{AHMMT20}. The authors explain by several examples that the result is in many ways optimal. It is a culmination of a long process involving many other people and articles, see \cite{AHMMT20} and \cite{HLMN17} for the history, references and other related results. In particular, that $A_{\infty}$ implies uniform rectifiability was proved earlier by Hofmann and Martell,  see \cite{HLMN17} where this is obtained even for the non-linear $p$-harmonic equation. 

\begin{thm}\label{AHMMT20} Let $\Omega$ satisfy the corkscrew condition and have AD-$(n-1)$-regular boundary. Then  $\partial\Omega$ is uniformly rectifiable and satisfies the weak local John condition if and only if $\omega_{\Omega}^p$ is locally in weak $A_{\infty}$.  
\end{thm}

The weak local John condition is a quantitative connectivity condition saying, roughly, that each point of $\Omega$ can be connected to many points in the boundary by rectifiable curves in $\Omega$ staying away from the boundary. That $\omega_{\Omega}^p$ is locally in weak  $A_{\infty}$ means that there is $s>0$ such that for every $x\in\partial\Omega$ and $0<r<d(\partial\Omega)/4$,
$$\omega_{\Omega}^p(A)\lesssim \left(\frac{\H^{n-1}(A)}{\H^{n-1}(\partial\Omega\cap B(x,r))}\right)^s\omega_{\Omega}^p(\partial\Omega\cap B(x,2r))$$
for all $p\in\Omega\setminus B(x,4r)$ and for all Borel sets $A\subset\partial\Omega\cap B(x,r)$.

The local weak $A_{\infty}$ condition is known to be equivalent to the quantitative solvability of the $L^p$ Dirichlet problem.
 
Let us now look at a qualitative rectifiability criterion. After many partial results by several people Azzam, Hofmann, Martell, Mayboroda, Mourgoglou,  Tolsa and Volberg proved in \cite{AHMMMTV16}

\begin{thm}\label{AHM16} Let $p\in\Omega$ and $E\subset\partial\Omega$ with $\H^{n-1}(E)<\infty$. 
\begin{itemize}
\item[(1)] If $\omega_{\Omega}^p\restrict E\ll\H^{n-1}\restrict E$, then $\omega_{\Omega}^p\restrict E$ is $(n-1)$-rectifiable.
\item[(2)] If $\H^{n-1}\restrict E\ll\omega_{\Omega}^p\restrict E$, then $E$ is $(n-1)$-rectifiable.
\end{itemize}
\end{thm}

By the Radon-Nikodym theorem it is easy to see that these statements are equivalent. For example, if (1) holds and $\H^{n-1}\restrict E\ll\omega_{\Omega}^p$, then $\H^{n-1}\restrict E = g\omega_{\Omega}^p$ for some non-negative function $g$ on $E$ for which $g(x)>0$ for $\H^{n-1}$ almost all $x\in E$. Then $\omega_{\Omega}^p\ll\H^{n-1}\restrict \{x\in E:g(x)>0\}$, and it follows from (1) that $E$ is $(n-1)$-rectifiable.

The key to the proof of Theorem \ref{AHM16} is the relation between the Green function and the Riesz transform and Theorem \ref{RieszL2} of Nazarov, Tolsa and Volberg. In classical potential theory the \emph{Green function} $G_{\Omega}:\Omega\times \Omega\setminus\{(p,x):p=x\}\to\R$ is a basic tool to study harmonic measure. It is harmonic in both variables separately. For a fixed $p\in\Omega$\ $G(p,\cdot)$ has zero boundary values. At $x=p$ it has the same singularity as the fundamental solution $\Gamma$ of the Laplacian, which is a constant multiple of $|x|^{2-n}$, if $n\geq 3$, and of $\log|x|$, if $n=2$. More precisely, $\Gamma(p-x)-G_{\Omega}(p,x)$ is a harmonic function of $x$ in $\Omega$.  For nice domains the harmonic measure is absolutely continuous with respect to the surface measure and its density is given by the normal derivative of the Green function. By \eqref{hm} the Green function can be written as 
$$G_{\Omega}(p,x)=\Gamma(p-x) - \int \Gamma(p-x)\,d\omega_{\Omega}^px,\ p,x\in\Omega, x\neq p.$$
Since $\nabla\Gamma = cR_{n-1}$ we have
\begin{equation}\label{green}
\nabla_xG_{\Omega}(p,x)=cR_{n-1}(p-x) - c\int R_{n-1}(p-x)\,d\omega_{\Omega}^px,\ p,x\in\Omega, x\neq p.
\end{equation}
As in (1) of Theorem \ref{AHM16} suppose that $\omega_{\Omega}^p\restrict E\ll\H^{n-1}\restrict E$, so that $\omega_{\Omega}^p\restrict E=g\H^{n-1}\restrict E$ for some non-negative $g$. Given $M>0$ it is enough to prove that $E_M:=\{x\in E:1/M<g(x)<M\}$ is $(n-1)$-rectifiable. One could then hope  that, similarly to the case of nice boundaries, the left hand side of \eqref{green} would have enough boundedness to give boundedness for the Riesz transform when $x\in E_M$. In this very general case this is not clear at all, but the authors of \cite{AHMMMTV16} managed to show something like this. A bit more precisely, they again used generalized dyadic cubes, now from \cite{DM00} since there is no doubling, and they showed using \eqref{green} that there are enough cubes $Q$ to cover large part of $E_M$ such that the truncated Riesz transform $R^{n-1}_{\omega_{\Omega}^p,r(Q)}(x)$ with a suitable $r(Q)$ is bounded for $x\in Q$, with a quantitative bound. 
This allows to apply a $T(b)$-theorem of Nazarov, Treil and Volberg to get the $L^2(\omega_{\Omega}^p)$-boundedness of the Riesz transform on a subset of $E_M$ with positive measure. From this an application of Theorem \ref{RieszL2} yields rectifiability.


There are also results on two-phase problems involving rectifiability. The following was proved by Azzam, Mourgoglou, Tolsa and Volberg in \cite{AMTV19}, the paper \cite{AMT17} contains a quantitative version:

\begin{thm}\label{2phase}
Let $\Omega_1$ and $\Omega_2$ be disjoint open connected subsets of $\Rn$ and $E\subset\partial\Omega_1\cap\partial\Omega_2$ a Borel set such that $\omega_{\Omega_1}^{p_1}$ and $\omega_{\Omega_2}^{p_2}, p_i\in\Omega_i,$ are mutually absolutely continuous on $E$. Then $E$ contains an $(n-1)$-rectifiable subset $F$ such that 
$\omega_{\Omega_i}^{p_i}(E\setminus F)=0$ and $\omega_{\Omega_1}^{p_1}$ and $\omega_{\Omega_2}^{p_2}$ are mutually absolutely continuous with respect to $\H^{n-1}\restrict F$.
\end{thm}

The proof uses an interesting blow-up argument involving tangent measures and the Alt-Caffarelli-Friedman monotonicity formula for pairs of subharmonic functions applied to the Green functions. This method was introduced by Kenig, Preiss and Toro in \cite{KPT09}, where a partial result  and deep information about harmonic measures was obtained. The proof also relies on Theorem \ref{AHM16}, and so on the Nazarov-Tolsa-Volberg Riesz transform theorems \ref{NToV} and \ref{RieszL2}, and on \cite{GT18}. 


\subsection{Elliptic measures in codimension 1}
The previous section dealt with the Laplace equation $\Delta u = 0$. It is natural to expect that the results would have analogues for more general elliptic equations, and this indeed is the case. Assuming that $\Omega$ is a uniform domain,  Hofmann, Martell, Mayboroda, Toro and Zhao \cite{HMMTZ20} characterized the $A_{\infty}$ property with uniform rectifiability for elliptic measures corresponding to equations with optimal conditions for the coefficients. Prat, Puliatti and Tolsa proved in \cite{PPT21} an analogue of Theorem \ref{AHM16} for such elliptic measures  with H\"older continuous coefficients using their singular integral results mentioned in Section \ref{L2UR}. Then the same arguments as for the Laplace equation work.

Cao, Hidalgo-Palencia and Martell \cite{CHM22} investigated corona decompositions associated to quite general elliptic measures. They showed, among other things, that these are equivalent to square function estimates as in Theorem \ref{squareharm}, also to a weaker form of them, which are equivalent to uniform rectifiability. The boundary of the domain is assumed to be AD-regular and to satisfy the cork-screw condition. 
\subsection{Elliptic measures in codimension bigger than one}

If $E\subset\Rn$ is a closed set with $\H^{n-2}(E)<\infty$, then it is polar. In particular, if $E\subset\partial\Omega$, then $\omega_{\Omega}^p(E)=0$ and the properties of the harmonic measure are in no way related to the geometric properties of $E$. The same is true for the elliptic equations as above.
But considering suitable degenerate elliptic equations a  rich theory can be developed. This was realized and done by David, Feneuil and Mayboroda and their co-authors in many papers, only some of which are listed in the references.

Let $E\subset\Rn$ be AD-$m$-regular for some integer $0<m<n-1$. The standard ellipticity conditions \eqref{elliptic} and \eqref{elliptic1} for a measurable $n\times n$ matrix valued function $A$ are now replaced by

\begin{equation}\label{elliptic3}
|\xi|^2 \lesssim d(E,x)^{n-m-1}A(x)\xi\cdot\xi\ \text{for all}\ \xi, x\in\Rn,
\end{equation} 

\begin{equation}\label{elliptic4}
d(E,x)^{n-m-1}A(x)\xi\cdot\zeta\lesssim |\xi||\zeta|\ \text{for all}\ \xi,\zeta, x\in\Rn.
\end{equation} 

The authors developed in \cite{DFM17} a comprehensive theory for the degenerate elliptic operators $L=-\di(A\nabla)$ analogous to the classical theory including among others solutions of the Dirichlet problem with continuous boundary data, the corresponding elliptic measure $\omega^p_L$ and its basic properties, and the existence and properties of the Green function and its relations to the elliptic measure. As remarked before, in the classical theory one often needs some conditions for the boundary to prevent it from being too massive. Here they are not needed. The AD-$m$-regularity with $m<n-1$ gives automatically the corkscrew and Harnack chain conditions.

To get absolute continuity of the elliptic measure on Lipschitz graphs or more general sets one needs stronger assumptions on $A$, also in the classical theory. I now restrict to a particular operator, which is the authors' replacement of the Laplacian. Again there is a weight like $d(E,x)^{n-m-1}$, but this seems to be too rough and it is replaced by a regularized distance $D_{\a,\mu}$, where $\a>0$ is a parameter and $\mu=\H^m\restrict E$ or some other AD-$m$-regular measure with support $E$:

$$D_{\a,\mu}(x)=\left(\int|x-y|^{-m-\a}\,d\mu y\right)^{-1/\a}.$$
AD-regularity easily implies that $D_{\a,\mu}(x)\sim d(E,x)$.

The elliptic operator attached to $\mu$ is given by
\begin{equation}
L_{\a,\mu}=-\di(D_{\a,\mu}^{m+1-n}\nabla).
\end{equation}
Denote now simply by $\omega^p$ the elliptic measure related to $L_{\a,\mu}$ and some base point $p\in\Rn\setminus E$.  In this setting David and Mayboroda proved in \cite{DM20a}

\begin{thm}\label{Davidmay}
If $E$ is uniformly $m$-rectifiable, $m\leq n-2$, then $\omega^p\ll\mu$. Moreover $\omega^p\in A_{\infty}(\mu)$, which means that for every $\e>0$ there is $\delta>0$ such that if $x\in E, r>0$ and $F\subset E\cap B(x,r)$ is a Borel set, then
\begin{equation}\label{Ainfty}
\frac{\omega^p(F)}{\omega^p(E\cap B(x,r))}<\delta \implies \frac{\mu(F)}{\mu(E\cap B(x,r))}<\e,
\end{equation}
where $p\in\Rn\setminus E$ is such that $d(E,p)\sim |p-x|\sim r$.
\end{thm}

In \cite{Fen20} Feneuil gave a different simpler proof. 

The converse of Theorem \ref{Davidmay} fails completely for $\a=n-m-2$, then $m<n-2$. In this case  $\omega\in A_{\infty}(\mu)$ for any AD-$m$-regular measure $\mu$, and $m$ need not even be an integer. I shall say a bit more about this soon. The above authors believe that this value of $\a$ is a unique exception and the converse should be true for the other values.

In addition to extending from codimension one to general dimensions the work of David, Feneuil, Mayboroda and their co-authors contains several interesting new results and phenomena also for codimension one. I present here two new characterizations of uniform rectifiability. Now $m$ can be any integer with $0<m<n$.

In \cite{DEM18} David, Engelstein and Mayboroda characterized uniform rectifiablity in terms of the distance $D_{\a,\mu}$ in the spirit of Theorems \ref{squarecauchy} and \ref{squareharm}. 

\begin{thm}
Let $E\subset\Rn$ be AD-$m$-regular. Then $E$ is uniformly rectifiable if and only if 
\begin{equation}\label{distgrad}
d(x,E)^{m+2-n}|\nabla(|\nabla D_{\a,\mu}|^2)(x)|^2dx\end{equation} is a Carleson measure on $\Rn\setminus E$.
\end{thm}

The 'only if' direction uses Tolsa's $\a$s, recall Theorem \ref{ATTthm}; \eqref{distgrad} vanishes if $\mu$ is a flat measure and the uniform rectifiability leads to approximation of $\mu$ with flat measures. For the converse direction the starting idea is that if \eqref{distgrad} vanishes then $\nabla D_{\a,\mu}$ is constant, from which it follows, but with a lot of work, that $\mu$ is $m$-flat.
They also characterized the rectifiability of $E$ in terms of the non-tangential limits of $\nabla| D_{\a,\mu}|$. 

Another characterization is in terms of the Green functions. Let $G_{\a,\mu}$ be the Green function corresponding to $L_{\a,\mu}$ with the pole at $\infty$. If $E$ is an $m$-plane, then $G_{\a,\mu}(x)=cd(E,x)$,  and the elliptic measure is a constant multiple of $\H^m\restrict E$. Uniform rectifiability means that $E$ is well approximated by $m$-planes and this turns out to be equivalent that $G_{\a,\mu}$ is well approximated by distances to $m$-planes:

\begin{thm}
Let $E\subset\Rn$ be an unbounded AD-$m$-regular set. When $m=n-1$ assume also the corkscrew and Harnack chain conditions. Then $E$ is uniformly rectifiable if and only if for every $\e>0$ and $M>1$ the set $E\times (0,\infty)\setminus \mathcal G(\e,M)$ is a Carleson set, where $\mathcal G(\e,M)$ is the set of $(x,r)\in E\times (0,\infty)$ such that there are an affine $m$-plane $V$ and $c>0$ for which
$$|d(y,V)-cG_{\a,\mu}(y)|\leq  \e r\ \text{for}\ y\in (\Rn\setminus E)\cap B(x,Mr).$$
\end{thm}

This is a special case of the results proved by David and Mayboroda in \cite{DM20b}. Again the proof of the 'only if' direction uses Tolsa's $\a$s. For the converse direction we have $G_{\a,\mu}(y)=0$ when $y\in E$, which immediately gives some approximation of $E$ by planes. So it is a good starting point, but it only gives weak geometric lemma and many more arguments are needed to get the bilateral approximation, recall Theorem \ref{bwgl}. 

The result in \cite{DM20b} is more general in that the authors considered much more general degenerate elliptic operators. 
They also have similar results with the regularized distance for $m=n-1$ where $d(y,V)-cG_{\a,\mu}(y)$ is replaced by $D_{\beta,\mu}(y)-cG_{a,\mu}(y), \beta>0$. Then they can start with any $n-2<m<n$ and the Carleson estimates force $m$ to be an integer. 

For $m < n-2$ and $\a=n-m-2$ the smallness of $|D_{\a,\mu}-cG_{a,\mu}|$ or the absolute continuity of the elliptic measure do not imply any kind of rectifiability, as was already stated after Theorem \ref{Davidmay}. The reason is that then by direct computation $L_{\a,\mu}D_{\a,\mu}=\tfrac{1}{\a}\Delta R_{\a,\mu}$ where, when $\a=n-m-2$, $R_{\a,\mu}(x)=\int|x-y|^{2-n}\,d\mu y$  is harmonic, so $L_{\a,\mu}D_{\a,\mu}=0$, from which it follows by the uniqueness of the Green function that $G_{\a,\mu}$ is a constant multiple of $D_{\a,\mu}$. This leads to $\omega_{L_{\a,\mu}}\sim\mu$. Further, if $E$ is $m$-rectifiable, using basic properties of rectifiable sets it follows that $\omega_{L_{\a,\mu}}=c\H^m\restrict E$. The details can be found in \cite{DEM18}.

\section{Sets of finite perimeter and functions of bounded variation}
This topic is well documented in many books, see \cite{Fed69}, Section 4.5, \cite{Sim83}, \cite{Giu84}, \cite{Zie89}, \cite{EG92}, \cite{AFP00}, \cite{LY02} and \cite{Mag12}. The structure theory of sets of finite perimeter is due to De Giorgi in the 1950s based on the earlier work of Caccioppoli, see, in particular, \cite{Deg55}. English translations of many of De Giorgi's papers can be found in \cite{Deg05}. Maggi \cite{Mag12} follows rather closely De Giorgi's original ideas. 

\subsection{Sets of finite perimeter}\label{Finper}

What is the perimeter of an arbitrary Lebesgue measurable set $E$ in $\Rn$? Even for open sets the right notion clearly is not the $\H^{n-1}$ measure of the topological boundary $\partial E$. For example, if $E$ is a countable union of balls $B_i\subset B(0,1)$ with $\sum_i\H^{n-1}(B_i)<\infty$, the topological boundary could be almost anything, in particular $B(0,1)$ if the centers are dense, but a more reasonable notion of perimeter would seem to be $\sum_i\H^{n-1}(B_i)$. The Gauss-Green theorem gives a hint how to define a good general notion of perimeter. If $E$ has smooth boundary, then for any compactly supported $C^1$ vectorfield $\phi$, $\phi\in C^1_c(\Rn)$,
\begin{equation}\label{perclass}
\int_E\di \phi = \int_{\partial E}\phi\cdot n_E\,d\H^{n-1},
\end{equation} 
where $n_E$ is the outer unit normal of $E$. Among $\phi$ with $|\phi|\leq 1$ the right hand side is maximized when $\phi=n_E$ on $\partial E$ and yields $\H^{n-1}(\partial E)$. The left hand side is defined for all measurable sets $E$, so we can define

\begin{df}\label{perimeter}
The \emph{perimeter} of a Lebesgue measurable set $E\subset\Rn$ is
$$P(E)=\sup\{\int_E\di\phi:\phi\in C^1_c(\Rn), |\phi|\leq 1\}.$$
If $P(E)<\infty$, we say that $E$ is a set of finite perimeter.
\end{df}

That $E$ is a set of finite perimeter means that the characteristic function $\chi_E$ is of bounded variation:

\begin{df}\label{BV}
A Lebesgue integrable function $u$ on $\Rn$ is of \emph{bounded variation}, $u\in BV(\Rn)$, if
$$\sup\{\int u\di\phi:\phi\in C^1_c(\Rn), |\phi|\leq 1\}<\infty.$$
\end{df}

Then $u\in BV(\Rn)$ if and only its distributional partial derivatives are finite Radon measures. That is, $Du$ is a vector valued Radon measure. I shall discuss them a bit more in the next section.

It is easy to see that the perimeter is lower semicontinuous: if $\chi_{E_i}\to \chi_E$ in $L^1(\Rn)$, then $P(E)\leq\liminf_{i\to\infty}P(E_i)$. From this we see that if, as in the beginning of this chapter, $E$ is a countable union of balls $B_i\subset B(0,1)$ with $\sum_i\H^{n-1}(B_i)<\infty$, then $E$ has finite perimeter. Moreover, so does $B(0,1)\setminus E$. 

By abstract arguments based on the Riesz representation theorem we have

\begin{thm}\label{perirepr}
Let $E\subset\Rn$ be a set of finite perimeter. There are $\mu_E\in\mathcal M(\Rn)$ and a Borel function $\nu_E:\Rn\to\Rn$ such that $|\nu_E(x)|=1$ for $\mu_E$ almost all $x\in\Rn$ and 
\begin{equation}\label{Gaussgreen}
\int_E\di\phi = \int\phi\cdot\nu_E\,d\mu_E\ \text{for}\ \phi\in C^1_c(\Rn).
\end{equation}
\end{thm}

Recalling \eqref{perclass} it is natural to call $\mu_E$ the generalized perimeter measure of $E$ and $\nu_E$ its generalized outer normal.

A fairly easy but important result is the compactness theorem:

\begin{thm}\label{pericomp}
If $E_j\subset B(0,1), j=,1,2,\dots,$ are Lebesgue measurable with\\ $\sup_jP(E_j)<\infty$, then there is a subsequence $(E_{j_i})$ and a set $E$ with finite perimeter such that $\chi_{E_{j_i}}\to \chi_E$ in $L^1(\Rn)$ and $\mu_{E_{j_i}}\to\mu_E$ weakly.
\end{thm}

\begin{df}\label{essbdry}
The \emph{reduced boundary} $\partial^{\ast}E$ of $E$ is the set of points $x\in\spt\mu_E$ such that $|\nu_E(x)|=1$ and
\begin{equation}\label{essbdryeq}
\lim_{r\to 0}\frac{1}{\mu_E(B(x,r))}\int_{B(x,r)}\nu_E\,d\mu_E
\end{equation}
exists and has norm 1.
\end{df}

It follows from the general theory of differentiation of measures that
\begin{equation}\label{essbdry1}
\mu_E(\Rn\setminus \partial^{\ast}E)=0.
\end{equation}

Formula \eqref{Gaussgreen} is a very general, but abstract, form of the Gauss-Green theorem. In order to make it more concrete we should understand better what $\mu_E$, $\nu_E$ and $\partial^{\ast}E$  really are. This is included in De Giorgi's structure theorem:

\begin{thm}\label{peristru}
Let $E\subset\Rn$ be a set of finite perimeter. Then $\partial^{\ast}E$ is $(n-1)$-rectifiable, $\mu_E=\H^{n-1}\restrict\partial^{\ast}E$, for $\H^{n-1}$ almost all $x\in\partial^{\ast}E$ the approximate tangent $(n-1)$-plane of $\partial^{\ast}E$ is $\{y:(y-x)\cdot\nu_E(x)=0\},$ and 
\begin{equation*}\label{Gaussgreen1}
\int_E\di\phi = \int_{\partial^{\ast}E}\phi\cdot\nu_E\,d\H^{n-1}\ \text{for}\ \phi\in C^1_c(\Rn).
\end{equation*}
\end{thm}

I say a few words about the main steps of the proof, which themselves are of independent interest.
First, there are the isoperimetric inequality
\begin{equation}\label{isoper}
\mathcal L^n(E)^{(n-1)/n}\lesssim P(E)=\mu_E(\Rn)
\end{equation}
and the local isoperimetric inequality for every ball $B\subset\Rn$,
\begin{equation}\label{isoper1}
\min\{\mathcal L^n(B\cap E)^{(n-1)/n}, \mathcal L^n(B\setminus E))^{(n-1)/n}\}\lesssim \mu_E(B).
\end{equation}

These follow from, and are in fact equivalent to, Sobolev and Poincar\'e inequalities for BV functions, which in turn follow from the classical inequalities and the fact that the smooth functions are dense in $BV(\Rn)$. 

The isoperimetric inequalities lead to density estimates, both for the Lebesgue measure and $\mu_E$. The key for deriving these is the identity
$$\int_{E\cap B(x,r)}\di\phi = \int_{B(x,r)}\phi\cdot\nu_E\,d\mu_E + \int_{E\cap \partial B(x,r)}\phi\cdot\nu\,d\H^{n-1}$$
for $\phi\in C^1_c(\Rn)$, where $\nu$ is the outward unit normal of $B(x,r)$. This follows applying \eqref{Gaussgreen} to a $C^1$ approximation of the characteristic function of $B(x,r)$. Then we have

\begin{lm}
There are positive constants $c$ and $C$ depending only on $n$ such that if $x\in\partial^{\ast}E$, then for all sufficiently small $r>0$,
$$\mathcal L^n(E\cap B(x,r))\geq cr^n,\ \mathcal L^n(B(x,r)\setminus E)\geq cr^n,$$
$$cr^{n-1}\leq\mu_E(B(x,r))\leq Cr^{n-1}.$$
\end{lm}

These ingredients can be used to prove the blow-up theorem:

\begin{thm}\label{periblow}
Let $x\in\partial^{\ast}E$ and $H(x)=\{y:y\cdot\nu_E(x)\leq 0\}$. Then $\chi_{r^{-1}(E-x)} \to \chi_{H(x)}$ as $r\to 0$ locally in $L^1(\Rn)$.
\end{thm}

To get this one first uses compactness to show that some subsequence converges to the characteristic function of a set $F$ with locally finite perimeter for which $\nu_F$ is constant $\mu_F$ almost everywhere. Moreover, the distributional derivatives of $\chi_F$ vanish in directions orthogonal to $\nu_F$, whereas the derivative in the direction of $\nu_F$ is non-zero. From this it is not trivial but not very difficult either to show that $F$ is a half-space. Moreover, it follows that $\nu_E(x)$ is the approximate normal of $E$ at $x\in\partial^{\ast}E$ in the following sense:
\begin{equation}\label{appnormal}
\lim_{r\to 0}r^{-n}\mathcal L^n({\{y\in E\cap B(x,r):(y-x)\cdot\nu_E(x)>0\}})=0,
\end{equation}
\begin{equation}\label{appnormal1}
\lim_{r\to 0}r^{-n}\mathcal L^n({\{y\in B(x,r)\setminus E:(y-x)\cdot\nu_E(x)<0\}}) =0.
\end{equation}

From the blow-up theorem one can proceed to show that $\{y:y\cdot\nu_E(x)=0\}$ is the approximate tangent plane of $\partial^{\ast}E$ at $x$, from which the rectifiability follows.

The \emph{essential boundary} of $E$, 
$$\partial_eE=\{x:\Theta^{\ast n}(E,x)>0\ \text{and}\ \Theta^{\ast n}(\Rn\setminus E,x)>0\}$$
gives a different view of the finite perimeter sets and the reduced boundary. We have, see, Theorems 4.5.6 and 4.5.11 in \cite{Fed69}, and also \cite[Theorem 16.2]{Mag12},

\begin{thm}
If $E\subset\Rn$ is a set of finite perimeter, then $\partial^{\ast}E\subset\partial_eE$ and $\H^{n-1}(\partial_eE\setminus\partial^{\ast}E)=0$.
\end{thm}

\begin{thm}
A measurable set $E\subset\Rn$ has finite perimeter if and only if $\H^{n-1}(\partial_eE)<\infty$.
\end{thm}

Lahti \cite{Lah20} and Eriksson-Bique \cite{Eri22} gave different proofs and metric space versions for the last theorem.

\subsection{Plateau type problems}
Sets of finite perimeter give a convenient setting to define and study codimension one generalized minimal surfaces. Other settings will be discussed later. The classical Plateau problem asks to find and describe the surface with minimal area among surfaces with a given boundary. Many variants of this, often of very general type, have been studied, some of them will be discussed later. Usually there are several non-trivial subproblems: what is surface, what is area, what is boundary?

In the case of finite perimeter sets boundary is not defined but instead one considers sets which agree with a given set outside a fixed set:

\begin{df}\label{perimin}
Let $A\subset\Rn$ and let $E_0,E\subset\Rn$ be sets of finite perimeter. We say that $E$ is \emph{perimeter minimizing} in $A$ with boundary data $E_0$ if  $E\setminus A = E_0\setminus A$ and $P(E)\leq P(F)$ for all sets $F$ of finite perimeter such that $F\setminus A = E_0\setminus A.$
\end{df}

The existence of perimeter minimizing sets follows by the usual direct method of calculus of variations: choose a minimizing sequence;
$P(E_i)\to\inf\{P(F): F\setminus A = E_0\setminus A\}$ and use compactness to select a converging subsequence, then the limit is a minimizer by lower semicontinuity:

\begin{thm}\label{periminex}
Let $A\subset\Rn$ be bounded and let $E_0\subset\Rn$ be a set of finite perimeter. Then there exists a perimeter minimizing set in $A$ with boundary data $E_0$.
\end{thm}

The real problem then is the regularity of the minimizers. I shall say something about this in Chapter \ref{Singularities}.

Sets of finite perimeter can be used to model many other geometric variational problems, too, see \cite{Mag12}.

\subsection{Functions of bounded variation}\label{BVfunctions}
Above we already defined these. The book \cite{AFP00} contains a lot of detailed information about them. Here I only discuss some properties related to rectifiability.

A function $u\in L^1(\Rn)$ has an \emph{approximate limit} $a$ at $x$ if
$$\lim_{r\to 0}r^{-n}\int_{B(x,r)}|u(y)-a|\,dy.$$
The set $S_u$ where $u$ does not have any approximate limit is called the approximate discontinuity set of $u$. At jump points the nature of the discontinuity is more specific: $x$ is called an \emph{approximate jump point} of $u$ if there are $a,b\in\R, a\neq b$, and $\nu\in S^{n-1}$ such that
$$\lim_{r\to 0}r^{-n}\int_{\{y\in B(x,r):(y-x)\cdot\nu>0\}}|u(y)-a|\,dy =\lim_{r\to 0}r^{-n}\int_{\{y\in B(x,r):(y-x)\cdot\nu<0\}}|u(y)-b|\,dy=0.$$
The set of approximate jump points is denoted by $J_u$. By \cite{Fed69} and \cite{Vol67} we have, see \cite[Theorem 3.78]{AFP00},

\begin{thm}\label{BVrect}
If $u\in BV(\Rn)$, then $S_u$ is $(n-1)$-rectifiable and $\H^{n-1}(S_u\setminus J_u)=0$.
\end{thm}

For the characteristic functions of sets of finite perimeter this follows from Theorem \ref{peridens} and \eqref{appnormal}, \eqref{appnormal1}. There is a coarea formula for BV-functions which implies that for almost all $t\in\R$ the sets $\{x:u(x)>t\}$ have finite perimeter. Letting $D$ be a countable dense set of such $t$, one can show that up to $\H^{n-1}$ measure zero $S_u$ is contained in $\cup_{t\in D}\partial^*\{x:u(x)>t\}$, from which the rectifiability of $S_u$ follows.

The book \cite{AFP00} contains much more about the structure of the derivative measure $Du$. First $|Du|\ll\H^{n-1}$, where $|Du|$ is the total variation measure of $Du$. For $u=\chi_E$ this follows from Theorem \ref{peristru} and after that for general $u$ from the coarea formula. Then, see the proof of \cite[Proposition 3.92]{AFP00}, 

\begin{thm}\label{BVrect1}
If $u\in BV(\Rn)$, then $Du\restrict\{x:\Theta^{\ast n-1}(|Du|,x)>0\}$ is $(n-1)$-rectifiable.
\end{thm}

A vector valued function $u:\Rn\to\R^k$ is in $BV(\Rn,\R^k)$ if its coordinate functions are of bounded variation. Then the derivative $Du$ is $k\times n$ matrix valued measure. Let $D_su$ denote its singular part in the Lebesgue decomposition. The following Alberti's rank one theorem, \cite{Alb93}, has applications in many areas:

\begin{thm}\label{Alberti}
If $u\in BV(\Rn,\R^k)$, then the Radon-Nikodym derivative\\ $D(D_su,|D_su|)(x)$ has rank 1 for $|D_su|$ almost all $x\in\Rn$.
\end{thm}

Massaccesi and Vittone \cite{MV19} have given a very simple and elegant proof using sets of finite perimeter, and De Philippis and Rindler derived it as a special case of their more general result in \cite{DR16}, see Section \ref{PDEsing}. But since Alberti's original proof gives a lot information about the structure, also related to rectifiability, of singular measures on $\Rn$, I shall briefly discuss it. As already indicated in Section \ref{Lebnull}, this is closely connected to the work of  Alberti, Cs\"ornyei and Preiss, \cite{ACP05} and \cite{ACP10}. They describe in \cite{ACP05} how to prove Theorem \ref{Alberti} by their tangent field results.

To prove Theorem \ref{Alberti} Alberti studied tangential properties of general measures. For $\mu\in\mathcal M(\Rn)$ and $x\in\Rn,$ let $E(\mu,x)$ be the set of vectors $v\in\Rn$ such that for some real valued $u\in BV(\Rn)$,
$$\lim_{r\to 0}\frac{|Du-v\mu|(B(x,r))}{\mu(B(x,r))}=0.$$
Then $E(\mu,x)$ is a linear subspace of $\Rn$. One can show  that if $\mu=\H^{n-1}\restrict E$, where $E$ is $\H^{n-1}$ measurable with $\H^{n-1}(E)<\infty$, then for $\mu$ almost all $x\in\Rn$, $E(\mu,x)= \apTan(E,x)^{\perp}$ if $E$ is $(n-1)$-rectifiable, and $E(\mu,x)= \{0\}$ if $E$ is purely $(n-1)$-unrectifiable. For $n=2$, $E(\mu,x)$ is closely related to the decomposition bundle $V(\mu,x)$ of Alberti and Marchese \cite{AM16}, which we discussed in Section \ref{Lebnull}.

For $v\in E(\mu,x), v\neq 0$,\ $\mu$ is kind of tangential to a derivative at $x$ of a BV-function. The key to the proof of Theorem \ref{Alberti} is that for any singular measure there is at most one direction  for which this can happen:

\begin{thm}\label{Alberti1}
If $\mu\in\mathcal M(\Rn)$ is singular, then the dimension of $E(\mu,x)$ is either 0 or 1 for $\mu$ almost all $x\in\Rn$.
\end{thm}

The proof of this theorem requires most of the effort. Alberti proved it first for $n=2$ with a calculus type argument. Then he used disintegration of $\mu$ with 2-dimensional slices. 

If $\mu\in\mathcal M(\Rn)$ is any singular Borel measure, we can decompose by Theorem \ref{Alberti1}
\begin{equation}\label{Alberti2}
\mu =\mu_1+\mu_0= \mu\restrict B_1 + \mu\restrict B_0\ \text{where}\ B_i=\{x:\dim E(\mu,x)=i\}, i=0,1.
\end{equation}

Alberti showed that $\dim E(\mu,x)>0$ for $\mu$ almost all $x\in\Rn$ if and only if $\mu=|Du|\restrict B$ for some $u\in BV(\Rn)$ and some Borel set $B\subset\Rn$. So $\mu_1=|Du|\restrict B_1$ for some $u\in BV(\Rn)$ and $\mu_0$ is orthogonal $Du$ for every $u\in BV(\Rn)$. From this it follows that for any singular measure $\mu$ and $u\in BV(\Rn,\R^k)$, the rank of the Radon-Nikodym derivative $D(Du,\mu)(x)$ is 1 for $\mu$ almost all $x\in\Rn$. Theorem \ref{Alberti} follows applying this to  $\mu=|D_su|$.

In the decomposition \eqref{Alberti2} $\mu_1$ has an Alberti representation (recall Definition \ref{rectrepr} and \eqref{Alberti3}) $\mu_1=\int\H^{n-1}\restrict E_t\,dt$ where each $E_t$ is $(n-1)$-rectifiable. 

A subclass SBV of BV, special functions of bounded variation, consists of functions $u\in BV(\Rn)$ for which $Du$ is a sum of an absolutely continuous and a rectifiable measure, the latter being concentrated on $J_u$. That is, $u\in SBV$ if the so-called Cantor part of $Du$ vanishes.  SBV and its many applications are extensively discussed in \cite{AFP00}

In \cite{Amb90} Ambrosio developed a theory of metric space valued functions of bounded variation. Let $X$ be  a locally compact metric space. A Borel function $u:\Rn\to X$ belongs to $BV(\Rn,X)$ if there is $\mu\in\mathcal M(\Rn)$ such that for every 1-Lipschitz function $\phi:X\to\R$,\ $\phi\circ u\in BV(\Rn)$ with $|D(\phi\circ u)|(A)\leq \mu(A)$ for $A\subset\Rn$. In particular, he proved an analogue of Theorem \ref{BVrect} in this setting.

Ambrosio, Coscia and Dal Maso investigated in \cite{ACD97} mappings $u:\Rn\to\Rn$ of bounded deformation, $BD(\Rn)$, that is, $Du+(Du)^T$ is a matrix valued Radon measure. Clearly, BV implies BD. Among other things, they proved the generalization of Theorem \ref{BVrect1}, the proof uses Besicovitch -- Federer projection theorem \ref{profed}:

\begin{thm}\label{BDrect}
If $u\in BD(\Rn)$, then $Du\restrict\{x:\Theta^{\ast n-1}(|Du|,x)>0\}$ is $(n-1)$-rectifiable.
\end{thm}

Various physically motivated PDE problems lead to functions of BV type for which analogous rectifiability results have been proven, see \cite{DO03}, \cite{DOW03}, \cite{DR03}, \cite{AKL02},  \cite{Mar20a}, \cite{Mar20b}.

\subsection{Perimeter in Heisenberg and Carnot groups}\label{perheis}
The initial motivation of Franchi, Serapioni and Serra Cassano in \cite{FSS01} to study rectifiability in Heisenberg groups was to develop De Giorgi's theory of sets of finite perimeter there. Recall the structure of the Heisenberg group $\hn$ and the notion of $(m,\h)$-rectifiable sets from Chapter \ref{heisenberg}, with $m=2n+1$ in the codimension one case.

Denote by $C^1_c(\hn,H\hn)$ the space of compactly supported continuous continuously Pansu differentiable functions with values in the horizontal sections of $\hn;\\ \phi(p) \in \tau_p(\h)$, where $\h=\{(z,t):t=0\}$ is the horizontal plane. Now the perimeter is defined in terms of the Heisenberg divergence: $\di_{H}\phi = \sum_{j=1}^n(X_j\phi_j + Y_j\phi_{n+j})$.

\begin{df}\label{perimeterheis}
The \emph{Heisenberg perimeter} of a Lebesgue measurable set $E\subset\hn$ is
$$P_H(E)=\sup\{\int_E\di_H\phi:\phi\in C^1_c(\hn,H\hn), |\phi|\leq 1\}.$$
If $P_H(E)<\infty$, we say that $E$ is a set of finite Heisenberg perimeter.
\end{df}

If $P_H(E)<\infty$ we again have by the Riesz representation theorem that there are $\mu_E\in\mathcal M(\hn)$ and a Borel function $\nu_E:\hn\to H\hn$  such that $|\nu_E(p)|=1$ for $\mu_E$ almost all $p\in\hn$ and 
\begin{equation*}
\int_E\di_H\phi = \int\phi\cdot\nu_E\,d\mu_E\ \text{for}\ \phi\in C^1_c(\hn,H\hn).
\end{equation*}

The reduced boundary $\partial^{\ast}E$ with $\mu_E(\hn\setminus \partial^{\ast}E)=0$ can then be defined as in the Euclidean case. 

This much is true in general Carnot groups. Franchi, Serapioni and Serra Cassano proved De Giorgi's structure theorem first in $\h^n$ in \cite{FSS01} and then in all step 2 Carnot groups in \cite{FSS03}. In particular, we have

\begin{thm}\label{peristruheis}
Let $E\subset\hn$ be a set of finite Heisenberg perimeter. Then $\partial^{\ast}E$ is $(2n+1,\h)$-rectifiable.
\end{thm}



The proof follows the same main lines as in the Euclidean case but it is technically much harder. Again, the blow-ups at the points of the reduced boundary converge to vertical subgroups. This statement is false in higher order groups, at least in the Engel group which is of step 3, and the analogue of De Giorgi's theorem is not known. Ambrosio, Kleiner and Le Donne \cite{AKL09} proved a partial result in general Carnot groups: some sequences of blow-ups  converge to vertical subgroups.

See the survey \cite{Ser16} of Serra Cassano for further comments and references.

\section{Currents and varifolds}\label{currents}

The classical reference for the overall theory of the currents in the Euclidean spaces is \cite{Fed69}. Currents and varifolds are very well presented also in \cite{Sim83}, \cite{LY02} and \cite{KP08}, and more informally in \cite{Mor88}. De Lellis's  survey \cite{Del21} enlightens the background and recent developments for this and next two chapters.

\subsection{Currents in Euclidean spaces}

Federer and Fleming begin their ground-breaking paper \cite{FF60} by the following quotation:

"Long has been the search for a satisfactory analytic and topological formulation of the concept '$k$ dimensional domain of integration in
 euclidean $n$-space.' Such a notion must partake of the smoothness of differentiable manifolds and of the combinatorial structure of polyhedral
 chains with integer coefficients. In order to be useful for the calculus of variations, the class of all domains must have certain compactness properties. All these requirements are met by the integral currents studied in this paper."
 
So the currents they introduced are generalized surfaces which, as they expected, have turned out to be extremely useful for the calculus of
 variations, and in many other topics, too. De Giorgi's theory of sets of finite perimeter already gave such a setting for codimension one surfaces. Currents can be of any dimension and they have many other advantages over sets of finite perimeter, but also some disadvantages.
 
Analytic theory of currents was developed by De Rham in the 1950s. They are just distributions over differential forms. Federer and Fleming introduced geometric aspects. The idea how they are related to smooth surfaces is simple. If $M$ is a smooth $m$-dimensional submanifold of $\Rn$, then one can integrate differential $m$-forms $\omega$ over it and define the linear functional $[M]$:
$$[M](\omega)=\int_M\omega.$$
By the Stokes theorem  $\int_{\partial M}\omega=\int_Md\omega$, where $d\omega$ is the exterior derivative of $\omega$. Thus if we define the boundary $\partial T$ of a general current $T$ by 
$$\partial T(\omega)=T(d\omega),$$
then $\partial [M]=[\partial M]$. 

The integral $\int_M\omega$ can be written as $\int_M\langle\omega(x),\xi(x)\rangle\,d\H^mx$, where $\xi(x)$ is an $m$-vector associated to the tangent plane of $M$ at $x$; differential $m$-forms can be defined as functions with values in the dual of the $m$-vectors. This means that $\xi(x)$ is of the form $v_1\wedge\dots\wedge v_m$ where $\{v_1,\dots,v_m\}$ is an orthonormal basis of the tangent $m$-plane of $M$ at $x$.

Let $\mathcal D^m(\Rn)$ be the space of  differential $m$-forms on $\Rn$ with compact support. They can be written as $\omega = \sum_{\a}\omega_{\a}dx^{\a}$, where $\a$ runs through the sequences $a(1)<\dots<\a(m), \a(i)\in \{1,\dots,n\}$, the $\omega_{\a}$ are smooth functions and $dx^{\a}=dx^{\a(1)}\wedge\dots\wedge dx^{\a(m)}$. 

By definition, an $m$-current $T, T\in\mathcal D_m(\Rn)$, is  a continuous linear functional on $\mathcal D^m(\Rn)$. The \emph{support} $\spt T$ of a current $T$ is the  smallest closed set such that $T(\omega)=0$ for every $\omega\in\mathcal D^m(\Rn)$ such that $\spt\omega\subset\Rn\setminus\spt T$. The \emph{mass}, generalizing area, is
$$M(T)=\sup\{|T(\omega)|:\|\omega\|_{\infty}\leq 1,\omega\in\mathcal D^m(\Rn)\}.$$ 
Differential forms can be pulled back, so currents can be pushed forward by maps (with proper conditions) $f; f_{\#}T(\omega)=T(f^{\#}\omega)$. 
We have the easy weak compactness theorem: if  $\sup_jM(T_j)<\infty, j=1,2,\dots,$, then there is a subsequence $(T_{j_i})$ and a current $T$ such that $T_{j_i}(\omega)\to T(\omega)$ for all $\omega\in\mathcal D^m(\Rn)$.

A consequence of the Riesz representation theorem is that if $M(T)<\infty$, then there is a Radon measure $\mu_T$ and an $m$-vector valued Borel function $\overrightarrow{T}$, the tangent vector field of $T$, such that
\begin{equation}\label{currriesz}
T(\omega)=\int \langle\omega(x),\overrightarrow{T}(x)\rangle\, d\mu_Tx\ \text{for}\ \omega\in\mathcal D^m(\Rn).\end{equation}
If $B\subset\Rn$ is a Borel set, then for $T$ as above we define the restriction of $T$ to $B, T\restrict B$, by integrating only over $B$, that is,  $\mu_{T\restrict B}=\mu_{T}\restrict B$. Similarly one defines $T\restrict f$ for functions $f$.

The currents with  $N(T):=M(T)+M(\partial T)<\infty$ are called \emph{normal}. Although our main interest is in $m$-currents with $m<n$, the $n$-currents in $\Rn$ too are interesting. In particular, the normal $n$-currents in $\Rn$ can be identified with BV-functions:

\begin{thm}\label{BVcurrent}
Let $T\in\mathcal D_n(\Rn)$ be normal. Then there is $g\in BV(\Rn)$ such that
$$T(fdx^{1}\wedge\dots\wedge dx^{n})=\int fg\,d \mathcal L^n\ \text{for}\ f\in C^{\infty}(\Rn).$$
\end{thm}

The proof can be done approximating $T$ by usual convolutions, cf. \cite[7.1.9]{LY02}. The converse also is true.

Based on this connection Federer presented almost the whole theory of BV-functions in \cite[Theorem 4.5.9]{Fed69} and its 31 statements.

Let $P_{\a}$ be the projection $P_{\a}(x_1,\dots,x_n)=(x_{\a(1)},\dots,x_{\a(m)})$. Then one checks for $T\in\mathcal D_m(\Rn)$ and $\omega= \sum_{\a}\omega_{\a}dx^{\a}\in\mathcal D^m(\Rn)$ that
$$T(\omega)=\sum_{\a}P_{\a\#}(T\restrict\omega_{\a})(dy^{1}\wedge\dots\wedge dy^m),$$ 
where $dy^{1},\dots,dy^m$ are the coordinate 1-forms on $\R^m$. With some extra work this yields
\begin{lm}\label{currproj}
If $T\in\mathcal D_m(\Rn), N(T)<\infty$ and $B$ is a Borel set with $\mathcal L^m(P_{\a}(B))=0$ for all $\a$, then $T\restrict B=0$. In particular this holds if  $\mathcal H^m(B)<\infty$ and $B$ is purely $m$-unrectifiable.
\end{lm}

The second statement follows from the Besicovitch-Federer projection theorem \ref{profed}, since that allows us to choose the appropriate coordinate axis.

Of the rich theory of currents I now only discuss rectifiable currents.

\begin{df}\label{currentrect}
An $m$-current $T$ in $\Rn$ with finite mass is called 
\emph{$m$-rectifiable} if there are an $m$-rectifiable $\H^m$ measurable set $E\subset\Rn$ and an $\H^m$ measurable positive  function $\theta$ on $E$ with $\int_E\theta\,d\H^m<\infty$ such that the values of $\overrightarrow{T}$ are simple $m$-vectors associated to the approximate tangent planes of $E$ and we have
\begin{equation}\label{currentrect1}
T(\omega)=\int_E\langle\omega(x),\overrightarrow{T}(x)\rangle\theta(x)\,d\H^mx\ \text{for}\ \omega\in\mathcal D^m(\Rn).\end{equation}  
If in addition the values of $\theta$ are integers  $T$ is called  \emph{integer multiplicity $m$-rectifiable current}. We say that $T$ is an \emph{integral current} if both $T$ and $\partial T$ are integer multiplicity rectifiable currents.\end{df}

We denote the set of $m$-rectifiable currents in $\Rn$ by  $R_m(\Rn)$, the set of integer multiplicity $m$-rectifiable currents in $\Rn$ by $\mathcal R_m(\Rn)$, and the set of integral $m$-currents in $\Rn$ by $\mathcal I_m(\Rn)$.

Again the condition on $\overrightarrow{T}$ means that for $\H^m$ almost all $x\in E$ the $m$-vector $\overrightarrow{T}(x)$ is of the form $v_1\wedge\dots\wedge v_m$ where $\{v_1,\dots,v_m\}$ is an orthonormal basis of the approximate tangent $m$-plane of $E$ at $x$. Then $\overrightarrow{T}(x)$ is uniquely determined up to sign. Choosing the sign means orienting $E$ and $T$.

The terminology and notation differ in different books and papers. 

Here are some of the main tools to study rectifiable currents:

The \emph{deformation theorem}: if $T$ is an integral current, then $T=P+Q+\partial S$, where $P$ is a polyhedral chain, $Q$ and $S$ are integral currents with small masses. The proof consists of carefully projecting $T$ into the skeletons of a cubical decomposition of $\Rn$. It gives a useful approximation of integral currents by polyhedral chains as well as the \emph{isoperimetric theorem}: if $T\in \mathcal I_{m-1}(\Rn)$ and $\partial T=0$, then there is $S\in \mathcal I_m(\Rn)$ such that $\partial S = T$ and $M(S)\lesssim M(T)^{m/(m-1)}$.

\emph{Slicing} is a very useful operation on currents. Let $T\in\mathcal D_m(\Rn)$ with $N(T)<\infty$ and  $f:\Rn\to\R^k, k\leq n,$ be a Lipschitz map. Then the slice of $T$ at $t\in\R^k$, $\langle T,f,t\rangle$, is an $(n-k)$-current with support in $\spt T\cap f^{-1}\{t\}$ such that $T$ is obtained as an integral of the $\langle T,f,t\rangle$. For simplicity I only consider slicing with real valued Lipschitz  maps $f:\Rn\to\R$. Then we can define
$$ \langle T,f,t\rangle=(\partial T)\restrict\{x:f(x)>t\}-\partial( T\restrict\{x:f(x)>t\}).$$
For almost all $t\in\R$, 
$$\spt\langle T,f,t\rangle \subset \spt T\cap f^{-1}\{t\},\ \partial\langle T,f,t\rangle=-\langle \partial T,f,t\rangle,\  
N(\langle T,f,t\rangle)<\infty\ \text{and}$$
$$\langle T,f,t\rangle\in \mathcal R_{m-1}(\Rn)\ \text{if}\ T\in\mathcal R_m(\Rn).$$
The first line is rather easy to prove, the second can be proven with the help of Theorem \ref{rectslice}, see \cite[Section 28]{Sim83}.

Here is a general, not very hard, rectifiability theorem, see \cite[Theorem 32.1]{Sim83}. 
Recall similar results in \ref{BVrect1} and \ref{BDrect}.

\begin{thm}\label{currentsrect} If $T\in\mathcal D_m(\Rn)$ is normal and $\Theta^{\ast m}(\mu_T,x)>0$ for $\mu_T$ almost all $x\in\Rn$, then $T\in R_m(\Rn)$.
\end{thm}

I say a few words about the proof. We should establish \eqref{currentrect1} starting with \eqref{currriesz}. By Theorem \ref{densmeas} $\Theta^{\ast m}(\mu_T,x)<\infty$ for $\H^m$ almost all $x\in\Rn$, and so also for $\mu_T$ almost all $x\in\Rn$ by Lemma \ref{currproj}. The set  
$$E=\{x:\Theta^{\ast m}(\mu_T,x)>0\}$$
has $\sigma$-finite $\H^m$ measure, and $\mu_T$ and $\H^m\restrict E$ are mutually absolutely continuous. Hence $\mu_T=\theta\H^m\restrict E$ for some $\theta$. 
By Lemma \ref{currproj} $\mu_T(B)=0$ for every purely $m$-rectifiable set $B\subset E$, whence $E$ is $m$-rectifiable. That $\overrightarrow{T}$ is associated to the approximate tangent planes of $E$ requires more work. It can be established  by a blow-up method, see \cite[Theorem 32.1]{Sim83}.

We shall discuss mass minimizing rectifiable currents in the Chapter \ref{Singularities}. For their existence we need the compactness theorem:

\begin{thm}\label{currentscomp} If $T_j\in\mathcal I_m(\Rn), j=1,2,\dots,$ with  $\sup_jN(T_j)<\infty$, then there is a subsequence $(T_{j_i})$ and a current $T\in\mathcal I_m(\Rn)$ such that $T_{j_i}(\omega)\to T(\omega)$ for all $\omega\in\mathcal D^m(\Rn)$.
\end{thm}

The main point here is that the limit current is rectifiable. By the lower semicontinuity of the mass, $N(T)<\infty$. To apply Theorem \ref{currentsrect} we need to know that the upper density of $\mu_T$ is positive, which roughly means that $T$ should not be scattered around a set of non-$\sigma$-finite $\H^m$ measure. Another thing that needs checking is that the weight function of $T$ is integer valued. The proof of the theorem is by induction on $m$, so we also ought to know to what $(m-1)$-dimensional currents we should apply the induction hypothesis. All this is dealt with in the following slicing criterion for rectifiablity: 

\begin{lm}\label{slicecrit}
If $T\in\mathcal D_m(\Rn)$ is normal, $\partial T=0$ and $\partial(T\restrict B(a,r))\in\mathcal R_{m-1}(\Rn)$ for every $a\in\Rn$ and for almost all $r\in (0,\infty)$, then $T\in\mathcal R_m(\Rn)$.
\end{lm}

To prove this, first the isoperimetric theorem with some covering arguments are used to prove that $\Theta^{\ast m}(\mu_T,x)>0$ for $\mu_T$ almost all $x\in\Rn$. Then by Theorem \ref{currentsrect} we know  that $T\in R_m(\Rn)$. That the weight function of $T$ is integer valued follows, for example, by a rather simple blow-up argument. Now we have reduced the problem one dimension lower and the compactness theorem follows by an induction argument. As an easy consequence it has the interesting boundary rectifiablity theorem:

\begin{thm}\label{currentsrect1} If $T\in\mathcal R_m(\Rn)$ with  $M(\partial T)<\infty$, then $\partial T\in\mathcal R_{m-1}(\Rn)$.
\end{thm}

The above essentially was a very rough sketch of the  original proof of Federer and Fleming in  \cite{FF60} (and also in \cite{Fed69}).  Later other proofs were given which avoided the use of the Besicovitch-Federer projection theorem \ref{profed},  
by Solomon \cite{Sol84} using multivalued functions, then by White \cite{Whi89} by more classical analysis, and later by Ambrosio and Kirchheim in \cite{AK00b} in the metric space setting relying on BV-functions, as we shall see in the next section. 

Above the coefficients of rectifiable currents have been real numbers or integers. Other coefficient groups also have been studied. In particular for  integers modulo $p$, where $p\geq 2$ is an integer, the same rectifiability and compactness theorems are valid. White characterized in \cite{Whi99} the normed coefficient groups for which they hold. For this he proved a rectifiability criterion  with 0-dimensional slices. A similar criterion in the metric space setting was proved independently by Ambrosio and Kirchheim, which we shall discuss in the next section.

The currents modulo $p$ exhibit many interesting new phenomena. They are extensively discussed in \cite{Del21}. Simple illustrative examples are presented in \cite[Section 11.1]{Mor88}. 

Instead of using mass to represent the area one can use size: for a rectifiable current $T$ as in \eqref{currentrect1} Size$(T)=\H^m(\{x\in E:\theta(x)\neq 0\})$. In some cases this is better than mass, but the existence of minimizers is harder to prove. See \cite{Dav19}.

\subsection{Currents in metric spaces}

Based on an idea of De Giorgi, Ambrosio and Kirchheim \cite{AK00b} developed the theory of currents in complete metric spaces. It might seem that currents, as linear forms on differential forms, would need a differential structure. But a differential form $\omega$ can be written as a linear combination of $f_0df_1 \wedge \dots \wedge df_m$
with smooth functions $f_i$. For the geometric theory of currents we could as well consider Lipschitz functions, and this would make sense in metric spaces.

Let $X$ be a complete metric space and let $\mathcal D^m(X)$ be the set of all $(m+1)$-tuples $(g,\pi_1,\dots,\pi_m)$ of real valued Lipschitz functions on $X$ with $g$ bounded. Then Ambrosio and Kirchheim defined an $m$-dimensional current on $X$ to be any multilinear, positively homogeneous and continuous (in a suitable weak sense) functional $T$ on $\mathcal D^m(X)$ such that $T(g,\pi_1,\dots,\pi_m)=0$ whenever some $\pi_i$ is constant on a neighbourhood of $\{g\neq 0\}$. In this simple setting they were able to develop and generalize the Euclidean theory to a surprising extent. In particular they proved compactness and boundary rectifiability theorems, results which seemed to rely heavily on Euclidean tools such as the deformation theorem. Thus this work gives much new insight also to the classical theory. I skip the rest of the definitions, but let us see how boundary is defined. First the exterior derivative of $\omega=(g,\pi_1,\dots,\pi_m)$ is $d\omega=(1,g,\pi_1,\dots,\pi_m)$ and then, as before, $\partial T(\omega)=T(d\omega)$.

The theory of rectifiable currents in metric spaces is based on the theory of rectifiable sets in \cite{AK00a}, recall Chapter \ref{metricspaces}. The proofs of the compactness and boundary rectifiability theorems are again by induction and  one of the main tools is slicing. But now the rectifiability criterion for $m$-dimensional currents is given in terms of the slices with Lipschitz maps $f:X\to\R^m$, which are zero-dimensional currents, that is, measures, and when rectifiable, sums of point masses. Another basic tool is provided by a theory of metric space valued BV-functions which Ambrosio developed in \cite{Amb90} and which was mentioned in Section \ref{BVfunctions}. A key fact is that  $y\mapsto\langle T,f,y\rangle, y\in\R^m,$ is a BV-function. In the Euclidean setting this was proved by Jerrard in \cite{Jer02}. Since the BV-functions are Lipschitz on large subsets, one can, essentially, conclude that if $T$ is a normal current, then the set of those $x\in\spt T$ which are atoms of $\langle T,f,f(x)\rangle$ is rectifiable. Using this one then characterizes rectifiable $m$-currents $T$ by the property that for every Lipschitz map $f:X\to\R^m$ for $\mathcal L^m$ almost all $y\in\R^m$  the slice $\langle T,f,y\rangle$ is a rectifiable $0$-current. The compactness and boundary rectifiability theorems then follow by induction arguments employing slicing with real valued functions.

As mentioned before this gives new proofs also in the Euclidean setting.

Lang developed in \cite{Lan11} another approach which applies to local currents, not necessarily having finite mass. In \cite{AK00b} it is essential that the currents have finite mass.

Currents in Heisenberg groups is a fairly complicated issue, even from the point of view of definitions, since they are based on a difficult (at least to me) concept of Rumin's complex. Their theory is much less developed than the Euclidean and metric theories, see \cite{Vit20}.

\subsection{Varifolds}\label{Varifolds}
Varifolds, like currents, are generalized surfaces, better in some aspects and worse in some. In particular, there is no concept of boundary and no need for orientation. In a way they are more general than currents; any current $T$ with finite mass induces a varifold via the formula \eqref{currriesz}. They were introduced by Almgren in the 1960s in unpublished notes and his little book \cite{Alm66}, and the basic results were presented by Allard in \cite{All72}.
They are discussed in \cite{Sim83} and \cite{LY02}.

For any $A\subset\Rn$ we set $G_m(A)=A\times G(n,m)$.

\begin{df}\label{varifold}
Any Radon measure on $G_m(\Rn)$ is called an $m$-\emph{varifold}. To each $\H^m$ measurable and $m$-rectifiable set $E$ and non-negative $\H^m$ measurable function $\theta$ on $E$ with $\int_E\theta\,d\H^m<\infty$ we associate the rectifiable $m$-varifold $v(E,\theta)$ defined by 
$$v(E,\theta)(B)=\int_{\{x\in E:(x,\apTan(E,x))\in B\}}\theta(x)\,d\H^mx,\ B\subset G_m(\Rn)\ \text{a Borel set}.$$ 
When $\theta=1$ we write $v(E)=v(E,\theta)$.
\end{df}

To any $m$-varifold $v$ we associate the Radon measure $\mu_v$ on $\Rn$ by $\mu_v(A)=v(G_m(A))$. The \emph{mass} of $v$ is $M(v)=\mu_v(\Rn)$.
The image $f_{\#}v$ of $v$ under a smooth map $f:\Rn\to\Rn$ is defined by
$$f_{\#}v(A)=\int_{F^{-1}(A)}J_Vf(x)\,dv(x,V),\ A\subset G_m(\Rn)\ \text{a Borel set},$$ 
where $F(x,V)=(f(x),df(x)(V))$ and $J_Vf$ is a Jacobian of $f$ along $V$. Then for a rectifiable $m$-varifold $v(E,\theta)$,\ $M(v(E,\theta))=\int_E\theta\,d\H^m$ and the image of $v(E,\theta)$ under a diffeomorphism $f:\Rn\to\Rn$ is the rectifiable $m$-varifold given by
$f_{\#}v(E,\theta)=v(f(E),\theta\circ f^{-1}).$

There is no natural boundary operator, but in order to study Plateau-type problems we can use the classical approach from the calculus of variations. Let $h_t:\Rn\to\Rn, t\geq 0$, be a one-parameter family of diffeomorphisms with $h_0$ the identity and with each $h_t$ the identity outside a fixed compact set $K$. Then a computation shows that
$$\frac{d}{dt}M(h_{t\#}(v\restrict G_m(K)))_{|t=0}=\int\di_{V}X(x)\,dv(x,V),$$
where $X(x)=\frac{d}{dt}h_t(x)_{|t=0}$ and $\di_{V}X(x)$ is the divergence along $V$. Motivated by this we define

\begin{df}\label{firstvar}
The \emph{first variation} $\delta v$ of a varifold $v$ is defined by
$$\delta v(X)=\int\di_{V}X(x)\,dv(V,x)$$
for any smooth vector-field $X:\Rn\to\Rn$ with compact support. If $\delta v(X)=0$ for all such $X$ with support in an open set $U$, then \ $v$ is called \emph{stationary} in $U$.
\end{df}


The first variation is a very interesting operator on vector fields for many reasons, not only that it defines stationary varifolds. If $v=v(M)$ is a rectifiable varifold over a smooth manifold $M$, then $\delta v$ can be expressed with the mean curvature of $M$, which leads to a concept of genealized mean curvature.

Allard  \cite{All72} proved the following rectifiability theorem.

\begin{thm}\label{Allard}
Suppose that the $m$-varifold $v$ satisfies $|\delta v(X)|\lesssim \|X\|_{\infty}$ for all smooth $X:\Rn\to\Rn$ with compact support, that is, (the total variation of) $\delta v$ is a Radon measure. If $\Theta^{\ast m}(\mu_v,x)>0$ for $\mu_v$ almost all $x\in \Rn$, then $v$ is rectifiable.
\end{thm}

The two main tools to prove this are the monotonicity formula and tangent cones. Both of these have analogues for currents and they are important in many ways.

The monotonicity formula for stationary varifolds, and analogously for mass minimizing currents, is the following: if $v$ is a stationary $m$-varifold in $U$, then for $x\in U$ and $0<r<s$ with $B(x,s)\subset U$,
$$s^{-m}\mu_v(B(x,s))-r^{-m}\mu_v(B(x,r)) = \int_{y\in B(x,s)\setminus B(x,r)}|y-x|^{-m-2}|P_{V^{\perp}}(y-x)|^2\,dv(y,V).$$
This is proved using test vector fields of the type $y\mapsto \varphi(|y-x|)(y-x)$ in the definition of the first variation. In particular, the finite density $\Theta^m(\mu_v,x)$ exists. Assuming its positivity, as in Theorem \ref{Allard}, we could conclude by Preiss's theorem \ref{preiss} the rectifiability of $\mu_v$. This was not available for Allard, but it would anyway be unnecessarily heavy machinery; the monotonicity formula gives much more than the mere existence of the density. A variant of the monotonicity formula also holds when $\delta v$ is Radon measure, not necessarily stationary, and gives the existence of density.

Tangent cones are important in particular for studying singularities of mass minimizing currents and stationary varifolds. Let, for example, $v$ be a  stationary $m$-varifold and, as in the case of tangent measures, $T_{x,r}(y)=(y-x)/r.$ We say that a varifold $C$ is a \emph{tangent cone} of $v$ at $x$ if $T_{0,r\#}C=C$ for $r>0$, that is, $C$ is a cone with vertex at $0$, and there is a sequence $r_i>0, \lim_{i\to\infty}r_i=0$, such that $\lim_{i\to\infty}T_{x,r_i\#}v=C$. The monotonicity of the density ratios together with an easy compactness theorem imply that such limits exist, and further, still based on the monotonicity formula, they are cones. But can there be more than one tangent cone at a point? The uniqueness of tangent cones is a central open problem and known only in some cases. White proved it for 2-dimensional currents in \cite{Whi83}, this paper also gives references to other cases. In \cite{Whi92} he constructed a counterexample for harmonic maps. See also \cite{Min18} and \cite{Del21}.

I now briefly explain how the monotonicity formula and tangent cones can be used to prove Theorem \ref{Allard} when $v$ is stationary. The same ideas work in the general case. As for rectifiable sets and measures, it suffices to show that for $\mu_v$ almost all $x$ the varifold $v$ has a unique tangent cone which is an $m$-plane. Let $C$ be some tangent cone at a typical point $x$. Then $C$ is stationary and $0\in\spt\mu_C$. As for tangent measures, $\mu_C$ is an $m$-uniform measure:
\begin{equation}\label{allardeq}
\mu_C(B(y,r))=\a(m)\Theta^m(\mu_v,x)r^m\ \text{for}\ y\in\spt\mu_C, r>0.
\end{equation}
From this we see by the monotonicity formula that $P_{V^{\perp}}(y-x)=0$ for $C$ almost all $(y,V)$. This implies that when $(0,V)\in\spt C$, then $\spt\mu_C\subset V=V_C$, and further, by stationarity and a constancy theorem, that $C=\Theta^m(\mu_v,x)v(V_C)$. So we have left to show that $V_C$ is unique, but this is now fairly easy: a general differentiation theorem implies that for any continuous function $\varphi$ on $G(n,m)$ and $v$ almost all $(x,V)$ the limit $\lim_{r\to 0}\int_{G_m(B(x,r))}\varphi(V)\,dv(x,V)/\mu_v(B(x,r))$ exists. But using the definition of the tangent varifold one quickly checks that this equals $\varphi(V_C)$; go to zero through the sequence defining $C$. This implies that $V_C$ is uniqe.

Allard's main result was a regularity theorem for stationary varifolds, and more generally for varifolds with $L^p$ conditions for the generalized mean curvature, see \cite{Sim83} and \cite{LY02}.

De Philippis, De Rosa and Ghiraldin extended Allard's rectifiability theorem in \cite{DDG18} replacing mass by more general integrands. Their proof was different. Instead of the monotonicity formula they used tangent measures and results of Preiss from \cite{Pre87}. 
Still another proof, based on PDE operators, and a more general result is provided by \cite{ADHR19}, see Section \ref{PDEsing}.

In \cite{AS97} Ambrosio and Soner introduced generalized varifolds and applied them to gradient flows; $G(n,m)$ considered as a class of matrices is replaced by a larger subclass of symmetric matrices. They proved a rectifiability theorem, but relying on Allard's theorem. Brakke \cite{Bra78} developed mean curvature flow with varifolds. We shall come to it in Section \ref{Flow}.

Moser proved a general result in \cite{Mos03} related to the above as well as to harmonic maps and Yang-Mills connections discussed in Chapter \ref{Singularities}. 

\section{Minimizers and quasiminimizers}

\subsection{Quasiminimizers}\label{Quasimin}
In addition to currents and varifolds there are several other ways to model minimal surfaces and related objects, see \cite{Dav19} and \cite{Del21}. Quasiminimizers provide a very natural and general setting for many variational problems. Let $E\subset\Rn$ be closed and unbounded and such that, for a fixed positive integer $m$,\ $0<\H^m(E\cap B(x,r))<\infty$ for $x\in E, r>0$. We say that $E$ is an \emph{$m$-quasiminimizer} if for some $M<\infty$,
$$\H^m(E\cap W)\leq M\H^m(f(E\cap W))$$
for all Lipschitz mappings $f: \Rn \to \Rn$ such that $W=\{x: f(x) \neq x\}$ is bounded. 
If this holds with $M = 1$, then $E$ minimizes $m$-dimensional Hausdorff measure.
The setting in the papers quoted below is more general. In particular, there is also a local, often very useful, version, but I skip it here. The quasiminimizers were introduced by Almgren in \cite{Alm76} under the name restricted sets. He proved that they are AD-$m$-regular and $m$-rectifiable. David and Semmes investigated them in \cite{DS00}. They re-proved  Almgren's results and went further. The following is a special case of their results:
\begin{thm}\label{quasimin}
If $E\subset\Rn$ is a closed $m$-quasiminimizer, then $E$ is AD-$m$-regular, uniformly $m$-rectifiable and it contains big pieces of Lipschitz graphs (recall Section \ref{Lipmap-appr}).
\end{thm}

Both Almgren's and David-Semmes's proofs use Lipschitz projections into $k$-dimensional cubical skeleta like in the Federer-Fleming proof of the deformation theorem of currents. First this gives AD-regularity. Then, by David and Semmes, via many complicated constructions the big pieces of Lipschitz graphs condition is verified. 

The codimension 1 case was studied by different methods in \cite{DS98} and \cite{JKV97}. All these papers contain many interesting results on and connections with various geometric variational problems. 

There is much later work along these lines, see David's long paper \cite{Dav19} for a very general setting, for discussion and references. It seems to give the most general rectifiability results. In particular, he there used sliding conditions; the deformations were required to preserve given boundary pieces but were allowed to slide along them.



\subsection{Mumford-Shah functional}\label{Mumsha}

Let $\Omega\subset\Rn$ be a domain and $g$ a bounded measurable function in $\Omega$. The \emph{Mumford-Shah functional} $J$ is then defined by
$$J(u,K)=\int_{\Omega\setminus K}(u-g)^2+\H^{n-1}(K)+\int_{\Omega\setminus K}|\nabla u|^2$$
for
$$(u,K)\in\mathcal A(\Omega):=\{(u,K):K\subset\Omega\ \text{relatively closed and}\ u\in W^{1,2}_{loc}(\Omega\setminus K)\}.$$
We assume that there are $(u,K)\in\mathcal A(\Omega)$ with $J(u,K)<\infty$, which is always true if $\Omega\subset\Rn$ is bounded. 
For many aspects of the Mumford-Shah functional, including applications to image segmentation and conjectures and results on minimizers, see the books \cite{AFP00} and \cite{Dav05}. Here I restrict to things related to rectifiability.

A minimizer for $J$ is a pair $(u,K)\in\mathcal A(\Omega)$ which gives the smallest value for $J$. They always exist, although it is far from obvious since Hausdorff measure is not lower semicontinuous. One way to prove the existence is to minimize first
$$\int_{\Omega}(u-g)^2+\H^{n-1}(S_u)+\int_{\Omega}|\nabla u|^2$$
for $u\in SBV(\Omega)$, recall Section \ref{BVfunctions}. Minimizers for this exist by the compactness properties of SBV, but the problem that $S_u$ need not be closed has to be dealt with. Here one cannot use the full $BV(\Omega)$, since it would give 0 for the infimum. Anyway, now $S_u$ is $(n-1)$-rectifiable by Theorem \ref{BVrect}. This approach is discussed in \cite{AFP00}. In \cite{Dav05} a different approach without SBV is explained.

For a minimizer $(u,K)$, $u$ is  in $C^1(\Omega\setminus K)$, which follows from the fact that it solves the PDE $\Delta u = u-g$. For $K$ there are conjectures which are only partially solved. David and Semmes proved the following in \cite{DS96}, see also \cite{Dav05}:

\begin{thm}\label{MS}
If $(u,K)$ is a minimizer for $J$ and $B(x,2r)\subset\Omega$, then $K\cap B(x,r)$ is contained in an AD-$(n-1)$-regular uniformly $(n-1)$-rectifiable set.
\end{thm}

The key to the proof is that the failure of the Poincar\'e inequality in the complement of an AD-$(n-1)$-regular set $E$ at most scales implies uniform rectifiability of $E$. This is understandable, because the validity of the Poincar\'e inequality requires that $E$ does not separate the space too much. More precisely: $E$ is uniformly $(n-1)$-rectifiable, if there exists a 
 positive number $c$ such for all $M\geq 1$ the set $F(E,c,M)$ of pairs $(x,r), x\in E, 0<r<d(E)$,  satisfying the following condition is a Carleson set: for all balls $B(x_i,r_i)\subset B(x,r)\setminus E, i=1,2,$ with $r_i>cr$ and for all 
$f\in W^{1,1}(B(x,Mr)\setminus E)$,
\begin{equation}\label{MS1}
\left|r_1^{-n}\int_{B(x_1,r_1)}f-r_2^{-n}\int_{B(x_2,r_2)}f\right| \leq Mr^{1-n}\int_{B(Mx,r)\setminus E}|\nabla f|.\end{equation}
David and Semmes proved this 
by showing that this condition implies the local symmetry of Theorem \ref{LSthm}. Another proof is described in \cite{Dav05}. The converse  is false; an example is the complement of the union of the balls with radius $1/10$ centered in the integer lattice of a coordinate hyperplane.

For slight simplicity, assume $\Omega=\Rn$. To prove that for a minimizer $(u,K)$ the set $F(K,c,M)$ is a Carleson set, one applies \eqref{MS1} with $u=f$ and constructs a competitor to get for some $p<2$, 
\begin{equation*}\label{MS2}
\omega_p(x,Mr)=r^{p/2-n}\int_{B(x,Mr)\setminus K}|\nabla u|^p>\e(M)>0.\end{equation*}
As $r^{1-n}\int_{B(x,r)\setminus K}|\nabla u|^2$ is bounded, it is not very difficult to prove that the set of $(x,r)$ such that $\omega_p(x,r)>\e$ satisfies a Carleson condition, from which it follows that $F(K,c,M)$ is a Carleson set.

Theorem \ref{MS} holds for a much larger class of quasiminimizers.

\subsection{Some free boundary problems}\label{Freebdry}

In \cite{DET18} David, Engelstein and Toro studied the following two-phase free boundary problem: Let $\Omega\subset\Rn$ be a bounded domain and $q_+$ and $q_-$ bounded continuous functions on $\Omega$. Let
$$J(u)=\int_{\Omega}(|\nabla u(x)|^2+q_+(x)^2\chi_{\{u>0\}}(x)+q_-(x)^2\chi_{\{u<0\}}(x))\,dx.$$
Among other things they proved that if $u$ is almost minimizer (I omit the definition) for $J$, then, under slight extra conditions, the sets $\Omega\cap\partial\{x\in \Omega:u(x)>0\}$ and $\Omega\cap\partial\{x\in \Omega:u(x)<0\}$ are locally AD-$(n-1)$-regular and uniformly $(n-1)$-rectifiable. The proof is a complicated mixture of potential theory and geometric measure theory. In particular, proving the AD-regularity is quite demanding and achieved with  estimates for the harmonic measure.

We shall return to the corresponding one-phase problem in Section \ref{Freebdry1}.

Rigot \cite{Rig00} proved the uniform rectifiability of sets almost minimizing  perimeter, recall Section \ref{Finper}. Let $g(0,\infty)\to (0,\infty)$ with $g(x)=o(x^{(n-1)/n})$. 

\begin{thm}
Let $E\subset\Rn$ be Lebesgue measurable. If 
$$P(E)\leq P(F)+g(\mathcal L^n((E\setminus F)\cup (F\setminus E))$$
whenever $F\subset\Rn$ is Lebesgue measurable and $F=E$ outside some compact set, then $E$ is equivalent to $E'$ for which $\partial E'$ is AD-$(n-1)$-regular and uniformly $(n-1)$-rectifiable.
\end{thm}

She proved this by showing that $\partial E'$ is a Semmes surface, recall Section \ref{UR}.

\section{Rectifiability of singularities}\label{Singularities}
There is a huge literature on singularities of solutions to geometric variational, and other, problems. Here I only briefly present some results and ideas related to rectifiability. I first discuss currents and varifolds and then harmonic maps. The methods for both are rather similar and I give a few more details in the latter case. These topics have pretty much in common with the other topics I then discuss.

\subsection{Mass minimizing currents and stationary varifolds}\label{mincur}
De Lellis's survey \cite{Del21} gives an excellent up-to-date view of this wide topic. Currents give a very convenient setting for the Plateau problem for orientable surfaces. Let $B\in\mathcal R_{m-1}(\Rn)$ with $\partial B=0$ and  with finite mass. Then there are currents $T\in\mathcal I_m(\Rn)$ with $\partial T=B$, for example cones over $B$. We say that such a $T$ is \emph{mass minimizing} if $M(T)\leq M(S)$ for all $S\in\mathcal I_m(\Rn)$ with $\partial S = \partial T$. For the important local and homological minimizers, see \cite{Fed69}, \cite{Sim83} or \cite{LY02}. The existence of mass minimizing currents follows from the compactness theorem \ref{currentscomp} together with the easy facts that mass is lower semicontinuous  and the boundary operator is continuous. How much can we say about their regularity?

\begin{df}\label{singset}
Let $T\in\mathcal D_m(\Rn)$. A point $x\in\spt T\setminus \spt\partial T$ is called a \emph{regular} point of $T$ if it has a neighbourhood $U$ such that $\spt T\cap U$ is an $m$-dimensional smooth  submanifold of $\Rn$. Otherwise $x$ is called a \emph{singular} point of $T$. The set of singular points of $T$ is denoted by $\sing(T)$.
\end{df}


\begin{thm}\label{singthm}
Let $T\in\mathcal I_m(\Rn)$ be mass minimizing.
\begin{itemize}
\item[(1)] If $m=1$ or $m=n-1\leq 6$, then $\sing(T)=\emptyset$.
\end{itemize}
\begin{itemize}
\item[(2)] If  $m=n-1\geq 7$, then $\dim\sing(T)\leq m-7$.
\end{itemize}
\begin{itemize}
\item[(3)] If  $m\geq 2$, then $\dim\sing(T)\leq m-2$.
\end{itemize}
The bounds in (2) and (3) are sharp.
\end{thm}

(1) was proved by Simons in \cite{Sim68}, the cases $m=2$, by Fleming, and $m=3$, by Almgren, were done earlier. An example showing the sharpness in (2) in $\R^8$ is the current induced by the cone $\{x\in\R^8:x_1^2+\dots +x_4^2=x_5^2+\dots +x_8^2\}$. It has a singularity at the origin and it was shown to be mass minimizing (locally) by Bombieri, De Giorgi and Giusti in \cite{BDG69}, see also \cite{Giu84}. Cartesian products with $\R^{n-8}$ give higher dimensional examples. Examples for (3) are obtained by complex analytic varieties, which too minimize area. For example, $\{(w,z)\in \C^2: w^2=z^3\}$ has a genuine branch point at the origin. The estimate (2) was proved by Federer in \cite{Fed70} and (3) by Almgren in his massive more than 1000 pages long paper, which appeared as a Princeton University preprint in the early 1980s and was published in \cite{Alm00}. One of the key tools, due to De Giorgi and Reifenberg, in codimension 1 regularity theory is the approximation by graphs of harmonic functions. In higher codimensions this is not possible, for example because of the complex analytic varieties. To overcome this Almgren developed a theory of multivalued functions minimizing Dirichlet integral. Much of Almgren's work has been simplified and extended by De Lellis and Spadaro, see the surveys \cite{Del16a}, \cite{Del21} and the references given there.

So the singular sets are small but their structure is a big open question. As above, in all known examples they are quite nice, but in general it is not known if they could be some kind of fractals. However, see the comment at the end of this section. Simon proved the following in \cite{Sim95b}:

\begin{thm}\label{singthm1}
If $n\geq 8$ and $T\in\mathcal I_{n-1}(\Rn)$ is mass minimizing, then $\sing(T)$  is $(n-8)$-rectifiable.
\end{thm}

The proof has similar ingredients as the proof of Simon's  Theorem \ref{enminthm} for harmonic maps, which we shall discuss soon. In particular, monotonicity formula, tangent cones and the general form of Theorem \ref{simrect} play decisive roles.

For discussion and results on the structure of the singular sets of minimizing currents modulo $p$, see \cite{DHMS20} and \cite{DHMSS22}.

Much less is known about the singular sets of stationary $m$-varifolds, not even if they have $\H^m$ measure zero. The best known result still is due to Allard \cite{All72} saying that the regular set is an open dense subset of the support. However, the stratification can be used to obtain relevant information about the singularities, both for currents and varifolds. Let $v$ be an integer multiplicity rectifiable $m$-varifold in $\Rn$. A varifold cone $C$ is said to be $k$-symmetric if there exists $V\in G(n,k)$ such that $C$ is $V\times C'$ for some cone $C'$. Then the $k$th stratum of $v$ is 
$$S^k(v)=\{x\in \spt v: \text{no tangent cone of}\ v\ \text{at}\ x\ \text{is}\ (k+1)-\text{symmetric}\}.$$
So $S^0(v)\subset\dots\subset S^{m-1}(v)\subset\sing(v)$. Almgren \cite{Alm00} had proved in the 1980s that $\dim S^k(v)\leq k$. Naber and Valtorta \cite{NV20} proved much more:

\begin{thm}\label{NV}
Let $v$ be an integer multiplicity rectifiable $m$-varifold in $\Rn$ whose first variation is a Radon measure. Then for $k=0,\dots,m-1$,\ $S^k(v)$ is $k$-rectifiable. Moreover, for $\H^k$ almost all $x\in S^k(v)$ there exists a unique $V\in G(n,k)$ such that every tangent cone of $v$ at $x$ is of the form $V\times C$ for some cone $C$.
\end{thm}

Note that the last statement does not mean that the tangent cones at $x$ would be unique. 

For mass minimizing integer multiplicity rectifiable $(n-1)$-currents $T$ in $\Rn$ Simon proved, in addition to Theorem \ref{singthm1}, that the whole singular set agrees with the top stratum: $\sing(T)=S^{n-8}(T)$. We shall see other results like this later in this chapter. For mass minimizing currents of codimension bigger than 1 $S^{m-1}(T)$ need not be the whole singular set. For example, for 
$\{(w,z)\in \C^2: w^2=z^3\}$ the origin is a singular point but there is a unique 2-plane tangent cone. However, with multiplicity $2 > 1$, which is the reason why this can happen. 

As for harmonic maps (see below) a key feature in Naber and Valtorta's method is to consider quantitative stratifications; no tangent cone in a ball is $\e$ close to a symmetric cone. The proof is rather similar to that of Theorem \ref{enminthm1}.  Their method also gives a new proof for Theorem \ref{singthm1}.

If the Euclidean metric is perturbed to a suitable $C^{\infty}$ metric $d$, the singular sets can be quite wild, as shown by Simon in \cite{Sim21a} and \cite{Sim21b}: if $N\geq 7, l\geq 1$ and $K\subset \R^l$ is compact, he then constructed minimal hypersurfaces in $(\R^{N+1+l},d)$ with singular set $\{0\}\times K$. The very complicated and technical proof is mainly based on PDE methods; singular solutions of the symmetric minimal surface equation, see \cite{FS20}, are the building blocks of the construction.

\subsection{Energy minimizing maps}\label{Harmmaps}

Let $\Omega$ be an open subset of $\Rn$ and $N$ a Riemannian submanifold of $\R^p$. 

\begin{df}\label{enmin}
A map $u:\Omega\to N$ in the Sobolev space $W^{1,2}(\Omega,N)$ (that is, $u\in W^{1,2}(\Omega,\R^p)$ and $u(x)\in N$ for almost all $x\in\Omega)$ is \emph{energy minimizing} if for every ball $B(x,r)\subset\Omega$,
\begin{equation}\label{enmineq}
\int_{B(x,r)}|Du|^2\leq \int_{B(x,r)}|Dw|^2
\end{equation}
for all $w\in W^{1,2}(\Omega,N)$ such that $w=u$ in a neighbourhood of $\partial B(x,r)$.

The \emph{regular set}, $\reg(u)$, of $u$ is the set of points $x\in\Omega$ such that $u$ is $C^{\infty}$ in some neighbourhood of $x$. The \emph{singular set}, $\sing(u)$, of $u$ is $\Omega\setminus$Reg($u$).
\end{df}
Then $x\in\reg(u)$ if $u$ is continuous in some neighbourhood of $x$, cf. \cite[p.309]{SU82}.

Energy minimizing maps satisfy a Laplace-type equation generated by the restriction that the target is $N$. If $n\leq 2$ these maps have no singularities, but already when $n=3$ and $N$ is a two-sphere they may occur: $u:\R^3\to S^2, u(x)=x/|x|$, is energy minimizing and it has a singularity at the origin. Then $u_n:\R^n\to S^2, n>3, u_n(x)=u(x_1,x_2,x_3),$ is energy minimizing in $\Rn$ with an $(n-3)$-dimensional singular set. This is sharp: Schoen and Uhlenbeck \cite{SU82} showed that the Hausdorff dimension of $\sing(u)$ is at most $n-3$. Simon proved in \cite{Sim95a}, see also \cite{Sim96}, the following much stronger result:

\begin{thm}\label{enminthm}
If $n\geq 3$, $N$ is real analytic and  $u:\Omega\to N$ is energy minimizing, then $\sing(u)$ is $(n-3)$-rectifiable.
\end{thm}
 
The proof is ingenious and very complicated. I try to give some flavour of it. The book \cite{Sim96} gives a detailed exposition, also of the background material. In the end the rectifiability is obtained by showing that essential parts of the singular set satisfy the assumptions of Theorem \ref{simrect}, that is, the $\e$ approximation by $(n-3)$-planes and the gap condition. Or rather, the assumptions of a generalization of Theorem \ref{simrect} where additional exceptional sets are allowed. But it is a long way to get there. 

Let $u:\Omega\to N$ be energy minimizing. For much of what is said below the real analyticity of $N$ is not needed and Simon proved some partial results in the general case. There is some PDE theory involved, which, in particular, gives that there is $\e(n,N)>0$ such that if $u:\Omega\to N$ is energy minimizing, $B(x,r)\subset\Omega$ and $r^{2-n}\int_{B(x,r)}|Du|^2<\e(n,N)$, then $x\in\reg(u)$. This gives almost immediately that the singular set of $u$ has locally finite $\H^{n-2}$ measure and with a little more work one finds that $\H^{n-2}(\sing(u))=0$. But the bound $n-3$ requires  further effort.

Two basic tools are analogous to those we met in Allard's rectifiablity theorem for varifolds; the monotonicity formula and  tangent maps. The monotonicity formula now takes the form
\begin{equation}\label{harmmon}
s^{2-n}\int_{B(x,s)}|Du|^2 - r^{2-n}\int_{B(x,r)}|Du|^2 = 2\int_{B(x,s)\setminus B(x,r)}|y-x|^{2-n}|\partial_{\tfrac{y-x}{|y-x|}}u(y)|^2\,dy
\end{equation}
when $0<r<s$ and $B(x,s)\subset\Omega$. In particular, the density
\begin{equation}\label{harmdens}
\Theta_u(x)=\lim_{r\to 0}r^{2-n}\int_{B(x,r)}|Du|^2
\end{equation}
exists. The proof of the monotonicity formula is based on the variational equation which $u$ satisfies.

A useful property of the density is that it is upper semicontinuous. We have also that $x\in \reg(u)$ if and only if $\Theta_u(x)=0$. One direction follows from the above $\e(n,N)$-property and the other is trivial.

Set $u_{x,r}(y)=u(x+ry)$ when $B(x,r)\subset\Omega$. A map $\varphi:\R^n\to N$ is called a \emph{tangent map} of $u$ at $x$ if there is a sequence $r_i>0$ tending to $0$ such that $u_{x,r_i} \to \varphi$ locally in $W^{1,2}(\Rn)$. The boundedness of the density ratios, due to the monotonicity formula, and a compactness theorem yield that tangent maps always exist. White showed in \cite{Whi92} that in general they need not be unique, but it is an open question whether they are unique when the target $N$ is real analytic. See \cite{Min18} for a discussion on this and related issues.

The tangent maps are energy minimizing and they have other very useful properties. First,
\begin{equation*}
\Theta_u(x)=\Theta_{\varphi}(0)=r^{2-n}\int_{B(0,r)}|D\varphi|^2\ \text{for all}\ r>0.
\end{equation*}
From this one concludes with the monotonicity formula  for $\varphi$ (since the left, and hence also the right, hand side now is 0) that 
\begin{equation}\label{HCM}
\varphi(\lambda x)=\varphi(x)\ \text{for all}\ x\in\Rn, \lambda>0,
\end{equation}
and that $x\in \reg(u)$ if and only if there is a constant tangent map at $x$. A little more calculus gives that $\Theta_{\varphi}(x)\leq \Theta_{\varphi}(0)$ for all $x\in\Rn$ and that
$$S(\varphi):=\{x\in\Rn:\Theta_{\varphi}(x)= \Theta_{\varphi}(0)\}$$
is a linear subspace of $\Rn$ such that $\varphi(x+y)=\varphi(x)$ for $x\in\Rn, y\in S(\varphi)$. Moreover, $S(\varphi)=\Rn$ if and only if $\varphi$ is constant, otherwise $S(\varphi)\subset\sing(\varphi)$. As $\H^{n-2}(\sing(\varphi))=0$ we have $\dim S(\varphi)\leq n-3$ for any non-constant tangent map $\varphi$.

Let $x\in \sing(u)$ and $\delta>0$.  Then there is $0<\e<\Theta_u(x)$ such that for every $0<r<\e$ there is an $(n-3)$-dimensional affine plane $V$ for which
\begin{equation}\label{harmeq}
\{y\in B(x,r): \Theta_u(y)\geq\Theta_u(x)-\e\}\subset \{y: d(y,V)\leq \delta r\}.
\end{equation}
Recall that for any tangent map $\varphi$ the set $S(\varphi)$ is contained in an $(n-3)$-plane. \eqref{harmeq} is not immediate from this, but not very difficult either. The upper semicontinuity of the density plays a role here.



\eqref{harmeq} gives an approximation of the singular set with $(n-3)$-planes and one can deduce from it that $\dim \sing(u) \leq n-3$. But to conclude  rectifiability  with something like Theorem \ref{simrect} we would still need the gap condition. That is the hardest part of the proof. Roughly speaking (very roughly), suppose we start with \eqref{harmeq} at some $x$ at some scale $r$ such that there are no gaps in a range of scales below $r$. Then many technical integral estimates on the derivatives of $u$ imply that at these scales $u$ is $L^2$ close to a map $\varphi$ as in \eqref{HCM}, which gives the required approximation at smaller scales with planes parallel to $V$.

Naber and Valtorta proved in \cite{NV17} with different methods further deep results on the structure of the singular sets, recall their similar results for varifolds from the previous section. Define 
$$S^k(u)=\{x\in \sing(u):\dim S(\varphi)\leq k\ \text{for all tangent maps}\ \varphi\ \text{of}\ u\ \text{at}\ x\}.$$
Then we have the stratification of the singular set
$$S^0(u)\subset S^1(u) \subset\dots\subset S^{n-3}(u)=S^{n-2}(u)=S^{n-1}(u)=\sing(u).$$
The equalities follow from the fact, essentially mentioned above, that $\dim S(\varphi)\leq n-3$ if and only $x\in\sing(u)$. 

By what was said about $S(\varphi)$, $x\in S^k(u)$ means that no tangent map $\varphi$ of $u$ at $x$ is $(k+1)$-symmetric. A tangent map $\varphi$ is $k$-symmetric if there is $V(\varphi)\in G(n,k)$ such that $\varphi(x+y)=\varphi(x)$ for $x\in\Rn, y\in V(\varphi)$. 

Schoen and Uhlenbck  proved in \cite{SU82} that $\dim S^k(u)\leq k$; similar arguments that gave $\dim \sing(u)\leq n-3$ apply. Naber and Valtorta proved in \cite{NV17}, without assuming real analyticity,

\begin{thm}\label{enminthm1}
If  $u:\Omega\to N$ is energy minimizing, then $S^k(u)$ is $k$-rectifiable for $k=0,1,,\dots,n-3$.
\end{thm}

This gives a new proof for Simon's theorem \ref{enminthm}. They also prove the result more generally for stationary maps, and in \cite{NV18} for a larger class of approximate harmonic maps. The latter paper has some simplifications of the proof of Theorem \ref{enminthm1}. In addition, these papers contain many other significant results. A key feature in their method is to consider quantitative stratifications; no tangent map in a ball is $\e$ close to a symmetric map. For these they proved volume estimates for the $r$-neighbourhoods of $S^k(u)$ which in particular yield the bound $k$ for the Minkowski dimension. They also showed that the plane $V(\varphi)$, mentioned above, is independent of $\varphi$. 

The density ratios are again one of the key factors. Setting\\ 
$\Theta_u(x,r)=r^{2-n}\int_{B(x,r)}|Du|^2$ it is shown for any finite Borel measure $\mu$ that if $u$ is not $\e$ close to any $(k+1)$-symmetric map in $B(0,8)$, then
\begin{equation}\label{NVeq}
\inf_{V\ k-\text{plane}} \int d(x,V)^2\,d\mu x \lesssim \int(\Theta_u(x,8)-\Theta_u(x,1))\,d\mu x,\end{equation}
Scaled versions of this are applied with $\mu$ a discrete approximation of the Hausdorff measure. They give volume estimates. Based on \eqref{NVeq} rectifiability is derived from Reifenberg type theorems, recall Section \ref{Reiftype} and Theorem \ref{ENVreif}. The proofs involve complicated covering and induction arguments.

Defect measures provide one of the tools: by compactness, any bounded sequence $(u_i)$ in $W^{1,2}$ has a subsequence  $(u_{i_j})$ converging weakly in $W^{1,2}$ to $u$ such that $|\nabla u_{i_j}|^2$ converges weakly to $|\nabla u|^2+\nu$. 

\begin{thm}\label{Lin}
If $(u_i)$ is a sequence of stationary maps in $\Rn$ converging weakly in $W^{1,2}$ to $u$ and $|\nabla u_{i}|^2$ converges weakly to $|\nabla u|^2+\nu$, then $u$ is stationary and the defect measure $\nu$ is  $(n-2)$-rectifiable.
\end{thm}

 This was  proved by Lin in \cite{Lin99a}, see also \cite{Lin99b}. The methods have some similar ingredients as those in Sections \ref{tanmeas} and \ref{Densities}, but Lin gives independent proofs relying on facts at hand.
 
De Lellis, Marchese, Spadaro and Valtorta \cite{DMSV18} and Hirsch, Stuvard and Valtorta \cite{HSV19} proved for Almgren's multivalued functions results analogous to Theorem \ref{enminthm}. Alper proved in \cite{Alp18} that the zero sets of harmonic maps from three-dimensional domains into a cone over the real projective plane are 1-rectifiable. In \cite{Alp20} he showed that the singular set of the free interface in an optimal
partition problem for the first Dirichlet eigenvalue in $\Rn$ is $(n-2)$-rectifiable.
 
\subsection{Mean curvature flow}\label{Flow}

A one parameter family $\{M_t, t\geq 0\}$ of compact surfaces moves by mean curvature if the normal velocity of $M_t$ equals the mean curvature vector at each point of $M_t$. I first discuss the case where the initial surface is a smooth hypersurface $M_0$ in $\Rn$. Then the flow is governed by a heat equation type partial differential equation: if $F:M_0\times[0,T]\to\Rn$ parametrizes part of the motion; $M_t=F(M_0\times\{t\})$, then $\partial_tF=H$. Here $H$ is  the mean curvature vector of $M_t$, which also is a Laplacian of $F$ in the metric of $M_t$. Then the area of $M_t$ is decreasing and singularities will appear. For instance, a sphere shrinks into a point and a cylinder into a line. See for example the survey \cite{CMP15} of Colding, Minicozzi and Pedersen and the book \cite{Man11} of Mantegazza for many interesting phenomena. Here I only briefly discuss the result of Colding and  Minicozzi \cite{CM16} on the rectifiability of singularities.

A point $(x,s), x\in M_s,s>0,$ is a singular point of the flow $(M_t)$ if $\{(y,t), y\in M_t,t>0\}$ is not a smooth manifold in any neighbourhood of $(x,s)$. Let $S\subset\Rn\times[0,\infty)$ be the set of singular points. The flow $(M_t)$ is called mean convex if $M_0$, and hence every $M_t$, has non-negative mean curvature. White \cite{Whi00} showed that for them the singular set has parabolic (and so also Euclidean) Hausdorff dimension at most $n-2$. Colding and  Minicozzi proved for such, and more general, flows that the singular set is $(n-2)$-rectifiable. In fact they proved much more:

\begin{thm}\label{CM16}
If the flow $(M_t)$ is mean convex, then the singular set $S$ can be covered with finitely many bi-Lipschitz images of  subsets of $\R^{n-2}$ together with a set of Hausdorff dimension $n-3$. 
\end{thm}

The paper \cite{CM16} contains many other facts about the structure of $S$. For example, the bi-Lipschitz images can be taken to be Lipschitz graphs with respect to the parabolic distance on $\Rn\times\R$, so $S$ is parabolic rectifiable, recall Section \ref{Parabolic}. Again, also the rectifiability of the stratification of $S$ is established.

The proof uses similar tools that we have seen above; monotonicity formula, tangent flows, and a parabolic Reifenberg theorem. The monotonicity formula is due to Huisken \cite{Hui90}; $t^{(1-n)/2}\int_{M_t} e^{-|y-x|^2/t}\,d\H^{n-1}y$ is non-decreasing. The tangent flows of $(M_t)$ are the weak limits of $\delta_i^{-1}M_{\delta_i^2t}, \delta_i\to 0.$ By an earlier result of Colding and  Minicozzi \cite{CM15} the tangent flows are unique for the mean convex flows, which is crucial in the proof of Theorem \ref{CM16}. In full generality the uniqueness of tangent flows is open. 
The parabolic Reifenberg theorem is now rather simple, because the approximating plane is assumed to be the same at all small scales. That this suffices is due to the strong information coming from the uniqueness of the tangent flows.

The tangent flows are not only unique but also cylindrical, that is, of the form $\R^k\times S^{n-k}$. The uniqueness means that $k$ and the direction of the $\R^k$ factor are unique. In fact, Colding and Minicozzi proved their results for all motions which have only cylindrical singularities. 
For $k<n$ the $k$th stratum $S_{k}$ consists of those points of $S$ where the Euclidean factor has dimension at most $k$. Then 
$S_0\subset S_1 \subset\dots\subset S_{n-2}=S$ and, by \cite{CM16}, each $S_{k}$ is $k$-rectifiable.

Much of the basics of the theory with smooth initial surfaces was established only in the 1980s and later. Surprisingly early, in the 1970s, Brakke \cite{Bra78} developed a very general theory with rectifiable $m$-varifolds $v_t, 0<m<n$. There are at least two good reasons to do this. Quite often the classical solutions develop singularities and the evolution in the classical sense stops. For varifolds singularities are allowed and the flow exists for all $t>0$. In many applications singularities are unavoidable. Brakke himself applied his flow to grain boundaries. For many further developments with a number of variants, connections and applications, see, for example, \cite{Ton19}, \cite{HL21} and the references given there.

In this general case the equation  for the flow is more complicated and in fact an inequality rather than equality is needed. This is necessary to have useful compactness. As mentioned in Section \ref{Varifolds} the first variation of a varifold leads to a concept of mean curvature. One also needs a weak formulation of the velocity. The first variations $\delta_{v_t}$ are assumed to be Radon measures and the curvatures are assumed to be in $L^2(\mu_{v_t})$. I skip the precise definitions.

The construction of Brakke's flow is a highly non-trivial matter. Brakke used a complicated approximation procedure. Ilmanen \cite{Ilm93} showed that  sequences of energy densities of the Allen-Cahn equation converge to rectifiable measures leading to a Brakke flow. In \cite{Ilm94} he used elliptic regularization. For other methods, also with prescribed boundary conditions, see \cite{QZ18}, \cite{Ton19}, \cite{Whi19}, \cite{ST21} and \cite{HL21} and their references. 

\subsection{Gromov-Hausdorff limits and related matters}\label{GH}
This section deals with one small but important part of the theory Riemannian manifolds with lower bound-ed Ricci curvature; structure of Gromov-Hausdorff limits of such manifolds. Since Gromov's fundamental work in the 1980s, see, for example, \cite{Gro81}, this topic has been intensively studied by many authors. Below I shall briefly discuss only matters directly related to rectifiability. For many other aspects, see, for example, the introductions of \cite{JN21} and \cite{CJN21} and the survey \cite{Nab20b}. There is a lot of similarity in results and some similarity in methods to those  for harmonic maps discussed in Section \ref{Harmmaps}. 

Recall the definitions related to Gromov-Hausdorff convergence from Section \ref{metrictan}. In addition, a \emph{tangent cone} of a 
metric space $X$ at a point $x\in X$ is any pointed limit $(X_{\infty},d_{\infty},x_{\infty})$ of $(X,r_i^{-i}d,x)$ where $r_i\to 0$.  

Let $(M_i)$ be a sequence of $n$-dimensional Riemannian manifolds and $m_i\in M_i$. We shall always assume that the Ricci curvatures and the volumes of the unit balls are bounded below: 
\begin{equation}\label{GHass}
\Ric_{M_i} \geq -(n-1)\ \text{ and}\ \H^n(B(m_i,1))> c > 0\ \text{for all}\ i.
\end{equation}
 Let $(X,d,x)$ be a pointed limit of $(M_i,m_i)$. Then, under the above conditions on $M_i$, tangent cones $Y$ exist at every point of $X$ by a compactness theorem of Gromov and they are metric cones, $Y=C(W)$, by a result of Cheeger and Colding \cite{CC97}. By definition, the metric cone $C(W)$ over a metric space $W$ is the completion of $(0,\infty)\times W$ with a particular metric, see \cite{CCT02}. They are also metric measure spaces where the measure is  a limit of volume measures. A point $x\in X$ is \emph{regular} if every tangent cone at $x$ is isometric to $\Rn$. Otherwise $x$ is \emph{singular}. We denote the set of singular points by $\sing(X)$. 
 
Colding and Naber \cite{CN13} showed that the tangent cones need not be unique under the above assumptions.

Cheeger and Colding proved in \cite{CC97} that $\dim \sing(X)\leq n-2$. Assuming the two-sided bound $|\Ric_{M_i}| \leq n-1$, Cheeger and Naber proved in \cite{CN15} that $\dim \sing(X)\leq n-4$. But we have more:

\begin{thm}\label{GHthm}
Let $n\geq 4$. Then $\sing(X)$ is $(n-4)$-rectifiable if $|\Ric_{M_i}| \leq n-1$.
\end{thm}

This was proved by Jiang and Naber  in \cite{JN21}. The proof is quite involved with many different kinds of techniques. Again monotonicity formulas play a decisive role, as they do for many other results mentioned below. Neck decompositions, recall Section \ref{sqgendim}, are central. Now instead of \eqref{NVeq} the $L^2$ integral of curvature is dominated by monotonic entities.

Cheeger and Colding introduced in \cite{CC97} the stratification of the singular set.
Let $S^k\subset\sing(X)$ be the set of points $x\in X$ for which no tangent cone at $x$ is of the form $C(W)\times\R^{k+1}$. Then 
$$S^0\subset S^1 \subset\dots \subset S^{n-1} =\sing(X),$$
and, by \cite{CC97}, $\dim S^k\leq k$,\ $S^{n-2}=\sing(X)$, and, by \cite{CN15}, $S^{n-4}=\sing(X)$ if the Ricci curvatures are also bounded above; $|\Ric_{M_i}| \leq (n-1)$.  Cheeger, Jiang and Naber proved the rectifiability in \cite{CJN21}:

\begin{thm}\label{GHthm1}
$S^k$ is $k$-rectifiable for all $k=0,1,\dots,n-2$. In particular, $\sing(X)$ is $(n-2)$-rectifiable.
\end{thm}

As for harmonic maps, quantitative stratification is considered in  \cite{CJN21} and quantitative volume estimates are obtained for them. The proof again involves monotonicity and neck decompositions.

With only the lower bound for the Ricci curvature rectifiability is essentially the best one can say; the singular sets can be Cantor sets, see \cite{LN20} and \cite{CJN21}.

Recall from the end of Section \ref{Harmmaps} Lin's results on defect measures related to sequences of stationary harmonic maps. Similar defect measures occur for sequences of connections and their Yang-Mills energies on $n$-dimensional manifolds. Tian proved their $(n-4)$-rectifiability in \cite{Tia00}. See also \cite{NV19} for further work involving neck type decompositions.


There are also many rectifiability results on metric measure spaces. I don't go into any details here, I just mention some of them.  Li and Naber proved in \cite{LN20} results analogous to the above for Alexandrov spaces with curvature bounded below. A special case of the results of Mondino and Naber \cite{MN19} shows  that metric measure spaces $(X,d,\H^n)$ with lower bounds for the Ricci curvature are $n$-rectifiable, see also \cite{BPS21} for another proof. 
Bru\'e, Naber and Semola showed in \cite{BNS20} that their boundaries are $(n-1)$-rectifiable. The boundary of $X$ is the closure of $S^{n-1}\setminus S^{n-2}$, which in this setting need not be empty. Here $S^{k}$ is a stratum as above. Bru\'e, Pasqualetto and Semola \cite{BPS19} developed De Giorgi's theory for sets of finite perimeter, including their rectifiability, in these and more general spaces. Lee, Naber  and Neumayer introduced in \cite{LNN20} rectifiable Riemann spaces to deal with some problems in Gromov-Hausdorff convergence. They are topological measure spaces, not necessarily metric. David \cite{Dav15} showed that the tangents of AD-regular Lipschitz differentiability spaces, recall Section \ref{Cheeger}, are uniformly rectifiable provided the space has charts of maximal dimension. Otherwise they are purely unrectifiable.

\subsection{Measure solutions of PDEs}\label{PDEsing}
Consider the system of constant coefficient linear partial differential equations on $\Rn$,
$$\mathcal A u  = \sum_{|\a|\leq k}A_{\a}\partial^{\a}u=0,\ u:\Rn\to\R^l,$$
where $\partial^{\a}=\partial_{\a_1}\dots\partial_{\a_n}, |\a|=\max\{\a_j:j=1,\dots,n\}$ for the multi-index $\a=\{\a_1,\dots,\a_n\}$ and the $A_{\a}:\R^l\to\R^p$ are linear maps. An $\R^l$ valued Radon measure on $\Rn$, $\mu\in\mathcal M(\Rn,\R^l)$, is said to be $\mathcal A$-\emph{free} if $\mathcal A\mu =0$ in the distributional sense.

Examples of $\curl$-free measures are the derivatives $Du$ of  BV-maps $u:\Rn\to\R^l$, see below.

By the Radon-Nikodym theorem any $\mu\in\mathcal M(\Rn,\R^l)$ can be written as $\mu=\mu_a+\mu_s=\mu_a+D(\mu,|\mu|)|\mu|_s$, where $\mu_a$ is absolutely continuous with respect to $\mathcal L^n$, $D(\mu,|\mu|)$ is the Radon-Nikodym derivative of $\mu$ with respect to the total variation measure $|\mu|$ and $|\mu|_s$ is the singular part of $|\mu|$ in its Lebesgue decomposition.

The structure of the singular parts $\mu_s$ is a much studied question with many applications, see \cite{DR16}, \cite{ADHR19} and the ICM survey \cite{DR18}. The wave cone
$$\Lambda_{\mathcal A}= \bigcup_{\xi\in\Rn\setminus\{0\}}\ker A^k(\xi)\subset\R^l\ \text{where}\ A^k(\xi)=\sum_{|\a|= k}\xi^{\a}A_{\a},\ \xi^{\a}=\xi_1^{\a_1}\dots \xi_n^{\a_n},$$
is central in this investigation, as well as in many other topics, see,  \cite{DR16} and \cite{DR18} and the references given there. I shall look at some examples below.

If $\mathcal A$ is homogeneous, $\mathcal A  = \sum_{|\a|= k}A_{\a}\partial^{\a},$ then $\lambda\in\R^l$  belongs to $\Lambda_{\mathcal A}$ if and only if there exists $\xi\in \R^n\setminus\{0\}$  such that $x\mapsto\lambda h(x\cdot \xi)$ is $\mathcal A$-free for all smooth functions $h: \R\to \R$. That is, in
the words of De Philippis and Rindler in \cite{DR18}: 'Roughly speaking, $\Lambda_{\mathcal A}$ contains all the amplitudes along which the system is not elliptic' and '“one-dimensional” oscillations and concentrations are possible only if the amplitude (direction) belongs to the wave cone'. 

The following theorem of De Philippis and Rindler gives a new proof of Alberti's rank one theorem \ref{Alberti} and extends it to BD-maps. It has many other consequences, too, see \cite{DR16}.
\begin{thm}\label{Deph}
If $\mu\in\mathcal M(\Rn,\R^l)$ is $\mathcal A$-free, then $D(\mu,|\mu|)(x)\in \Lambda_{\mathcal A}$ for $|\mu|_s$ almost all $x\in\Rn$.
\end{thm}

The proof uses tangent measures and pseudo-differential calculus. Partially based on similar ideas Arroyo-Rabasa, De Philippis, Hirsch and Rindler proved in \cite{ADHR19} several interesting rectifiability results. They are formulated in terms of $m$-dimensional wave cones
$$\Lambda^m_{\mathcal A}= \bigcap_{V\in G(n,m)} \bigcup_{\xi\in V\setminus\{0\}}\ker A^k(\xi).$$
 Clearly, these sets increase when $m$ increases.

\emph{Examples}
(1) The divergence operator on $\Rn$ is $\mathcal A=\sum_{i=1}^nA_i\partial_i$ where $A_ix=e_i\cdot x$, with $e_i$ the standard basis vectors. Then 
$A^1(\xi)x=\xi\cdot x$, so $\Lambda_{\mathcal A}=\Lambda^m_{\mathcal A}=\Rn$ for $2\leq m\leq n$ and $\Lambda^1_{\mathcal A}=\{0\}$.

(2) The curl is defined for $m\times n$-matrix valued measures $\mu$ on $\Rn$ by
$$\curl\mu=(\partial_i\mu^k_j-\partial_j\mu^k_i), i,j=1,\dots,n,k=1,\dots,m.$$
This can be written in the form $\curl\mu=\sum_{i=1}^nA_i\partial_i\mu$ for some $A_i$. Further it can be checked that for every $\xi\in\Rn\setminus\{0\}$ the kernel of $\curl(\xi)$ consists of the matrices $a\otimes \xi,\ a\in\R^m$, see \cite{FM99}, Remark 3.3(iii), for these facts.  It follows that $\Lambda^{n-1}_{\curl}=\{0\}$.

(3) Also the second order operator $\curl \curl$ on $n\times n$-matrix valued measures is of the above type with $\Lambda^{n-1}_{\curl\curl}=\{0\}$

Recall from Section \ref{projections} the integral geometric measure $\mathcal I^m_1$, whose null-sets are those which project to measure zero on almost all $m$-planes. By \cite{ADHR19}

\begin{thm}\label{Deph1}
Let $\mu\in\mathcal M(\Rn,\R^l)$ be $\mathcal A$-free and let $E\subset\Rn$ be a Borel set with $\mathcal I^m_1(E)=0$. Then 
$D(\mu,|\mu|)(x)\in \Lambda^m_{\mathcal A}$ for $|\mu|$ almost all $x\in E$. In particular, if 
$\Lambda^m_{\mathcal A}=\{0\}$, then $|\mu|(E)=0$, whence $|\mu|\ll\mathcal I^m_1\ll\mathcal H^m$. 
\end{thm}

Notice first that if $m=n$, and so $\mathcal I^m_1=\mathcal L^m$, this is the same as Theorem \ref{Deph}. In particular, if for $|\mu|$ almost all $x\in \R^m$,\  $A^k(\xi)(D(\mu,|\mu|)(x))\neq 0$ for $\xi\in \R^m\setminus\{0\}$, then $|\mu|$ is absolutely continuous.

The proof of Theorem \ref{Deph1} is by contradiction, the following sketch is rather imprecise and very incomplete. The counter-assumption leads to a subset $F$ of $E$ with $|\mu|(F)>0$, a point $x\in F$  and a plane $V\in G(n,m)$ such that $\H^m(P_V(F))=0$ and 
\begin{equation}\label{Arroyo}
A^k(\xi)(D(\mu,|\mu|)(x))\neq 0\ \text{for}\ \xi\in V\setminus\{0\}.\end{equation} 
Moreover, $x$ can be chosen so that for a sequence $\mu_j\in\mathcal M(B(0,1))$ of normalized  blow-ups at $x$ of $\mu\restrict F$ the total variations $|\mu_j|$ converge to a non-zero measure $\sigma\in\mathcal M(B(0,1))$. Setting $\nu_j=P_{V\#}(\mu_j)$ and $F_j=\spt|\nu_j|$, one has $\H^m(F_j)=0$. By delicate analysis based on \eqref{Arroyo} there is $\theta\in L^1(V\cap B(0,1))$ such that $\lim_{j\to\infty}||\nu_j|-\theta \H^m|(V\cap B(0,1))=0$. This leads to the contradiction
$$0<\sigma(B(0,1))\leq \liminf_{j\to\infty}\left(\int_{F_j}\theta\,d\H^m + ||\nu_j|-\theta \H^m|(F_j\cap B(0,1))\right)=0.$$

As a corollary to Theorem \ref{Deph1} we have

\begin{thm}\label{Deph2}
Let $\mu\in\mathcal M(\Rn,\R^l)$ be $\mathcal A$-free and suppose that $\Lambda^m_{\mathcal A}=\{0\}$. Then $\Theta^{\ast m}(|\mu|,x)<\infty$ for $|\mu|$ almost all $x\in \Rn$ and $|\mu|\restrict \{x\in\Rn:\Theta^{\ast m}(|\mu|,x)>0\}$ is $m$-rectifiable.
\end{thm}

By Theorem \ref{densmeas} and the last statement of Theorem \ref{Deph1} the upper density is finite $|\mu|$ almost everywhere. The set where $0<\Theta^{\ast m}(|\mu|,x)<\infty$ has $\sigma$-finite $\H^m$ measure and $|\mu|$ and $\H^m$ are mutually absolutely continuous on it by Theorem \ref{densmeas}. Hence the rectifiability follows from Theorem \ref{Deph1} and the Besicovitch-Federer projection theorem \ref{profed}. 

Here are some of the applications. 
Let $u=(u^1,\dots,u^l)\in BV(\Rn,\R^l)$. Then $\mu=Du=(\mu^k_i)=(\partial_iu^k),i=1,\dots,n,k=1,\dots,l,$ is $\curl$-free. 
By example 2 above $\Lambda^{n-1}_{\curl}=\{0\}$ and one can use Theorem \ref{Deph2} to get a new proof for the fact that $|Du|\restrict \{x:\Theta^{\ast m}(|Du|,x)>0\}$ is $(n-1)$-rectifiable, recall Theorem \ref{BVrect1}. In the same way, example 3 yields Theorem \ref{BDrect} for BD-maps.


As another application of Theorem \ref{Deph2} the authors of \cite{ADHR19} gave a new proof and extensions of Allard's rectifiability theorem \ref{Allard} and the related results of \cite{AS97} and \cite{DDG18}. Then $\mathcal A$ is a divergence operator on matrix valued measures. Lin's defect measure theorem \ref{Lin} also follows from Theorem \ref{Deph2}.

\subsection{A free boundary problem}\label{Freebdry1}

Recall from Section \ref{Freebdry} the work of David, Engelstein and Toro \cite{DET18}. Edelen and Engelstein \cite{EE19} studied the analogous one-phase   problem: Let $q$ be positive and H\"older continuous in $\Omega$. Minimize
$$J(u)=\int_{\Omega}(|\nabla u(x)|^2+q(x)^2\chi_{u>0}(x))\,dx$$
among $u\in W^{1,2}(\Omega)$ with given non-negative boundary values. This problem has an interesting history starting from Alt and Caffarelli in 1981 and there are similarities to codimension 1 minimal surface theory, see \cite{EE19}. In particular, let $k^{\ast}$ be the smallest $k$ such that the above problem admits a non-linear, one-homogeneous solution with $\Omega=\R^k, q=1$. Then $k^{\ast}$ is 5, 6 or 7, but otherwise the value is unknown. By a result of Weiss, we have for a minimizer $u, \dim \sing(u)\leq n-k^{\ast}$, meaning that $\sing(u)=\emptyset$, if $n<k^{\ast}$. The singular set now is the subset of $\partial\{u>0\}\cap\Omega$ where $u$ fails to be $C^{1,\a}$ for some $\a>0$. Edelen and Engelstein proved
\begin{thm}
If $u$ is a minimizer for $J$, then $\sing(u)$ is $(n-k^{\ast})$-rectifiable.
\end{thm}

They also considered stratification in the spirit we have seen before and they proved the rectifiability of each stratum. Their methods are influenced by those of Naber and Valtorta \cite{NV17}. They also apply to the two-phase problem and to almost minimizers.

For further work along these lines by De Philippis, Engelstein, Spolaor and Velichkov, see \cite{DESV21}.

\section{Miscellaneous topics related to rectifiability}
 
Here I only briefly present some other topics related to rectifiability.

\subsection{Curvature measures}\label{curvmeas}

Federer introduced in \cite{Fed59} sets with positive reach and curvature measures. This paper has had and is still having huge impact on matters related to  convexity and integral geometry. Sets with positive reach include both convex sets and $C^2$ submanifolds. 
A closed set $F\subset\Rn$ has positive reach if there is $r>0$ such that if $d(x,F)\leq r$ then there is a unique $\pi_F(x)\in F$ with $|\pi_F(x)-x|=d(x,F)$. 
 
Federer proved that there are signed Borel measures $\mu_i(F,\cdot)$, called \emph{curvature measures}, such that for every Borel set $B\subset\Rn$ and for $r>0$ as above
$$\mathcal L^n(\{x:d(x,F)\leq r, \pi_F(x)\in B\}) = \sum_{i=0}^n\a(n-i)r^{n-i}\mu_i(F,B).$$
For convex sets $F$ this is Minkowski's well-known Quermassintegrale formula. 

In addition Federer proved generalizations of the Gauss-Bonnet formula of differential geometry and the principal kinematic formula of integral geometry. His methods were partially based on his 1947 rectifiability theory. This paper provided one of the first applications of this theory.


In \cite{Zah86} Z\"ahle gave a representation of curvature measures in terms of rectifiable currents supported by the unit normal bundle of $F$, which is an $(n-1)$-rectifiable subset of $\R^{2n}$. Rataj and Z\"ahle continued this work in several papers, see, e.g., \cite{RZ05} and their monograph \cite{RZ19}, which also gives a wider presentation of the topic.


\subsection{Dynamical systems} For many dynamical systems the following dichotomy is typical: either the limit set is a fractal or something very special, for example, a piece of a plane, a sphere, or a real or complex analytic set. Such limit sets include self-similar and self-conformal sets, Julia sets of rational functions and limit sets of Kleinian groups. Fractal here could mean how Mandelbrot at one point defined fractal: Hausdorff dimension is strictly bigger than topological dimension. Mayer and Urba\'nski \cite{MU03} and Das, Simmons and Urba\'nski \cite{DSU17} used rectifiability to prove this type of results. In the paper \cite{DSU17} the setting is very general including the above cases as special cases, even in infinite-dimensional Hilbert spaces.

I describe the procedure vaguely (and with errors). Let $K$ be such a limit set. In all cases there is some self-similarity present: there are appropriate maps (similarities, conformal, etc.) which map small subsets of $K$ onto its large subsets. One should show that if the Hausdorff dimension $m=\dim K$ equals the topological dimension $\dim_TK$, then $K$ is something very special. First, the system offers some type of invariant measure, which can be related to $\H^m$ by scaling properties and then one has $\H^m(K)<\infty$. By \cite{Fed47b} the assumption  $\dim_TK=m$ implies that many projections of $K$ on $m$-planes have positive $\H^m$ measure. Hence by the Besicovitch-Federer projection theorem \ref{profed}, $K$ is not purely $m$-unrectifiable, so it has approximate tangent planes at many points $x$ by Theorem \ref{tanthm}. The appropriate maps send the planes to the very special sets we are after. Using these maps to blow-up small neighbourhoods of $x$ to large sets and taking a limit concludes the proof.

This approach does not work in infinite-dimensional spaces, in particular the projection theorem is false (recall Section \ref{Batepro}). The authors of \cite{DSU17} deal with this extending the family of rectifiable sets to what they call pseudorectifiable. A set $E$ with  $\H^m(E)<\infty$ is pseudorectifiable if there are $m$-planes $T_E(x)$ and a measure $\mu_E\sim\H^m\restrict E$ which, by definition, satisfy the area formula
$$\int\card\{x\in A:P_V(x)=y\}\,d\H^my = \int_A\det(P_V|T_E(x))\,d\mu_Ex$$
 for Borel sets $A\subset E$ and $m$-planes $V$. For rectifiable sets in finite dimensions the planes $T_E(x)$ are just the approximate tangent planes. Then any set $A$ with  $\H^m(A)<\infty$ splits into $m$-pseudorectifiable and purely $m$-unpseudorectifiable part. The latter is now defined by the property that it is a countable union of sets such that for each of them there exists a finite-dimensional linear subspace $V$ such that all projections into finite-dimensional linear subspaces containing $V$ are purely $m$-unrectifiable, in the classical sense. With these notions the above procedure can be followed, but with notable complications.

K\"aenm\"aki, Sahlsten and Shmerkin \cite{KSS15} used ergodic theoretic methods to investigate geometric properties of very general measures, involving also  rectifiability, see also the survey \cite{Kae18}. Each tangent measure of a measure $\mu$ at $x$ tells us something about $\mu$ around $x$ only at scales that generate it, which can be very sparse. In order to have a more complete picture we should look at all tangent measures at $x$. In-between there is a way to look at average behaviour via scenery flows and tangent 
distributions. These have been studied by many people, mainly in connection of fractal properties such as various dimensions. 

A tangent distribution is a measure on a space of measures. Let $\mu\in\mathcal M(B^n(0,1))$ and define for $x\in\spt\mu$ and $t\geq 0$ the probability measure $\mu_{x,t}$ on $B^n(0,1)$ by 
\begin{equation*}
\mu_{x,t}(A)=\frac{\mu(e^{-t}A+x)}{\mu(B(x,e^{-t}))},\ A\subset B^n(0,1).\end{equation*}
Then $(\mu_{x,t})_{t\geq 0}$ is called the scenery flow of $\mu$ at $x$. Letting $\delta_a$ denote the Dirac measure at $a$ set for $T>0$,
$$\langle\mu\rangle_{x,T} = \frac{1}{T}\int_0^T\delta_{\mu_{x,t}}\,dt.$$
Then any weak limit $P$ of a sequence $\langle\mu\rangle_{x,T_i}, T_i\to\infty,$ is called a tangent distribution of $\mu$ at $x, P\in\mathcal T\mathcal D(\mu,x)$. They are probability measures on $\mathcal M(B^n(0,1))$ and they enjoy very strong translation and scaling invariance properties by a result of Hochman. The support of a tangential distribution at $x$ is contained in the set of the  restrictions to the unit ball of the tangent measures at $x$.

The authors of \cite{KSS15} studied in particular conical density and porosity properties. By Theorem \ref{unrectsect}, and its higher dimensional analogues, for a purely $m$-unrectifiable set there is much measure in small cones around $(n-m)$-planes. The same is true for sets of Hausdorff dimension bigger than $m$. In \cite{KSS15} similar results are proven for general measures. In particular, the authors introduced a concept of average unrectifiability. A special case of their results states that if $\mu\in\mathcal M(B^n(0,1)), 0\leq p<1,$ and for every $P\in\mathcal T\mathcal D(\mu,x)$,
$$P(\{\nu\in\mathcal M(B^n(0,1)): \spt\nu\ \text{is not}\ m-\text{rectifiable}\})>p,$$
then for every $0<s<1$ there exists $0<\e<1$ such that
$$\liminf_{T\to\infty}\frac{1}{T}\mathcal L^1(\{t\in[0,T]: \inf_{V\in G(n,n-m)}\frac{\mu(X(x,V,s)\cap B(x,,e^{-t}))}{\mu(B(x,e^{-t}))}>\e\})>p$$
for $\mu$ almost all $x\in\Rn$. A converse holds if $\mu$ has positive lower and finite upper density almost everywhere.

Hovila, E. and M. J\"arvenp\"a\"a and Ledrappier \cite{HJJL14} proved that on a class of Riemann surfaces the union of complete geodesics has Hausdorff dimension 2 and $\H^2$ measure zero. To prove the second statement they used their generalization of the Besicovitch-Federer projection theorem for transversal families, recall Section \ref{projections}. The useful family was produced by investigating the geodesic flow on the tangent bundle.

Fuhrmann and Wang \cite{FW17} proved  that ergodic measures of certain dynamical systems on the 2-torus are 1-rectifiable.

\subsection{Higher order rectifiability}\label{Ckrect}

Anzellotti and Serapioni \cite{AS94} introduced higher order rectifiable sets. Let us say that $E\subset\Rn$ is $(m,k,\a)$-rectifiable if there are $m$-dimensional $C^{k,\a}$ submanifolds $M_i$ of $\Rn$ such that $\H^m(E\setminus \bigcup_iM_i)=0$. Here $k$ and $m$ are positive integers and $0\leq\a\leq 1$, with $C^{k,0}$ meaning $C^{k}$. Of course, $(m,1)$-rectifiable is then the same as $m$-rectifiable. Then, among $m$-rectifiable sets Anzellotti and Serapioni characterized $(m,2)$- and $(m,1,\a)$-rectifiable sets in terms of non-homogeneous blow-ups. That is, the blow-up maps are rotations of $(x,y)\mapsto (r^{-1}x,r^{-1-\a}y), x\in\R^m,y\in\R^{n-m}$.

As observed in \cite{AS94} the $(m,k,1)$-rectifiable sets are the same as the $(m,k+1)$-rectifiable sets due to the Lusin type theorem \cite[Theorem 3.1.15]{Fed69}. The question whether the upper limit in that theorem can be replaced by the approximate upper limit was solved in the negative by Kohn in \cite{Koh77}.

Santilli \cite{San19} characterized all $(m,k,\a)$-rectifiable sets with approximate differentiability. This means that he introduced a notion of approximate differentiability of higher order for subsets of $\Rn$ in analogy to the corresponding notion of functions (see \cite[Section 3.1]{Fed69} and \cite{Isa87}) and he proved that an $\H^m$ measurable set with finite $\H^m$ measure is $(m,k,\a)$-rectifiable if and only if it is almost everywhere approximately differentiable with parameters $m,k,\a$.  Santilli's definition of approximate differentiability of $E$ at $a$ involves conditions like 
$$\lim_{r\to 0}r^{-m}\H^m(\{x\in E\cap B(a,r): d(x,G)>sr^{k+\a}\})=0,$$
where $G$ is a graph of a polynomial of degree at most $k$ over an $m$-plane. I omit the precise definition, but I state the special case $k=1$ more explicitly: $E$ is $(m,1,\a)$-rectifiable if and only for $\H^m$ almost all $a\in E$ there exist $V\in G(n,m)$ and $s>0$ such that 
\begin{equation}\label{idu}
\lim_{r\to 0}r^{-m}\H^m(\{x\in E\cap B(a,r): |P_{V^{\perp}}(x-a)| > s|P_{V}(x-a)|^{1+\a}\})=0.\end{equation}    
Del Nin and Idu \cite{DI19} used this formulation to give a different proof and a slightly more general result in this case, such as Theorem \ref{unrectsect}. When $\a>0$, they also proved an analogue of Theorem \ref{tanmthm1}. That is, assuming positive lower density they gave a sufficient condition in terms of 'rotating' planes. The proof is much simpler than that of Theorem \ref{tanmthm1} since \eqref{idu} forces the approximating planes to converge at a geometric rate. In \cite{IM21} Idu and Maiale characterized in $\hn$ the $(m,1,\a)$ rectifiability, $n+2\leq m\leq 2n+1,$ in terms of approximate tangent paraboloids, following the scheme of \cite{DI19}.

Recently there have been many other results related to higher order rectifiability. Here are some:

Menne \cite{Men19} defined higher order differentiability of a set $A\subset\Rn$ by the approximability of $d(x,A)$ by polynomials over $m$-planes. He proved the higher order rectifiability of the sets where this happens. In \cite{Men21} he proved analogous results for distributions. 

It is easy to see that for any subset $E$ of $\Rn$ the set of points in $E$ that can be touched by a ball from the complement is $(n-1)$-rectifiable. Menne and Santilli \cite{MS19} showed that the set where a closed subset of $\Rn$ can be touched by a ball from at least $n - m$ linearly independent directions is $(m,2)$-rectifiable. The authors informed me that this also follows from Zajicek's results in \cite{Zaj79}. Hajlasz had related results in \cite{Haj20}.

Menne \cite{Men13} proved the $(m,2)$ rectifiability of integer multiplicity rectifiable varifolds whose first variation is a Radon measure and Santilli  \cite{San21} proved the same for general rectifiable varifolds with some extra conditions.

Bojarski, Hajlasz and Strzelecki \cite{BHS05} proved the $C^k$ rectifiability of level sets of $k$th order Sobolev mappings.

Kolasinski \cite{Kol17} and Ghinassi and Goering \cite{GG20} gave for higher order rectifiability sufficient conditions in terms of a Menger type curvature (recall \eqref{menger}), and Ghinassi \cite{Ghi20} and Del Nin and Idu \cite{DI19} in terms of square functions (recall Chapter \ref{Rectmeas}). 

Delladio has had many results on higher order rectifiability, for example of sets related to differentiability properties of functions, see \cite{Del13} and \cite{Del16b} and the references given there.

\subsection{Fractal rectifiability}
Let $0<s<m$ and $E\subset\Rn$ with $0<\H^s(E)<\infty$. Then we could consider $E$ as $(s,m)$-rectifiable if $\H^s\restrict E$ is $m$-rectifiable according to Definition \ref{m-rectmeas}. In \cite{MM88} it was studied how much of the theory of $m$ rectifiable sets could be extended to this setting. The paper contains several fairly easy positive results and many counter-examples.

Another possibility for $(s,m)$-rectifiability for non-integral $s$ would be to use $s/m$-H\"older maps from $\R^m$. Some fractal curves of positive and finite $s$-dimensional Hausdorff measure, for example, the von Koch snowflake, can be parametrized by $1/s$-H\"older maps from $\R$ and they would be rectifiable in this sense. On the other hand, many standard self-similar Cantor sets meet such H\"older curves in measure zero, see \cite{MM00}, and they would be purely unrectifiable. Investigation along these lines was made by Badger and Vellis in \cite{BV19}. In \cite{BNV19} Badger, Naples and Vellis gave an analyst's traveling salesman type condition  which implies that a set is contained in a $(1/s)$-H\"older curve.

\vspace{1cm}
\begin{footnotesize}
{\sc Department of Mathematics and Statistics,
P.O. Box 68,  FI-00014 University of Helsinki, Finland,}\\
\emph{E-mail address:} 
\verb"pertti.mattila@helsinki.fi" 

\end{footnotesize}

\end{document}